%% file: Imex_Practice_Builder_Generic.tex
\documentclass[11pt]{article}
\usepackage{amssymb}
\usepackage{geometry}
\usepackage{graphicx}
\usepackage{amsmath}
\usepackage{amssymb}
\usepackage{amsthm}
\usepackage{amsfonts}
\usepackage{enumerate}
\usepackage{bbm}
\usepackage{epsfig}
\usepackage{color}
\usepackage{setspace}
\usepackage{rotating}
\usepackage{subfig}
\usepackage{sidecap}
\usepackage{array}
\usepackage{calc}
\usepackage{tcolorbox}

\newlength{\mylength}

\numberwithin{equation}{section}

\newtheorem{theorem}{Theorem}[section]
\newtheorem{proposition}[theorem]{Proposition}

\newtheorem{definition}[theorem]{Definition}
\newtheorem{remark}{Remark}

\newtheorem{example}{Example}

\newtheorem{recipe}{Unconditional Stability Recipe}
\newtheorem{condition}{Condition}

\renewcommand{\vec}[1]{\mbox{\boldmath$#1$}}

\definecolor{orange}{rgb}{1,0.5,0}
\definecolor{dgreen}{rgb}{0,0.5,0}

\newcommand{\Srm}{\textrm{\S}}

\newcommand{\du}{\, \mathrm{d}}

\newcommand{\mat}[1]{\mbox{\boldmath$#1$}}

\newcommand{\param}{\delta}
\newcommand{\delt}{k}

\newcommand{\imex}{ImEx }
\newcommand{\timen}{n}
\newcommand{\sizeN}{N}
\newcommand{\diag}{\textrm{diag}}

\newcommand{\sfact}{\sigma}
\newcommand{\Wp}{W_p(\mat{A}, \mat{B})}

\def\XXint#1#2#3{{\setbox0=\hbox{$#1{#2#3}{\int}$}
     \vcenter{\hbox{$#2#3$}}\kern-.5\wd0}}

\begin{document}

\title{Unconditional Stability for Multistep ImEx Schemes: Practice} 

\author{
Benjamin Seibold\thanks{Department of Mathematics, Temple University, Philadelphia, PA 19122, seibold@temple.edu.}, 
David Shirokoff\thanks{Corresponding author. Department of Mathematical Sciences, NJIT, Newark, NJ 07102, \newline \indent \hspace{2mm} david.g.shirokoff@njit.edu.}, 
Dong Zhou\thanks{Department of Mathematics, California State University, Los Angeles, 5151 State University Dr. 
 \newline \indent \hspace{2mm} Los Angeles, CA 90032, dzhou11@calstatela.edu.}
}


\maketitle
\input{Sec_Abstract}

\medskip\noindent



\input{Main_Body}

\section*{Acknowledgments}
	The authors wish to acknowledge support by the National Science Foundation through grants DMS--1719640 (B. Seibold and D. Zhou) and DMS--1719693 (D. Shirokoff). D. Shirokoff was supported by a grant from the Simons Foundation ($\#359610$).

\appendix

\input{SecAppendix}

{\small
\bibliographystyle{siam}
\bibliography{../references_complete}
}

\end{document}

%% file: Sec_Abstract.tex
%
\begin{abstract}
This paper focuses on the question of how unconditional stability can be achieved via multistep \imex schemes, in practice problems where both the implicit and explicit terms
are allowed to be stiff. For a class of new \imex multistep schemes that involve a free parameter, strategies are presented on how to choose the \imex splitting and the time stepping parameter, so that unconditional stability is achieved under the smallest approximation errors. These strategies are based on recently developed stability concepts, which also provide novel insights into the limitations of existing semi-implicit backward differentiation formulas (SBDF).
For instance, the new strategies enable higher order time stepping that is not otherwise possible with SBDF.
 With specific applications in nonlinear diffusion problems and incompressible channel flows, it is demonstrated how the unconditional stability property can be leveraged to efficiently solve stiff nonlinear or nonlocal problems without the need to solve nonlinear or nonlocal problems implicitly.

\end{abstract}

%% file: Main_Body.tex

\section{Introduction}\label{Sec_Intro}

This paper builds on the theoretical work \cite{RosalesSeiboldShirokoffZhou2017} on the unconditional stability of linear multistep methods (LMMs). While \cite{RosalesSeiboldShirokoffZhou2017} introduced a new unconditional stability theory for implicit--explicit (ImEx) methods, and presented a novel class of \imex LMMs that involve a stability parameter, this paper develops strategies on how to select the time stepping parameter and the \imex splitting in an optimal fashion. The key focus is on problems for which an \imex splitting is warranted in which both the implicit and the explicit terms are stiff, for example because the stiff terms are difficult to treat implicitly.

Conventional \imex splittings often treat all stiff terms implicitly to ensure that one does not encounter a stiff time step restriction (one usually accepts a time step restriction from the non-stiff explicit part). However, as demonstrated in \cite{RosalesSeiboldShirokoffZhou2017}, this is not always required: one may treat stiff terms explicitly and nevertheless avoid a stiff time step restriction, provided the implicit term and the scheme are properly chosen. This paper provides strategies on how to make these choices (splitting and scheme) in practical problems. We do so through the use of the unconditional stability theory from \cite{RosalesSeiboldShirokoffZhou2017}, which is based on geometric diagrams that play a role analogous to the absolute stability diagram in conventional ordinary differential equation (ODE) stability theory. Specifically, we present strategies on how to achieve unconditional stability via (i) choosing the splitting for a given scheme; (ii) modifying a time-stepping scheme for a given splitting; and (iii) designing the splitting and the scheme in a coupled fashion. In addition, we employ the stability theory to provide new insights on the limitations of popular semi-implicit backward differentiation formulas (SBDF). In fact, we show that the new \imex LMMs generalize SBDF methods, in a way that they overcome some of their fundamental stability limitations.

\subsection{Problem setting}

We are concerned with the time-evolution of linear ODEs of the form
\begin{align}\label{Eq:ODE}
	\vec{u}_t = \mat{L} \vec{u} + \vec{f}, \quad \vec{u}(0) = \vec{u}_0.
\end{align}
Here $\vec{u}(t) \in \mathbb{R}^{\sizeN}$,
$\mat{L} \in \mathbb{R}^{\sizeN \times \sizeN}$, and
$\vec{f}(t) \in \mathbb{R}^{\sizeN}$ is an external forcing.
We assume that $\mat{L}$ gives rise to asymptotically stable
solutions --- i.e.\ solutions to the homogeneous ODE
$\vec{u}_t = \mat{L}\vec{u}$ decay in time (the
eigenvalues of $\mat{L}$ are in the strict left-half-plane). This assumption can be relaxed; however then additional
caveats are required (see \Srm\ref{Sec_NumExperiments} for when $\mat{L}$ has a zero
eigenvalue, or \Srm\ref{Sec_Conclusions} for when $\mat{L}$ has purely imaginary eigenvalues).

For the right hand side $\mat{L}$, an \emph{\imex splitting} $(\mat{A}, \mat{B})$ is conducted \cite{AscherRuuthWetton1995, Crouzeix1980, Varah1980, HundsdorferVerwer2003}, i.e.\ $\mat{L}$ is split into two parts, $\mat{L} = \mat{A} + \mat{B}$, where $\mat{A}$ is treated implicitly (Im) and $\mat{B}$ is treated explicitly (Ex). Clearly, the splitting $(\mat{A}, \mat{B})$ is non-unique, and in fact, any matrix $\mat{A}$ defines a splitting by choosing $\mat{B} := \mat{L} - \mat{A}$.  For this \imex splitting, we now require the time-stepping scheme to be unconditionally stable. This is a stringent, but very practical property (especially when $\mat{L}$ is stiff) as it allows one to choose a time step as large as accuracy requirements permit.

Note that the theory in this paper is developed for linear ODEs, as this assumption allows for a rigorous geometric stability theory involving unconditional stability diagrams. However, we then extend the results, in an ad-hoc but rather natural fashion, to nonlinear problems as well.

\subsection{Examples from partial differential equations}

A crucial source of stiff problems is the method-of-lines (MOL) semi-discretization of a partial differential equation (PDE). In that situation, rather than having one single right hand side $\mat{L}$, one faces a family $\mat{L}_h$ (with $h$ the mesh size) that approximates a spatial differential operator $\mathcal{L}$. A key property of the time-stepping strategies studied here is that for many PDE problems, the choice of \imex splitting and scheme can in fact be conduced on the level of differential operators, or equivalently, to hold for the family $\mat{L}_h$, uniformly in $h$ (see Sections~\ref{Sec_NumExperiments} and~\ref{Sec:NavierStokes}).

An important PDE situation in which unconditional stability is important is the MOL discretization of diffusion. A fully explicit treatment of diffusion gives rise to a stiff time step restriction $\delt \leq C h^2$. Hence, for problems in which diffusion represents the highest spatial derivative, a common approach is to include all of the discretization of $\frac{\partial^2}{\partial x^2}$ into the implicit part $\mat{A}_h$, and leave $\mat{B}_h$ as the remaining non-stiff terms. Such an approach will then avoid a stiff time step restriction. However, treating all stiff terms of $\mat{L}$ implicitly may in general be costly (see \Srm\ref{Sec_NumExperiments}); and in fact it is not always necessary. Having new approaches that allow one to treat (some of the) stiff terms explicitly, without incurring a stiff time step restriction, can be a significant practical benefit. In problems where $\mat{L}$ is stiff and costly to treat fully implicitly, this opens the door for designing a well-chosen \imex splitting where $\mat{A}$ contains only part of the stiff components of $\mat{L}$, and is much more efficient to treat implicitly.

\subsection{Background and relation to other works}

\imex unconditional stability has been studied in numerous theoretical and practical works.
On the theoretical side, general abstract sufficient conditions for unconditional stability and arbitrary multistep schemes are stated in \cite{Akrivis2013,AkrivisCrouzeixMakridakis1998,AkrivisCrouzeixMakridakis1999, AkrivisKarakatsani2003}. Although these conditions have the advantage of incorporating nonlinear terms
(i.e.\ $\mat{B}$ is allowed to be a nonlinear operator), they have the drawback that they require the implicit matrix $\mat{A}$ be larger (in the sense of an appropriate norm) than $\mat{B}$, and are overly restrictive for the problems we consider (e.g.\ they do not apply to Example~\ref{Ex:SimpleCase}).

Generally speaking, in the context of multistep methods, proofs for unconditional stability are commonplace for first and second order methods. Meanwhile, for higher order schemes, unconditional stability is usually only studied numerically, and in limited settings. This gap is likely due to the limitations that existing high-order methods encounter (see \Srm\ref{Sec_SBDF}). Important works in which unconditional stability is proved for first or second order methods, or numerically observed in higher order schemes, are the following papers (and references therein). Some of the first applications involving unconditional stability originated in the 1970s, with alternating direction implicit (ADI) methods \cite{DouglasDupont1970b}. Others include magneto-hydrodynamics \cite{HeisterOlshanskiiRebholz2017}; unconditional stability (also referred to as unconditionally energy stable, or as convex-concave splitting methods) for phase-field models \cite{Eyre1998, BertozziJuLu2011,ElseyWirth2013,ShengWangDuWangLiuChen2010,Smereka2003, GlasnerOrizaga2016, GuanLowengrubWangWise2014, YanChenWangWise2015, ElseyWirth2013, BadalassiCenicerosBanerjee2003}; applications to fluid-interface problems \cite{DucheminEggers2014}; incompressible Navier-Stokes equations \cite{JohnstonLiu2004, LiuLiuPego2007, HeisterOlshanskiiRebholz2016, KarniadakisIsraeliOrszag1991, KimMoin1985, LiuLiuPego2010}; Stokes-Darcy systems \cite{LaytonTrenchea2012}, compressible Navier-Stokes equations \cite{BrunoCubillos2016, BrunoCubillos2017}, and PDEs with the explicit treatment of non-local terms \cite{AnitescuLaytonPahlevani2004, Trenchea2016}. One disadvantage of low (i.e.\ first or second) order methods is that
they can also have large error constants for dissipative PDEs~\cite{ChristliebJonesPromislowWettonWilloughby2014} and dispersive PDEs~\cite{Ceniceros2002}, thus further reducing their applicability for the long-time numerical simulations. We differ from these previous works in several ways:
\begin{enumerate}[ 1.]
	\item We include higher order schemes as part of the study.
	\item Whereas many existing works use von-Neumann analysis or energy estimates that are tailored to a specific problem, we make use of recently introduced unconditional stability diagrams \cite{RosalesSeiboldShirokoffZhou2017}. The diagram approach simplifies the design of high-order unconditionally stable schemes and is applicable to a wider range of applications.
	\item We include variable \imex time stepping coefficients.  It may be surprising that stability considerations for \imex schemes do not require all stiff terms to be included in $\mat{A}$. In fact, \imex schemes can even go far beyond such a restriction: not only can $\mat{B}$ be stiff, it can (in some sense) even be \emph{larger} than $\mat{A}$, while still retaining unconditional stability. The underlying mechanism is that $\mat{A}$ is chosen in a way that stabilizes the numerical instabilities created by the explicit treatment of $\mat{B}$ with a suitable (simultaneous) choice of a splitting and time stepping scheme.
\end{enumerate}
It should also be stressed that there are numerous time stepping approaches (not strictly multistep methods) for specific application areas that possess good stability properties. Recently, high order unconditional stable methods for ADI applications have been obtained by combining second order multistep schemes with Richardson extrapolation \cite{BrunoJimenez2014, BrunoLyon2010I, BrunoLyon2010II}. For PDE systems that have a gradient flow structure, new conditions \cite{ShinLeeLee2017} allow for the design of third order, unconditionally energy-stable Runge-Kutta (RK) methods. Other techniques include: semi-implicit deferred correction methods \cite{Minion2003}; semi-implicit matrix exponential schemes where the linear terms are treated with an integrating factor \cite{KassamTrefethen2005, JuZhangZhuDu2014, MilewskiTabak1999}; and explicit RK schemes with very large stability regions for parabolic problems \cite{AbdulleMedovikov2001}.

\subsection{Outline of this paper}

This paper is organized as follows. After introducing the key notation and definitions (\Srm\ref{Sec2:Notation}), a self-contained review of the employed unconditional stability theory is provided that takes a different viewpoint than \cite{RosalesSeiboldShirokoffZhou2017} by placing a practical emphasis on the eigenvalues of $\mat{A}$, $\mat{B}$.
Section~\ref{SecRecipe} and onward (including \ref{Sec:AppendixProof}) contain new results.
Section~\ref{SecRecipe} provides recipes for designing (optimal) unconditionally stable
\imex schemes that minimize the numerical error.
Section~\ref{Sec_SBDF} characterizes the limitations of the well-known SBDF methods.
Section~\ref{Sec_NumExperiments} uses insight from \Srm\ref{Sec_SBDF} to overcome the limitations of SBDF and devise optimal high order (i.e.\ beyond 2nd order) unconditionally stable schemes for the variable-coefficient and non-linear diffusion problems. This section includes new formulas for \imex splittings and schemes (accompanied by rigorous proofs in \ref{Sec:AppendixProof}); as well as computational examples.
Section~\ref{Sec:NavierStokes} studies an application example that is motivated by incompressible Navier-Stokes flow in a channel and provides general insight into stability issues in computational fluid dynamics.
Section~\ref{Sec_Conclusions} provides an outlook and conclusions,
and \ref{Sec_ImexCoefficientsTables}
lists the specific \imex coefficients to be used in practice.



\section{Introduction to the ImEx schemes and unconditional stability property}\label{Sec2:Notation}
This section introduces the assumptions, notations, and \imex schemes
used throughout the paper.  As discussed above, we are interested in unconditional
stability for \imex splittings
$\mat{L} = \mat{A} + \mat{B}$ of equation \eqref{Eq:ODE}
where in general both the implicit
matrix $\mat{A}$, and the explicit matrix $\mat{B}$ are allowed
to be stiff.

We restrict to splittings in which $\mat{A}$ is Hermitian (symmetric in the real case) negative definite, i.e.~$\mat{A}$ has strictly negative eigenvalues:
\begin{align} \label{Assumpt_1}
	\mat{A}^{\dag} = \mat{A}, \quad \textrm{and}
	\quad \langle \vec{u}, \mat{A} \vec{u} \rangle < 0,
	\quad \textrm{for all } \vec{u} \neq 0, \; \vec{u} \in \mathbb{C}^{\sizeN}.
\end{align}
Here we have adopted the standard notation on vectors
$\vec{x}, \vec{y} \in \mathbb{C}^\sizeN$ (or $\mathbb{R}^\sizeN$):
\[
	\langle \vec{x}, \vec{y} \rangle = \sum_{j = 1}^\sizeN \overline{x}_j y_j,
	\quad
	\Vert \vec{x} \Vert^2 = \langle \vec{x}, \vec{x} \rangle,
	\quad
	\mat{A}^{\dag} = \overline{\mat{A}}^T,
	\quad
	\vec{x} = \begin{pmatrix}
		x_1, x_2, \cdots, x_{\sizeN}
	\end{pmatrix}^T.
\]
Note that $\mat{L}$ itself is not assumed symmetric/Hermitian or negative definite. Furthermore, assumption \eqref{Assumpt_1} on $\mat{A}$ is not overly restrictive, because for any given $\mat{L}$ one can choose $\mat{A}$ symmetric negative definite, and then set $\mat{B} = \mat{L} - \mat{A}$.
Note that spectral methods (for the spatial discretization of PDEs) may give rise to a complex matrix $\mat{A}$, which is why we do not restrict $\mat{A}$ to be real.
It is also worth noting that much of the theory we present still persists even
when $\mat{A}$ is not Hermitian and negative definite (see Section~\ref{Sec_Conclusions}).

Finally, we remark that the implicit treatment of a matrix $\mat{A}$ in multistep methods (or even Runge-Kutta methods), requires one to solve linear systems with coefficient matrices of the form $(\mat{I} - \gamma \delt \mat{A})$, where $\gamma > 0$ is a constant and $\delt > 0$ is the time step. For $\mat{A}$ symmetric negative definite, those system matrices are positive definite and thus favorable for fast solvers (chapter IV, lecture 38,
\cite{TrefethenBau1997}).


We will generally assume that the problem gives rise to a preferred/natural matrix structure $\mat{A}_0$ (symmetric, negative definite) that one wishes to treat implicitly; however, its overall magnitude is up to choice. In other words, the user fixes $\mat{A}_0$ and would accept any implicit matrix of the form $\mat{A} = \sfact \mat{A}_0$ (with the \emph{splitting parameter} $\sfact > 0$), provided that such an $\mat{A}$ yields unconditional stability.
This is in a spirit similar to \cite{DouglasDupont1970b}.
For example, in spatial discretizations of a variable coefficient
diffusion PDE where $\mat{L}\vec{u} \approx ( d(x) u_x )_x$, the user may prefer
an implicit treatment of the constant coefficient Laplacian
$\mat{A}_0 \vec{u} \approx u_{xx}$, however, would accept any constant
multiple as well, i.e.~$\mat{A}\vec{u} \approx \sfact u_{xx}$.
Writing $\mat{A} = \sfact \mat{A}_0$ where $\mat{A}_0$ is fixed, introduces
the scalar $\sfact$ as a key parameter. This paper shows how to choose $\sfact$
in a systematic fashion to obtain unconditionally stability.

We restrict our attention to \imex versions of linear multistep
methods (LMMs) \cite{AscherRuuthWetton1995, Crouzeix1980};
however it is worth noting that some of the concepts developed here may extend to other time stepping schemes as well, such as Runge-Kutta (multi-stage)
\imex schemes. The general form of an $r$-step LMM applied to the
ODE \eqref{Eq:ODE} with a splitting $(\mat{A}, \mat{B})$ is:
\begin{align} \label{Eq:fulltimestepping}
 \frac{1}{\delt}\sum_{j = 0}^{r} a_j \; \vec u_{\timen+j} =
    \sum_{j = 0}^{r} \Big( c_j \; \mat A \vec u_{\timen+j} +
    b_{j} \; \mat B \vec u_{\timen+j} + b_j  \vec f_{\timen+j} \Big)\/.
\end{align}
Here $\delt > 0$ is the time step,
the variable $b_r = 0$ (so that $\mat{B}$ is explicit in \eqref{Eq:fulltimestepping}),
$\vec{u}_{\timen} = \vec{u}(\timen \delt)$ is the numerical
solution $\vec{u}(t)$ (with a slight abuse of notation)
evaluated at the $\timen$-th time step, and
$\vec{f}_{\timen} = \vec{f}(\timen \delt)$.  We refer to the values
$(a_j\/,\,b_j\/,\,c_j)\/$, with $0 \leq j \leq r\/$
as the \imex (time stepping) coefficients. The LMMs of the form
\eqref{Eq:fulltimestepping} require $r$ initial conditions
$\vec{u}_0, \vec{u}_1, \ldots, \vec{u}_{r-1}$. The computation
of these initial conditions to sufficient accuracy is a separate problem
(chapter 5.9.3, \cite{LeVeque2007}), and is not considered here.
When discussing stability it will be useful
to define the polynomials $a(z)$, $b(z)$, $c(z)$,
using the \imex coefficients in \eqref{Eq:fulltimestepping}:
\begin{align} \label{Eq:PolynomialCoeff}
	a(z) = \sum_{j = 0}^r a_j z^j, \quad
	b(z) = \sum_{j = 0}^{r-1} b_j z^j, \quad
	c(z) = \sum_{j = 0}^{r} c_j z^j.
\end{align}
Given the \imex coefficients, one may write down
the polynomials $a(z)$, $b(z)$, $c(z)$, or alternatively,
given polynomials $a(z)$, $b(z)$, $c(z)$, one may read off
the different coefficients in front of $z^j$ to obtain the time stepping
coefficients $(a_j, b_j, c_j)$.

In this work we utilize a one-parameter family of \imex coefficients,
introduced in \cite{RosalesSeiboldShirokoffZhou2017}, that have desirable
unconditional stability properties.
The new \imex coefficients are characterized by a parameter
$0 < \param \leq 1$, i.e.~they are functions of a single
\emph{\imex parameter} $\param$, and are defined for
orders $r = 1$ through $r = 5$.
Formulas for the new coefficients $(a_j\/,\,b_j\/,\,c_j)\/$,
in terms of $\param$, may be found in Table~\ref{Table:ImexCoeff};
and substituting different values
of $0 < \param \leq 1$ into these formulas
yields different \imex schemes. For example, the new one-parameter
\imex schemes for first ($r=1$) and second order ($r = 2$)
take the form:
\begin{align} \label{Eq:1storderCoeff}
	&\textrm{1st order:} \quad	
	\frac{1}{\delt}\big( \param \; \vec{u}_{\timen+1} - \param \;
	\vec{u}_{\timen} \big)	=  \mat{A} \vec{u}_{\timen+1}
	+ (\param -1)\mat{A} \vec{u}_{\timen} +
	\param \; \mat{B} \vec{u}_{\timen}, \\ \label{Eq:2ndorderCoeff}
	&\textrm{2nd order:} \quad
	\frac{1}{\delt} \Big( \big(2\param - \frac{1}{2}\param^2\big) \vec{u}_{\timen+2}
	+\big(-4\param + 2\param^2\big) \vec{u}_{\timen+1} + \big(2\param - \frac{3}{2}
	\param^2\big) \vec{u}_{\timen} \Big) = \\ \nonumber
	&\mat{A} \vec{u}_{n+2} + 2\big(\param -1)\mat{A}\vec{u}_{\timen+1}
	+ \big(\param -1)^2 \mat{A} \vec{u}_{\timen}
	+ 2\param \; \mat{B} \vec{u}_{\timen+1}
	+ \big( (\param -1)^2 - 1\big) \mat{B} \vec{u}_{\timen},
\end{align}
	For brevity we have set $\mat{f} = 0$ in the formulas
	\eqref{Eq:1storderCoeff}--\eqref{Eq:2ndorderCoeff},
	however one may include it
	in the explicit term $\mat{B}\vec{u}$ (or even the implicit term) as in
	equation \eqref{Eq:fulltimestepping}.
	Although the formulas for the coefficients might appear unruly,
	they have simple polynomial expressions.
	\begin{remark} \label{Rmk:NewImExCoeff}
		(ImEx coefficients from Table~\ref{Table:ImexCoeff}
		written in polynomial form)
		For orders $1 \leq r \leq 5$, and $0 < \param \leq 1$, the ImEx
		coefficients $(a_j, b_j, c_j)$, for $0 \leq j \leq r$ from
		Table~\ref{Table:ImexCoeff} correspond to the following polynomials:
		\begin{align} \label{NewImex_a}
			 a(z) &= \sum_{j = 1}^r \frac{f^{(j)}(1)}{j!} (z- 1)^j,
			 \quad \textrm{where} \quad
			 f(z) = (\ln z) (z - 1 + \param)^r,	\\ \label{NewImex_b}
			 b(z) &= (z - 1 + \param)^r - (z- 1)^r, \quad \quad
			 c(z) = (z - 1 + \param)^r. 		
		\end{align}
		The relationships between the polynomials, i.e.
		$b(z) = c(z) - (z-1)^r$, and $a(z)$ as the $r$-th order
		Taylor polynomial of $\ln(z) c(z)$ ensure that the \imex
		coefficients satisfy the order conditions required to define
		an $r$-th order scheme.
	\end{remark}
	Note that in the Remark~\ref{Rmk:NewImExCoeff},
	the polynomial $c(z)$ has roots that approach $1$ as
	$\param \rightarrow 0$.  This is not an accident, and it is this property
	that will eventually lead to good unconditional stability properties for
	the new schemes.
		
	Equations \eqref{Eq:1storderCoeff}--\eqref{Eq:2ndorderCoeff},
    as well as the 3rd, 4th, 5th order schemes in Table~\ref{Table:ImexCoeff},
	define families of time-stepping schemes.
	When the value $\param = 1$
	is substituted into the coefficient formulas in
	equations \eqref{Eq:1storderCoeff}--\eqref{Eq:2ndorderCoeff}, one
	obtains the well-known backward differentiation
	formulas for the coefficients of $\mat{A}$, also
	referred to as semi-implicit backward differentiation
	formulas (SBDFr, where $r$ denotes the order of the scheme):
	\begin{align}\nonumber
		&\mathrm{SBDF1}\; (\param = 1): \quad	
	\frac{1}{\delt}\big( \vec{u}_{\timen+1} - \vec{u}_{\timen} \big)
	=  \mat{A} \vec{u}_{\timen+1} + \mat{B} \vec{u}_{\timen}, \\  \nonumber
		&\mathrm{SBDF2}\; (\param = 1): \quad
		\frac{1}{\delt} \Big( \frac{3}{2} \vec{u}_{\timen+2}
	-2\vec{u}_{\timen+1} + \frac{1}{2} \vec{u}_{\timen} \Big) =
	\mat{A} \vec{u}_{n+2} + 2 \mat{B} \vec{u}_{\timen+1} - \mat{B} \vec{u}_{\timen}.	
	\end{align}
	Choosing values $\param \neq 1$ yields different (new) schemes.
	We have only displayed orders $r = 1, 2$ in the above expressions,
	however coefficients are also given for orders
	$r = 3, 4, 5$ in Table~\ref{Table:ImexCoeff}.
	Lastly we note that the new \imex schemes are zero-stable for any value
    $0<\param \le 1$,
    and the coefficients satisfy the order
	conditions \cite{RosalesSeiboldShirokoffZhou2017}
	to guarantee that they
	define an $r$-th order scheme (i.e.~solving \eqref{Eq:fulltimestepping}
	using the coefficients approximates the solution to \eqref{Eq:ODE}
	with an error that scales like
	$\mathcal{O}(\delt^r)$ as $\delt \rightarrow 0$).

	Each fixed set of \imex coefficients, such as SBDF ($\param = 1$),
	or \imex versions of Crank-Nicolson, or even schemes not considered in this
	paper, provide unconditional stability for
	only a certain set of matrix splittings $(\mat{A}, \mat{B})$ ---
	and these may not include a practitioner's desired splitting
	for a given problem. Introducing
	the one-parameter family of \imex schemes (parameterized by $\param$)
	provides the flexibility needed to attain unconditional stability for new
	classes of matrices $(\mat{A}, \mat{B})$
	beyond the capabilities of what is possible using
	a fixed set of coefficients.
    This point becomes particularly apparent in \Srm\ref{Sec_SBDF}, in
    the discussion of the limitations of SBDF methods.
	This gain in unconditional stability offered by the parameter $\param$
    may come with a trade-off of increasing the numerical approximation
    error constants. Thus, an important discussion (see \Srm\ref{SecRecipe})
    is how to choose an \imex scheme (i.e.~how to choose $\param$)
	for a given problem splitting (i.e.~$(\mat{A},\mat{B})$)
	to balance the trade off of gaining unconditional stability while minimizing
	the numerical error. Or, even better, how to choose the splitting and scheme
    in a coupled fashion.


Our goal is to avoid unnecessarily small time step restrictions in the
numerical scheme \eqref{Eq:fulltimestepping}.  To do this we examine
when \eqref{Eq:fulltimestepping} is \emph{unconditionally stable} --- i.e.~the
numerical scheme \eqref{Eq:fulltimestepping} with $\vec{f} = 0$
remains stable regardless of how large one chooses the time step $\delt > 0$.
Formally, we adopt the following definition:
%
\begin{definition} (Unconditional stability)
 A scheme \eqref{Eq:fulltimestepping} is unconditionally stable if: when
 $\vec f = 0\/$, there exists a constant $C$ such that
 \[
      \|\vec{u}_{\timen}\| \leq C\max_{0 \leq j \leq r-1} \|\vec{u}_j\|\/,
      \quad \textrm{for all}\;\;\timen \geq r, \; \delt > 0
      \quad \textrm{and}\;\; \vec{u}_j \in \mathbb{R}^{\sizeN}\/,
      \;\;\textrm{where}\;\;  0 \leq j \leq r-1\/.
 \]
 Note that $C\/$ may depend on the matrices $\mat{A}$, $\mat{B}\/$, and the
 coefficients $(a_j, b_j, c_j)$, but is independent of the time step $\delt\/$,
 the time index $\timen\/$, and the initial vectors
 $\vec{u}_j$, $0 \leq j \leq r-1$.
\end{definition}
It is important to note that unconditional stability of an \imex LMM
like \eqref{Eq:fulltimestepping} can be difficult
to determine in practice, as this question depends
simultaneously on the choice of coefficients
$(a_j\/,\,b_j\/,\,c_j)\/$ and the splitting $(\mat{A}\/,\,\mat{B})$.
The purpose of introducing a new stability theory
in \cite{RosalesSeiboldShirokoffZhou2017} was to
remedy this difficulty and formulate unconditional stability (or
failure thereof) in terms of two separate computable quantities:
one quantity that depends only on the coefficients
$(a_j\/,\,b_j\/,\,c_j)\/$, and one that depends only on the
splitting $(\mat{A}\/,\,\mat{B})$.
The theory then allows for a variety of possibilities:
\begin{enumerate}[(i)]
	\item Given a fixed splitting $(\mat{A}, \mat{B})$,
	design coefficients $(a_j\/,\,b_j\/,\,c_j)\/$
	(by choosing $0 < \param \leq 1$)
	that achieve unconditional stability ---
	see \Srm\ref{SecRecipe}, Recipe~\ref{Recipe:ChooseParameter}.
	\item Given a fixed set of coefficients $(a_j\/,\,b_j\/,\,c_j)\/$
	(such as SBDF when $\param = 1$), determine how to choose a splitting
	$(\sfact \mat{A}_0, \mat{B})$ (i.e. choose $\sfact > 0$) that
	guarantees unconditional stability ---
	see \Srm\ref{SecRecipe}, Recipe~\ref{Recipe:ChooseSplitting}.
	\item Offer the most flexibility by simultaneously choosing both
	the coefficients $(a_j\/,\,b_j\/,\,c_j)\/$ and the splitting
	$(\sfact \mat{A}_0, \mat{B})$ to achieve unconditional stability.
	This will involve the simultaneous choice of
	$(\sfact, \param)$ and is discussed in
	\Srm\ref{SecRecipe}, Recipe~\ref{Recipe:ChooseBothImexAndSplit}.
\end{enumerate}

\section{The unconditional stability theory}\label{Sec_Review}

In this section we review the unconditional stability
theory from \cite{RosalesSeiboldShirokoffZhou2017} --- which
imposes conditions on $(\mat{A}, \mat{B})$ and the time-stepping
coefficients $(a_j, b_j, c_j)$ that (when satisfied)
ensure the unconditional stability of \eqref{Eq:fulltimestepping}.
The stability theory will then provide a guide for choosing
the \imex coefficients $(a_j, b_j, c_j)$ and/or
splitting $(\mat{A}, \mat{B})$ that guarantee
unconditional stability for a given problem (i.e.\ $\mat{L}$).
The unconditional stability theory
is somewhat analogous to the classical absolute
stability theory (chapter 7, \cite{LeVeque2007}),
as it relies on a stability diagram --- and
we highlight the parallels with an example here:
\begin{example}\label{Ex:AST} (Absolute stability theory)
	Given an ODE of the form
	$\vec{u}_t = \mat{A} \vec{u}$, the absolute
	stability diagram $\mathcal{A}$ is defined as
	\begin{align}\label{Eq:MultistepAST}
		&a_r \vec{u}_{n+r} + \ldots + a_0 \vec{u}_n =
		\delt \;
		( c_r \mat{A}\vec{u}_{n+r} + \ldots + c_0 \mat{A}\vec{u}_n), \\ \nonumber
		&\mathcal{A} = \big\{ \mu \in \mathbb{C} : a(z) = \mu c(z),
		\textrm{ has stable solutions z} \big\}.
	\end{align}	
	The scheme \eqref{Eq:MultistepAST} is stable with time
	step $\delt$, if and only if every
	eigenvalue $\lambda$ of $\mat{A}$
	(i.e.\ $\mat{A}\vec{v} = \lambda \vec{v}$)
	satisfies $\delt \lambda \in \mathcal{A}$ (with the
	possible exception of repeated eigenvalues $\lambda$, and time steps
	$\delt$ that lie on the boundary $\delt \lambda \in \partial \mathcal{A}$).
\end{example}
	A key feature of the absolute stability theory is that it
	\emph{decouples} the stability criteria into
	(i) a property of the matrix $\mat{A}$ only
	(i.e.\ the eigenvalues), in relation to (ii) a property of the
	time stepping scheme only (i.e.\ $\mathcal{A}$).
	Decoupling the stability theory is extremely useful;
	for instance, it allows one to determine which
	matrices $\mat{A}$ can be solved using a given time stepping scheme.
	The unconditional stability theory in this section will parallel that
	of the absolute stability theory, and:
\begin{itemize}
	\item Introduce the unconditional stability diagram
	(defined solely by $(a_j\/,\,b_j\/,\,c_j)\/$); and provide
	formulas for the diagrams to the schemes corresponding to Table~\ref{Table:ImexCoeff}.	
	\item Provide computable quantities in terms of $(\mat{A}, \mat{B})$
	that are analogous to the eigenvalues of $\mat{A}$ in Example~\ref{Ex:AST}.
	Unconditional stability will then be framed in terms of the computable
	quantities lying inside the unconditional stability region.
\end{itemize}

\subsection{The unconditional stability diagram $\mathcal{D}$}

The absolute stability theory in Example~\ref{Ex:AST} was obtained
by replacing the matrix $\mat{A}$ with one of its eigenvalues
$\lambda$ --- resulting in a (simpler) stability analysis of a scalar ODE.
In a similar spirit, if $(\mat{A}, \mat{B})$ can be simultaneously diagonalized
(for instance when they are commuting and diagonalizable matrices), then
$(\mat{A}, \mat{B})$ may be replaced by their
eigenvalues --- resulting likewise in a scalar ODE. The unconditional stability
diagram can then be derived from this scalar ODE.
We stress that although the diagram is derived here assuming $(\mat{A}, \mat{B})$
are simultaneously diagonalizable, the diagram is also applicable to
general matrices $(\mat{A}, \mat{B})$ (i.e.\ that do not commute),
as outlined below.
Suppose $\vec{v}$ is a simultaneous eigenvector
to $\mat{A}$ and $\mat{B}$ and satisfies
\begin{align} \label{EigenValueDef}
	-\mat{A}\vec{v} = \lambda \vec{v}, \quad
	\mat{B}\vec{v} = \gamma \vec{v}, \quad
	-\mu \mat{A} \vec{v} =  \mat{B} \vec{v}, \quad
	\textrm{where } \mu = \frac{\gamma}{\lambda}.
\end{align}
Here $\lambda > 0$ (and real)
since $\mat{A}$ is symmetric/Hermitian and positive definite.
Substituting $\vec{u}(t) = v(t) \, \vec{v}$ into the
ODE \eqref{Eq:ODE} yields the scalar equation
\begin{align}\label{Mode_V}
	v_t = -\lambda v + \gamma v.
\end{align}
One can then examine stability for the \imex scheme
\eqref{Eq:fulltimestepping}, applied to
equation \eqref{Mode_V} (with the $\lambda$ term treated implicitly
and the $\gamma$ term explicitly), in the usual way: set
$v_n = z^n v_0$, to obtain a polynomial equation for
the growth factors $z$
\begin{align}\label{Eq:ModelEquation}
	\delt^{-1} a(z) = -\lambda c(z) + \gamma b(z).
\end{align}
Here $a(z), b(z), c(z)$ are the polynomials defined in \eqref{Eq:PolynomialCoeff}.
Note that the polynomial equation \eqref{Eq:ModelEquation} was used in
\cite{AscherRuuthWetton1995} for the purpose of determining CFL-type
time step stability restrictions for advection-diffusion problems; and also in
\cite{FrankHundsdorferVerwer1997} in the context of computing absolute stability-type
diagrams for \imex schemes (see also \cite{Koto2009} for a treatment of delay differential equations).
In both cases, the matrices $(\mat{A}, \mat{B})$
were assumed to be simultaneously diagonalizable, and neither study
was focused on unconditional stability. Equation \eqref{Eq:ModelEquation} is
also sometimes used as a (non-rigorous) model for stability in the case
when $(\mat{A}, \mat{B})$ are not
simultaneously diagonalizable. The study of unconditional stability
digresses from prior work by
re-parameterizing equation \eqref{Eq:ModelEquation} with the substitution
$y = -\delt \lambda$ and $\mu = \gamma \lambda^{-1}$:
\begin{align} \label{Eq:ReparamModelEq}
	a(z) = y \Big( c(z) - \mu b(z) \Big).
\end{align}
Note that $y$ takes on all values $y < 0$ as $\delt$
varies between $0$ and $+\infty$; and that $\mu \in \mathbb{C}$.
For a fixed mode, i.e.\ fixed $\lambda$ and $\gamma$,
unconditional stability demands that the growth factors $z$ solving
equation \eqref{Eq:ModelEquation}
are stable for all $\delt > 0$.
Viewed in the context of \eqref{Eq:ReparamModelEq}, this requirement
leads to the definition of the unconditional stability diagram
$\mathcal{D}$: the values $\mu \in \mathbb{C}$ for which
the growth factors $z$ to \eqref{Eq:ReparamModelEq}
are stable for all $y < 0$ (including $y \rightarrow -\infty$)
\begin{align}\nonumber
 \mathcal{D} &:= \Big\{ \mu \in \mathbb{C} :
 \textrm{Solutions } z \textrm{ to }
 \eqref{Eq:ReparamModelEq} \textrm{ are stable for all } y < 0
 \Big\}\/.
\end{align}
Here we say that $z$ is stable if $|z| < 1$; and for
technical convenience we exclude (non-repeated) values of $|z| = 1$.
Thus far, the definition for $\mathcal{D}$ is very general and may be
computed for any set of \imex LMM coefficients $(a_j, b_j, c_j)$.
It is also crucial to note that $\mathcal{D}$ is defined \emph{only} in
terms of the \imex scheme coefficients.

\begin{figure}[htb!]
	\includegraphics[width = .98\textwidth]{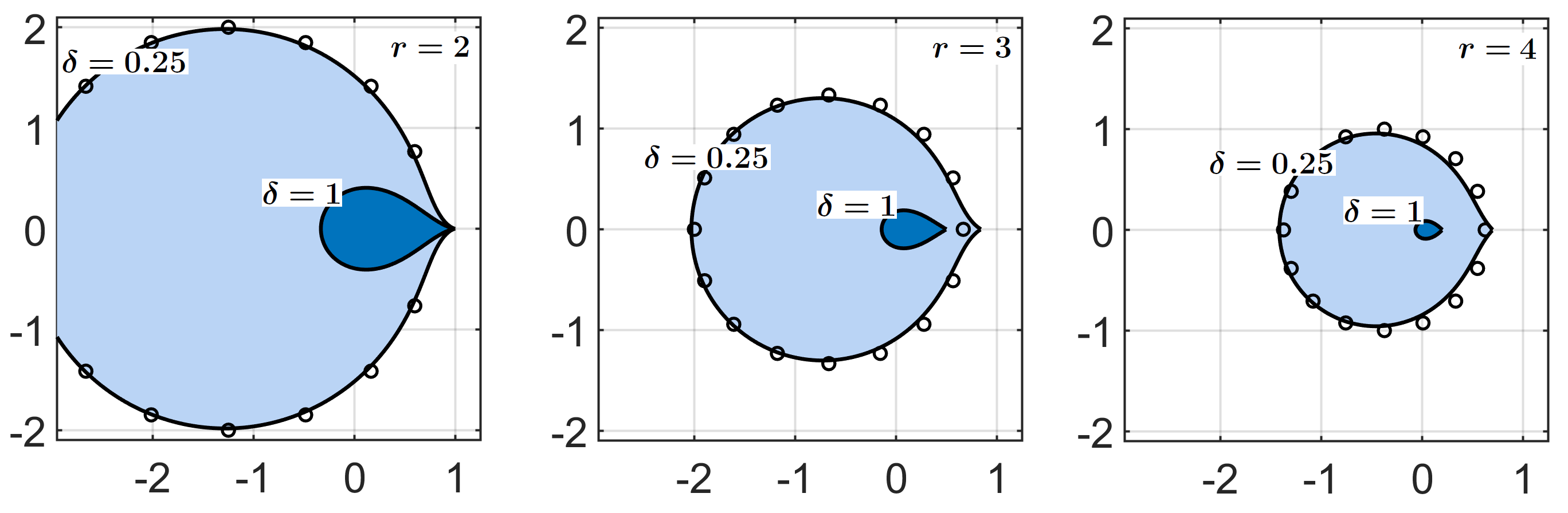}
	\caption{The sets $\mathcal{D}$ for orders (left to right) $r \in \{2,3,4\}$.
		The set $\mathcal{D}$ for $\param = 1$ (SBDF) (dark blue) is much smaller than
		$\mathcal{D}$ for $\param = 0.25$ (light blue). The
		asymptotic circle $C$ in formula (B3) for $\param = 0.25$ is shown in dots ($\circ$).
		The stability regions also decrease in size with increasing $r$. The orders $r = 1, 5$ (not plotted) exhibit a similar behavior.
	}
\label{Fig:NewImExStabilityRegions}
\end{figure}

It was proved (Thm.~8, Prop.~9 \cite{RosalesSeiboldShirokoffZhou2017}) that for 
the schemes in Table~\ref{Table:ImexCoeff}, the
value of $y \rightarrow -\infty$ (i.e.\ requiring stability for
large time steps, $\delt\rightarrow \infty$)
imposes the most severe restriction on the growth factors
in equation \eqref{Eq:ReparamModelEq}. This theoretical result has the
consequence that the set $\mathcal{D}$
is completely determined by setting $y \rightarrow -\infty$ in
\eqref{Eq:ReparamModelEq}, leading to the simplification
\begin{align}\label{Eq:ThmStabDiagram}
	\mathcal{D} = \Big\{ \mu \in \mathbb{C} :
	c(z) - \mu b(z) \textrm{ has stable roots} \Big\}
	\quad \textrm{(For schemes in Table~\ref{Table:ImexCoeff})}
\end{align}
We stress that \eqref{Eq:ThmStabDiagram} is not necessarily a general
property of \imex LMM --- but it holds for the schemes in Table~\ref{Table:ImexCoeff}.
Equation~\eqref{Eq:ThmStabDiagram} is useful as it allows
one to compute $\mathcal{D}$ (Thm.~8 \cite{RosalesSeiboldShirokoffZhou2017})
in terms of a \emph{boundary locus} formulation
(chapter 7.6, \cite{LeVeque2007}) with the polynomials $b(z)$ and $c(z)$
introduced in Remark~\ref{Rmk:NewImExCoeff}:
\begin{enumerate}
	\item[(B1)] The set $\mathcal{D}$ (for orders $1 \leq r \leq 5$)
	includes the origin (i.e.\ $0 \in \mathcal{D}$) and has the boundary
\begin{align}\label{ExactBoundary}
 &\partial \mathcal{D} =
   \Big\{ \frac{(z-1 + \param)^r}{(z-1 + \param)^r-(z-1)^r} : |z| = 1, \;
   \mathrm{arg}\; z_0 \leq \mathrm{arg} \; z \leq \;
   2\pi - \mathrm{arg}\; z_0 \Big\},\\  \nonumber
   &\textrm{with: } z_0 = 1, \textrm{ for } r = 1,
  \textrm{ and }
   z_0 = \frac{2-\param - 2(1-\param)\cos(\pi/r)
   e^{\imath \pi/r}}{2-\param - 2\cos(\pi/r) e^{\imath \pi/r}},
	\textrm{ for } 2 \leq r \leq 5.
\end{align}
	\item[(B2)] The right-most point $m_r$ and left-most point $m_l$ of
	$\partial \mathcal{D}$ are on the real axis with:
	\begin{align*}
		m_l = \Big( 1 - (1 - \param/2)^{-r} \Big)^{-1},
		\quad
		m_r =
		\left\{
		   \begin{array}{cl}
		     	1, & r = 1, \\
		     	\Big( 1 + \big( (1 - \param/2)\sec(\pi/r) \big)^{-r} \Big)^{-1}, & 2 \leq r \leq 5.
		   \end{array}
		\right.	
	\end{align*}
	\item[(B3)] In the asymptotic limit $\param \ll 1$, the set
	$\mathcal{D}$ approaches the circle $C$,
	where
	\begin{align*}
	  C &= \Big\{ z \in \mathbb{C} : \Big|z + \frac{1}{r\param} -
       \frac{r+1}{2r}\Big| \leq \frac{1}{r \param} \Big\}.	
	 \end{align*}
	Note that $C$ has a center at $\sim -\frac{1}{r \param}$ and radius
	$\sim \frac{1}{r \param}$; and hence becomes arbitrarily large
	as $\param \rightarrow 0$. 	Therefore, $\mathcal{D}$ becomes large
	as $\param \rightarrow 0$.
\end{enumerate}
Figure~\ref{Fig:NewImExStabilityRegions} plots the stability
diagrams $\mathcal{D}$ for different orders
and $\param$ values --- and also shows that
$\mathcal{D}$ asymptotically approaches (as $\param \rightarrow 0$) the
large circle $C$.
Having formulas for the shape and size of $\mathcal{D}$ as functions
of $\param$ will be important for designing
unconditionally stable schemes \eqref{Eq:fulltimestepping}, and for
characterizing the limitations of well-known schemes such as SBDF.
Lastly, we note that the \imex schemes parameterized by $\param$ bare
some similarity to the non-\imex schemes with large regions
of absolute stability originally examined in
\cite{JeltschNevanlinna1981, JeltschNevanlinna1982}.
However, we will eventually choose the
parameter value $\param$ to be as large as possible
(to minimize the error), while maintaining unconditional
stability. This is of a fundamentally different nature than the
non-\imex study carried out in
\cite{JeltschNevanlinna1981, JeltschNevanlinna1982}.

We now come to our first condition for unconditional
stability --- which is stated in terms of the
generalized eigenvalues
\begin{align}
	\Lambda( \mat{A}, \mat{B} ) := \big\{ \mu \in \mathbb{C} :
	-\mu \mat{A} \vec{v} = \mat{B}\vec{v}, \vec{v} \neq \vec{0} \big\}.
\end{align}
Note that a negative sign was added, for convenience, in the definition of
$\Lambda(\mat{A}, \mat{B})$ to make $(-\mat{A})$ positive
definite; and that $\Lambda(\mat{A}, \mat{B})$
is equivalent to the eigenvalues of $(-\mat{A})^{-1} \mat{B}$.
%
%
\begin{condition} \label{Cond:CommMatrices}
	(Unconditional stability when $(\mat{A}, \mat{B})$
	are simultaneously diagonalizable)
	Given time stepping coefficients $(a_j, b_j, c_j)$ with diagram $\mathcal{D}$, and
	simultaneously diagonalizable matrices $(\mat{A}, \mat{B})$ with
	generalized eigenvalues $\Lambda(\mat{A}, \mat{B})$,
	we have the following\dots
	\begin{itemize}
		\item[$\mathrm{(SC)}$] Sufficient conditions:
		The scheme \eqref{Eq:fulltimestepping} is unconditionally stable
	if every generalized eigenvalue $\mu \in \Lambda(\mat{A}, \mat{B})$ lies in
	$\mathcal{D}$, i.e.\ $\mu \in \mathcal{D}$.
		\item[$\mathrm{(NC)}$] Necessary conditions: If a generalized
	eigenvalue $\mu \in \Lambda(\mat{A}, \mat{B})$
	is not in $\mathcal{D}$, i.e.\ $\mu \notin \mathcal{D}$, then the
	scheme \eqref{Eq:fulltimestepping} is not unconditionally
    stable.\footnote{Strictly speaking, the
	precise theorem (Proposition~10, \cite{RosalesSeiboldShirokoffZhou2017})
	is that if $\mu \notin \mathcal{D}$, or $\mu \notin \varGamma$
	where $\varGamma = \{c(z)/b(z) : |z| = 1\}$ is the boundary locus of $\mathcal{D}$,
	then the scheme is not unconditionally stable.
    However, for practical purposes, the boundary locus
	can be ignored since it is a curve.}
	\end{itemize}
\end{condition}
%
%
In Condition~\ref{Cond:CommMatrices}, the (NC) and (SC) are essentially
identical and give a sharp characterization of unconditional
stability. Although Condition~\ref{Cond:CommMatrices} is useful when $(\mat{A}, \mat{B})$
are simultaneously diagonalizable, we also wish to consider
matrices $\mat{A}$ and $\mat{B}$ that do not commute.
The results in \cite{RosalesSeiboldShirokoffZhou2017} generalize the (SC) in
Condition~\ref{Cond:CommMatrices} to arbitrary matrices $(\mat{A}, \mat{B})$
($\mat{A}$ still symmetric positive definite)
by replacing the set $\Lambda(\mat{A}, \mat{B})$
with a (somewhat larger) set defined in terms of a
\emph{numerical range} (also known as the \emph{field of values}).
Specifically, let $p \in \mathbb{R}$ be any real number (different
values of $p$ will eventually be useful for different problem matrices
$\mat{L}$), and introduce the following sets:
\begin{align}\label{Def:Wp}
   \Wp := \Big\{ \langle \vec v, (-\mat A)^{p-1} \mat B \vec v\rangle :
          \langle \vec v, (-\mat A)^{p} \vec v\rangle = 1,
          \vec{v} \in \mathbb{C}^{\sizeN} \Big\}\/.
\end{align}
The set $\Wp$ can also be written, using
a change of variables $\vec{v} = (-\mat{A})^{\frac{p}{2}} \vec{x}$, as:
\begin{align}\label{Eq:NumRange}
	\Wp &=  W\Big( (-\mat{A})^{\frac{p}{2} - 1} \; \mat {B} \;
    (-\mat{A})^{-\frac{p}{2}} \Big), \quad \text{where} \\
    W(\mat X) &:= \big\{ \langle \vec x, \mat X \vec x \rangle :
            \| \vec x \| = 1,\; \vec x \in \mathbb{C}^{\sizeN} \big\}\/.
\end{align}
Here $W(\mat{X})$ is the definition of the \emph{numerical range}
of a matrix; and is a well-known set (chapter 1, \cite{HornJohnson1991})
that may be computed using a sequence of eigenvalue computations
\cite{Johnson1978}.  Note that $\Wp$ depends only on
the matrix splitting $(\mat{A}, \mat{B})$ and
is independent of the time stepping
coefficients. Condition~\ref{Cond:CommMatrices} may then be modified as follows.
\begin{condition} \label{Cond:NoCommMatrices}
	((Theorem~5, \cite{RosalesSeiboldShirokoffZhou2017})
	Unconditional stability for a general splitting $(\mat{A}, \mat{B})$)
	\begin{itemize}
		\item[$\mathrm{(SC)}$] Sufficient conditions:
		The scheme \eqref{Eq:fulltimestepping} is unconditionally stable
		if there is a value of $p \in \mathbb{R}$ for which
		the set $\Wp$ is contained
		in $\mathcal{D}$, i.e.\ $\Wp \subseteq \mathcal{D}$.
		\item[$\mathrm{(NC)}$] Necessary conditions: If a generalized
			eigenvalue $\mu \in \Lambda(\mat{A}, \mat{B})$
			is not in $\mathcal{D}$, i.e.\ $\mu \notin \mathcal{D}$, then the
			scheme \eqref{Eq:fulltimestepping} is not unconditionally stable.			
	\end{itemize}	
\end{condition}
Note that in Condition~\ref{Cond:NoCommMatrices} the (NC) are the same
as in Condition~\ref{Cond:CommMatrices}, however the (SC) are no long the same ---
due to the non-commuting matrices.
In Conditions~\ref{Cond:CommMatrices}--\ref{Cond:NoCommMatrices}
the (SC) provide a target criterion that will ensure unconditional stability; while the
(NC) will provide insight into when a scheme may fail to be unconditionally stable.

We provide a brief explanation here for why one should
replace $\Lambda(\mat{A}, \mat{B})$ with the sets $\Wp$ in Condition~\ref{Cond:NoCommMatrices}.
If one seeks
an eigenvector solution to \eqref{Eq:fulltimestepping}
of the form $\vec{u}_n = z^n \vec{v}$; and then
multiplies equation \eqref{Eq:fulltimestepping} from the left by
$(-\mat{A})^{p-1} \vec{v}$, then one obtains equation
\eqref{Eq:ReparamModelEq}, with the modification that the
value $\mu$ is no longer a generalized eigenvalue, but is
given by a general Rayleigh quotient:
\[
	 \mu = \frac{\langle \vec v, (-\mat A)^{p-1} \mat B \vec v\rangle }{\langle
      \vec v, (-\mat A)^{p} \vec v\rangle}\/ \subseteq \Wp.
\]
Hence, ensuring $\Wp \subseteq \mathcal{D}$ guarantees that the value of
$\mu$ in \eqref{Eq:ReparamModelEq} lies within the unconditional stability
region.

%
%

\begin{remark}\label{Rmk:NumRange} (Properties of the numerical range and $\Wp$)
	Since the sets $\Wp$ can be written in terms of a numerical range,
	they exhibit all the well-known properties of a numerical range.
	The numerical range $W(\mat{X})$ for a matrix $\mat{X}$ is
	convex (Hausdorf-Toeplitz theorem), bounded, and always contains the
	eigenvalues $\mu$ of $\mat{X}$, i.e.\ $\mu \in W(\mat{X})$.
	In the case when $\mat{X}$ is a normal matrix, $W(\mat{X})$ is the
	convex hull of the eigenvalues. Hence, the convex hull of
	$\Lambda(\mat{A}, \mat{B})$ is contained in $\Wp$
	(for all $p \in \mathbb{R}$).
\end{remark}

\begin{remark}
 Different values of $p$ may modify the size of $\Wp$ in the complex plane.
 Condition~\ref{Cond:NoCommMatrices} only requires one
 value of $p$ to satisfy $\Wp \subseteq \mathcal{D}$ (even if other
 values of $p$ violate $\Wp \subseteq \mathcal{D}$).
\end{remark}


\section{How to choose the \imex parameter $\param$ and splitting $(\mat{A}, \mat{B})$}\label{SecRecipe}

In this section we provide general recipes for choosing the \imex
parameter $\param$ and the matrix splitting $(\sfact \mat{A}_0, \mat{B})$
for a problem matrix $\mat{L}$.
The recipes are based on minimizing a proxy for the numerical error
while ensuring that the sufficient conditions (SC) are satisfied.

Solely based on the formulas for $\mathcal{D}$, one could think that one should use \imex coefficients with very large unconditional stability region $\mathcal{D}$, by taking $\param \ll 1$. After all, such a choice would increase the chance of unconditional stability by ensuring that $\Wp$ fits inside $\mathcal{D}$ thereby satisfying the (SC) in Conditions~\ref{Cond:CommMatrices}--\ref{Cond:NoCommMatrices}.

However, choosing $\param$ small without any regard for the
error is not a good strategy. Specifically, there is a
trade-off between schemes with good unconditional stability
properties (i.e.\ small $\param$ and large $\mathcal{D}$) and the
resulting numerical accuracy. Ideally, one would choose $\param$
so that the scheme's numerical approximation error is minimized,
while still guaranteeing unconditional stability. However,
because the true error is generally not accessible, we use $\param$
as a proxy for the approximation quality, which is justified by
the following remark.

\begin{remark}\label{Rmk:GTE_Err}
(Dependence of the global truncation error constant on $\param$)
The \emph{global truncation error} (GTE) at time $t_n = n\delt$ is defined by
$\max_{1\leq j \leq \sizeN} \left| \vec{u}_n - \vec{u}^*(n\delt)\right|_j$.
Because the \imex schemes in Remark~\ref{Rmk:NewImExCoeff} are formally of $r$-th order, for any fixed $0 < \param \leq 1$, the GTE scales (for $\delt$ small) like $C_r \delt^r$.
The error constant $C_r$ depends on $\mat{A}$, $\mat{B}$, $\vec{f}$, and the time stepping coefficients. Formulas for the behavior of the GTE error constants in a LMM may be computed in terms of the polynomials (see equation (2.3), p.~373, in \cite{HairerNorsettWanner1987}) $b(z)$ and $c(z)$. In particular, one may compute two separate error constants. One error constant is obtained when the \imex scheme is applied as a fully implicit scheme (i.e.\ $\mat{A} = \mat{L}$, $\mat{B} = 0$) as $C_r \propto 1/c(1) = \param^{-r}$. A second error constant may be computed when the \imex scheme is applied to a fully explicit splitting (i.e.\ $\mat{B} = \mat{L}$, $\mat{A} = 0$), where $C_r \propto 1/b(1) = \param^{-1}$. In general, for a fixed splitting $(\mat{A}, \mat{B})$, one then has a GTE that scales like
\begin{align}\label{Eq:AsymError}
	\text{GTE} \sim \mathcal{O}( \param^{-r} \delt^r ).
\end{align}
A more detailed description, along with numerical error
tests verifying the asymptotic formula \eqref{Eq:AsymError} may be found
in \cite{RosalesSeiboldShirokoffZhou2017}.
\end{remark}

Remark~\ref{Rmk:GTE_Err} indicates that for a fixed splitting $(\mat{A}, \mat{B})$, the GTE error is (asymptotically) minimized by taking a maximum value of $\param$. Moreover, as a secondary trend, if a family of \imex splittings $(\sfact \mat{A}_0, \mat{B})$ is considered, then it is generally observed that smaller values of $\sfact$ yield a smaller GTE. Hence, one should generally choose $\param$ as \emph{large} as possible and $\sfact$ \emph{small}, while still satisfying the (SC) constraint in Conditions~\ref{Cond:CommMatrices}--\ref{Cond:NoCommMatrices}.

We now provide recipes for three different scenarios that may arise in practice. Recipe~\ref{Recipe:ChooseParameter} specifies how to choose the \imex parameter $\param$ to achieve unconditional stability when a fixed matrix splitting $(\mat{A}, \mat{B})$ is specified (i.e.\ this a special case where $\sfact = 1$ and $\mat{A} = \mat{A}_0$).


\begin{tcolorbox}
\begin{recipe}\label{Recipe:ChooseParameter}
	How to choose the \imex parameter $\param$ for a fixed
	matrix splitting $(\mat{A}, \mat{B})$.
    \vspace{-.5\topsep}
	\begin{enumerate}[1.]\itemsep0em\parskip.3em\setcounter{enumi}{-1}
		\item Choose an order $1 \leq r \leq 5$; and retrieve the
		formulas for $\mathcal{D}$ in equation \eqref{ExactBoundary}.
		\item Compute/plot the generalized eigenvalues
   	 	$\Lambda(\mat{A}, \mat{B})$
   	 	and the sets $\mathcal{D}$ for different $\param$.
        Then check whether $(\mat{A}, \mat{B})$ can satisfy the (NC),
        either graphically or via the formulas
   	 	in \eqref{ExactBoundary}: is there
   	 	an admissible range of $\param$ values that guarantees
   	 	$\Lambda(\mat{A}, \mat{B}) \subseteq \mathcal{D}$?
   	 	(If not, then
   	 	unconditional stability is not possible for $(\mat{A}, \mat{B})$.)
   	 	\item Now use the sufficient conditions (SC) to determine $\param$.
            \vspace{-.5\topsep}
   	 		\begin{itemize}\itemsep.2em\parskip0em
	   	 		\item Choose a $p \in \mathbb{R}$, (try first $p = 1$).
			  		Compute $\Wp$ from equation \eqref{Eq:NumRange},
			  		for instance, using a software such as
			  		\textsf{Chebfun} \cite{DriscollHaleTrefethen2014}.
			  	\item By varying $0 < \param \leq 1$, find the largest
                    $\param$ that ensures $\Wp \subseteq \mathcal{D}$,
				  	and guarantees unconditional stability
                    ($\mathcal{D}$ becomes larger as $\param$ decreases).
				  	Call this parameter $\param^*$.
   	 		\end{itemize}   	 		
  		\item If no value $0 < \param \leq 1$ can be found in Step~2,
		or $\param^*$ is prohibitively small (leading to a large error
		constant), try and repeat Step~2 with a different $p$.
		\item Choose a $\param < \param^*$
        (e.g.\ $\param = 0.95 \, \param^*$, with 0.95 for robustness),
        and substitute it into
		Table~\ref{Table:ImexCoeff} to obtain the \imex
		coefficients for the ODE solver.
	\end{enumerate}
\end{recipe}
\end{tcolorbox}


\begin{example}\label{Ex:SimpleCase} (Simple example using the
	Recipe~\ref{Recipe:ChooseParameter})
    Consider the ODE $u_t = -10u$, with implicit part $\mat{A}u = -u$
	and explicit part $\mat{B} u = -9u$ (this ODE splitting was also examined in
    \cite{RosalesSeiboldShirokoffZhou2017}), for which we wish to devise
    a 3rd order ($r = 3$) unconditionally stable scheme.
    For this splitting,
	the matrices $\mat{A}$ and $\mat{B}$ are (trivially)
	simultaneously diagonalized with $\Lambda(\mat{A}, \mat{B}) = \{ -9\}$.
	Condition~\ref{Cond:CommMatrices} then requires $\{-9\} \in \mathcal{D}$
	for both the (NC) and (SC). For a 3rd order scheme, $r = 3$, we use the formulas
	for $m_r$ and $m_l$ in (B2) so that the constraint reads:
	\begin{align}\label{Eq:Example}
		m_l < -9 < m_r \quad \Longrightarrow \quad
	    -\frac{ (2-\param)^3}{8 - (2-\param)^3} < -9 <
		\frac{(2-\param)^3}{(2-\param)^3 + 1}.
	\end{align}
	The largest $\param$ value that satisfies the inequality
	\eqref{Eq:Example} (with $<$ replaced by $\leq$) is:
	$\param^* = 2 - (7.2)^{1/3}$.  Any value
	$0 < \param < \param^*$ will guarantee unconditional
	stability --- i.e.\ one could take a fraction $\param = 0.95\, \param^*$ so that
	$\param \approx 0.0656$. Substituting
	this value into the formulas in Table~\ref{Table:ImexCoeff}
	yields the \imex coefficients.
\end{example}


In situations where one is using a pre-programmed ODE or black-box solver,
it may not be possible to modify the time stepping coefficients
$(a_j, b_j, c_j)$. Instead, one may have the ability to modify the
matrix splitting $(\sfact \mat{A}_0, \mat{B})$ by varying the parameter
$\sfact$. Recipe~\ref{Recipe:ChooseSplitting} outlines how one may choose
the parameter $\sfact$ when the scheme and the matrix $\mat{A}_0$ are fixed.
The recipe uses the sets $\Lambda(\sfact \mat{A}_0, \mat{B})$ and
$W_p(\sfact \mat{A}_0, \mat{B})$, whose dependence on $\sfact$
is characterized by the following remark.

\begin{remark}\label{Rmk:DependenceOnSigma}
	(Dependence of $W_p(\sfact \mat{A}_0, \mat{B})$
	and $\Lambda(\sfact \mat{A}_0, \mat{B})$ on $\sfact$)
	The sets $W_p(\sfact \mat{A}_0, \mat{B})$ and
	$\Lambda(\sfact \mat{A}_0, \mat{B})$ are simple transformations of the
    $\sfact$-independent sets $W_p(\mat{A}_0, \mat{L})$ and
	$\Lambda(\mat{A}_0, \mat{L})$:
	\begin{align} \label{Eq:EVScaling}
		\Lambda(\sfact \mat{A}_0, \mat{B})
		= 1 + \sfact^{-1} \Lambda(\mat{A}_0, \mat{L}),
		\quad  \quad \quad
		W_p(\sfact \mat{A}_0, \mat{B}) = 1 + \sfact^{-1} W_p(\mat{A}_0, \mat{L}).
	\end{align}
	Here the identities \eqref{Eq:EVScaling} follow
	from a direct calculation using $\mat{B} = \mat{L} - \sfact \mat{A}_0$:
	\begin{align}\label{Eqn:EVProperty}
		(-\sfact \mat{A}_0)^{-1} \mat{B} &=
		\mat{I} + \sfact^{-1} (-\mat{A}_0)^{-1} \mat{L}.\\ \label{Eqn:WpProperty}
		(-\sfact \mat{A}_0)^{\frac{p}{2}-1} \mat{B} (-\sfact \mat{A}_0)^{-\frac{p}{2}}
		&= \mat{I} + \sfact^{-1} (-\mat{A}_0)^{\frac{p}{2}-1}
		\mat{L} (-\mat{A}_0)^{-\frac{p}{2}}.
	\end{align}
	Due to properties \eqref{Eq:EVScaling}, one can,
    for fixed $\mat{A}_0$ and $\mat{L}$, pre-compute the sets
	$\Lambda(\mat{A}_0,\mat{L})$ and $W_p(\mat{A}_0, \mat{L})$.
    The range $W_p(\sfact \mat{A}_0, \mat{B})$ and generalized eigenvalues
	$\Lambda(\sfact \mat{A}_0, \mat{B})$ are then simply rescaled versions
	(w.r.t.~the point 1 in the complex plane) of the corresponding range
	and eigenvalues using $\mat{A}_0$ and $\mat{L}$, where $\sfact$ yields
    the scaling parameter. This becomes important in \Srm\ref{Sec_SBDF}
    when we examine and overcome the fundamental limitations of SBDF.	
\end{remark}


\begin{tcolorbox}
\begin{recipe}\label{Recipe:ChooseSplitting}
    Given a fixed \imex scheme and matrix $\mat{A}_0$, how to choose
    the splitting parameter $\sfact$ for the splitting $(\sfact \mat{A}_0, \mat{B})$.
    \vspace{-.5\topsep}
	\begin{enumerate}[1.]\itemsep0em\parskip.3em\setcounter{enumi}{-1}
		\item Choose an order $1 \leq r \leq 5$; and retrieve the formulas
		for $\mathcal{D}$ in equation \eqref{ExactBoundary}. If the time stepping
		scheme being used is not included as one from Table~\ref{Table:ImexCoeff},
		then an unconditional stability diagram $\mathcal{D}$ will need to be
		computed.
		\item Compute/plot
   	 	the generalized eigenvalues $\Lambda(\sfact \mat{A}_0, \mat{B})$
		for different $\sfact$ (see Remark~\ref{Rmk:DependenceOnSigma}).
        Then check the (NC), either graphically or via the formulas
   	 	in \eqref{ExactBoundary}:
        is there an admissible range of $\sfact$ that guarantees
   	 	$\Lambda(\sfact \mat{A}_0, \mat{B}) \subseteq \mathcal{D}$?
        (If not, then unconditional stability is not possible, and
        a different \imex scheme or matrix $\mat{A}_0$ must be used.)
   	 	\item Now use the sufficient conditions (SC) to determine $\sfact$.
            \vspace{-.5\topsep}
   	 		\begin{itemize}\itemsep.2em\parskip0em
   	 			\item Choose a $p \in \mathbb{R}$ (try first $p = 1$)
			  		and compute $W_p(\sfact \mat{A}_0, \mat{B})$
					(see Remark~\ref{Rmk:DependenceOnSigma}).
				\item Vary $\sfact$ to find the smallest $\sfact > 0$
                    that ensures
                    $W_p(\sfact \mat{A}_0, \mat{B}) \subseteq \mathcal{D}$
                    and guarantees unconditional stability.
                    ($W_p(\sfact \mat{A}_0, \mat{B})$ becomes larger as
                    $\sfact$ decreases).
			\end{itemize}
  		\item If no value of $\sfact > 0$ can be found in Step~2,
			repeat Step~2 with a different $p$.
	\end{enumerate}
\end{recipe}
\end{tcolorbox}
Section~\ref{Sec_SBDF} provides examples that illustrate Recipe~\ref{Recipe:ChooseSplitting}. Recipes~\ref{Recipe:ChooseParameter} and~\ref{Recipe:ChooseSplitting} are in line with a common perspective on \imex schemes. Either, one has to determine the \imex parameter $\param$ when the matrix splitting is fixed; or choose the splitting parameter $\sfact$ when the scheme is fixed. In practice, there may be cases in which neither of these two approaches is able to achieve unconditionally stability.

We therefore advocate, whenever possible, to allow to simultaneously vary the \imex parameter $\param$ and the splitting parameter $\sfact$. It turns out that this yields an enormous amount of flexibility when designing unconditionally stable schemes. Many splittings of the form $(\sfact \mat{A}_0, \mat{B})$, where
$\mat{A}_0$ and $\mat{L}$ are chosen and predetermined from the problem (see Sections~\ref{Sec_NumExperiments}--\ref{Sec:NavierStokes} for specific PDE applications), can be stabilized this way.


\begin{tcolorbox}
\begin{recipe}\label{Recipe:ChooseBothImexAndSplit}
	Given a matrix $\mat{A}_0$, how to simultaneously choose
	both the \imex and splitting parameters $(\param, \sfact)$
	(with $0 < \param \leq 1$, $\sfact > 0$).
    \vspace{-.5\topsep}
	\begin{enumerate}[1.]\itemsep0em\parskip.3em
		\item Repeat Steps~0--1 in
		Recipes~\ref{Recipe:ChooseParameter}--\ref{Recipe:ChooseSplitting}
		to ensure that there is a range of values $(\param, \sfact)$
		that satisfy the necessary conditions (NC)
		$\Lambda(\sfact \mat{A}_0, \mat{B}) \subseteq \mathcal{D}$. Note:
		$\Lambda(\sfact \mat{A}_0, \mat{B})$ depends solely on $\sfact$,
		while $\mathcal{D}$ depends solely on $\param$.
   	 	\item Use the sufficient conditions (SC) to determine $(\param, \sfact)$.
            \vspace{-.5\topsep}
   	 		\begin{itemize}\itemsep.2em\parskip0em
   	 			\item Choose a $p \in \mathbb{R}$ (try first $p = 1$)
			  		and compute $W_p(\sfact \mat{A}_0, \mat{B})$
					(see Remark~\ref{Rmk:DependenceOnSigma}).
				\item The sufficient condition
				$W_p(\sfact \mat{A}_0, \mat{B}) \subseteq \mathcal{D}$
				provides a constraint on the parameters $(\param, \sfact)$
				that achieve unconditional stability. Within this constrained
                set, determine the points $(\param^*, \sfact^*)$
				that maximize $\param^*$. If there is more than one solution,
                choose $\sfact^*$ small.
			\end{itemize}
  		\item If no value of $(\param, \sfact)$ can be found in Step~2,
			repeat Step~2 with a different $p$.
	\end{enumerate}
\end{recipe}
\end{tcolorbox}
Sections~\ref{Sec_NumExperiments}--\ref{Sec:NavierStokes} provide specific applications of Recipe~\ref{Recipe:ChooseBothImexAndSplit} in PDE problems.

\subsection{Additional details for PDEs: choosing $\mat{A}_0$}\label{Subsec_PDEDetails}

When $\mat{L}_h$ arises as the spatial discretization
of a PDE with meshsize $h$,
one does not have a fixed matrix splitting $(\mat{A}, \mat{B})$, or
$(\sfact \mat{A}_0, \mat{B})$,
but rather a family of splittings parameterized by $h$:
$(\mat{A}_h, \mat{B}_h)$, or $(\sfact \mat{A}_{0,h}, \mat{B}_h)$.
In this situation, it is crucial to be able to choose the
parameters $(\param, \sfact)$ independent of the
meshsize $h$ --- i.e.\ to have one and the same \imex scheme be unconditionally
stable for an entire family of splittings $(\mat{A}_h, \mat{B}_h)$, or
$(\sfact \mat{A}_{0,h}, \mat{B}_h)$.
If, for example, unconditional stability required one to choose
the \imex parameter $\param$ as a function of the grid size
$h$ (i.e.\ such as $\param = h$),
then such a choice would have a deleterious effect on the GTE
(GTE $\sim \mathcal{O}( h^{-r} \delt^r)$), and limit the benefits of
unconditional stability.


To be able to choose a single set of parameters $(\param, \sfact)$
that stabilizes the family of splittings $(\sfact \mat{A}_{0,h}, \mat{B}_h)$
for all $h$, some care must be taken to
ensure the matrix $\mat{A}_{0,h}$ is properly chosen relative
to $\mat{B}_h$.  Once a suitable choice of $\mat{A}_{0, h}$ is fixed,
one may use the Recipe~\ref{Recipe:ChooseBothImexAndSplit}
to simultaneously choose $(\sfact, \param)$ for unconditional stability.
\begin{remark}	\label{Rmk:GuidelinesA0}
	(Guidelines for choosing $\mat{A}_{0, h}$ when $\mat{L}_h$ is the
	spatial discretization of a PDE)
	Generally speaking, it is a good idea to ensure that
	$\mat{A}_{0,h}$ has the same derivative order as
	$\mat{L}_h$, as backed up by the following heuristic scaling argument.
	Suppose
	\[
		\mat{L}_h \approx C(x) \frac{\partial^q}{\partial x^q} +
		\text{(lower order derivatives)}.
	\]
	A natural choice for $\mat{A}_{0,h}$ might be
	\[
		\mat{A}_{0,h} \approx \frac{\partial^s}{\partial x^s}
	\]
	(one could include a variable coefficient approximation as well).
    Many spatial approximation methods yield the scaling
    $\frac{\partial}{\partial x} \propto h^{-1}$,
    hence one may expect some of the
	eigenvalues of $(\mat{A}_{0,h})^{-1} \mat{L}_h$ to scale like $\mathcal{O}(h^{s-q})$.
	The re-scaling formulas in Remark~\ref{Rmk:DependenceOnSigma} then imply that
	there may be generalized eigenvalues $\mu \in \Lambda(\sfact \mat{A}_{0,h}, \mat{B}_h)$
	that scale like $\mu = 1 - \sfact^{-1} \mathcal{O}(h^{s-q})$.
	This gives rise to three cases for choosing $s$:
	\begin{itemize}
		\item If $s < q$, some of the generalized eigenvalues diverge
		$\mu = 1 - \sfact^{-1} \mathcal{O}(h^{s-q}) \rightarrow \infty$
		as $h\rightarrow 0$. Using formula (B3) for the asymptotic
		behavior of $\mathcal{D}$, the
		\imex parameter $\param$ would then have to scale like
		$\param \sim h^{q - s}$ as $h\rightarrow 0$ (and fixed $\sfact$)
		to ensure that these large eigenvalues remain inside $\mathcal{D}$
		(to satisfy the (NC)).
		Hence, $(\param, \sfact)$ cannot be chosen independent of the mesh $h$.
		\item If $s > q$, some of the generalized eigenvalues
		$\mu = 1 - \sfact^{-1} \mathcal{O}(h^{s-q}) \rightarrow 1$
		as $h\rightarrow 0$ (and fixed $\sfact$).
		In this case, the formulas in (B2) show that only order $r = 1, 2$
		schemes contain the point $1 \in \mathcal{D}$
		(see Figure~\ref{Fig:NewImExStabilityRegions}).
		Hence, $s > q$ is generally not a good
		choice if one is looking for a scheme with orders $r > 2$.
		\item If $s = q$, then all generalized eigenvalues $\mu$ have a chance
		(based solely on the scaling of $h$)
		to be uniformly bounded
		(i.e.\ do not become arbitrarily large) as $h \rightarrow 0$;
		and also remain strictly bounded away from $1$ as $h \rightarrow 0$.
		In this	case, there is a chance to obtain high order by means
        of choosing the parameters $(\param, \sfact)$ independent of $h$.
	\end{itemize}
\end{remark}

%
\section{Limitations of unconditional stability for SBDF schemes}\label{Sec_SBDF}

In \Srm\ref{Sec_Review}, the unconditional stability region $\mathcal{D}$
was used to derive sufficient (SC) and necessary (NC) conditions for
unconditional stability.
Using these conditions, this section illustrates how the
geometrical properties of $\mathcal{D}$ can be used to understand the
fundamental limitations that classical SBDF methods possess with regards to
unconditional stability. Specifically, two significant qualitative
transitions occur:
(i) moving from 1st to 2nd order schemes for non-symmetric matrices $\mat{L}$; and
(ii) moving from 2nd to 3rd order for symmetric matrices $\mat{L}$.
Guided by Recipe~\ref{Recipe:ChooseSplitting}, we discuss
under which circumstances a choice of $\sfact$ exists so that a splitting
$(\sfact \mat{A}_0, \mat{L})$ is unconditionally stable with SBDF.

\begin{figure}[htb!]
\centering
\includegraphics[width = .98\textwidth]{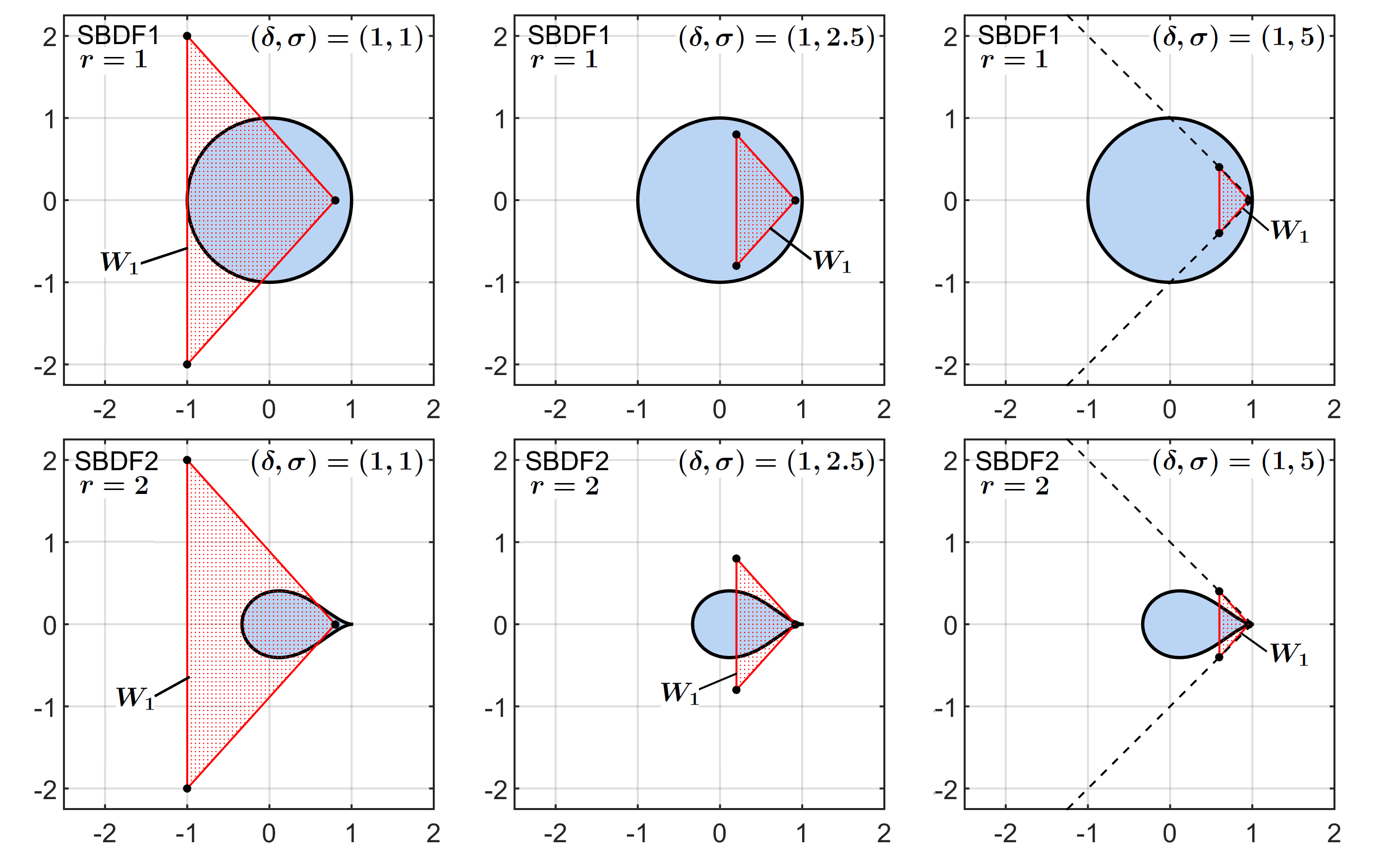}
\caption{Example for non-symmetric $\mat{L}$ in \eqref{Ex_SkewSym_L}.
   The figures show the SBDF1 (top row) and SBDF2 (bottom row) stability
   diagrams $\mathcal{D}$ (blue shaded region) in relation to the
   sets $W_1(\sfact\mat{A}_0, \mat{B})$ (red shaded region,
   abbreviated as $W_1$)
   and $\Lambda(\sfact \mat{A}_0, \mat{B})$ (black dots) for
   (left to right) $\sfact \in \{1, 2.5, 5\}$.  Note that
   $W_1(\sfact\mat{A}_0, \mat{B}) \subseteq \mathcal{D}$
   for SBDF1 with $\sfact \in \{2.5, 5\}$ guaranteeing the (SC) for unconditional
   stability.
   The bottom row highlights the fundamental limitation for SBDF2:
   no $\sfact > 0$ exists that can ensure
   $\Lambda(\sfact \mat{A}_0, \mat{B}) \subseteq \mathcal{D}$.
   Dashed lines show the effect of the rescaling by $\sfact$,
   outlined in Remark~\ref{Rmk:DependenceOnSigma},
   on the set $W_1(\sfact \mat{A}_0, \mat{B})$. }
\label{Ex_SkewSym_L_diagram}
\end{figure}

\medskip
\noindent\textbf{Case 1: $\mat{L}$ non-symmetric.}\quad
	Let $\mat{L}$ be a non-symmetric matrix that, together with $\mat{A}_0$,
	has both a range
	$\text{Re}\big( W_p(\mat{A}_0, \mat{L}) \big) < 0$
	and eigenvalues $\text{Re}\big( \Lambda(\mat{A}_0, \mat{L}) \big) < 0$
	with negative real part, i.e.~they lie strictly in the left-half plane,
	but are not necessarily contained on the real line. Such a situation occurs for
	instance in discretizations of advection--diffusion PDEs (with an implicit 
	diffusion, and explicit advection). The following
	transition arises between first and second order SBDF when the \imex
	splitting is taken as $(\sfact \mat{A}_0, \mat{B})$:
	\begin{enumerate}
	\item SBDF1 can always be made unconditionally stable, by choosing $\sfact$
		suitably large. This is due to the fact that $\mathcal{D}$ for SBDF1 is a
		circle with its right-most point at $1$. Hence one can always rescale
		$W_p(\mat{A}_0, \mat{L})$ (see Remark~\ref{Rmk:DependenceOnSigma}) so that
		$W_p(\sfact \mat{A}_0, \mat{B}) \subseteq \mathcal{D}$.
	\item SBDF2 can, in general, \emph{not} be made unconditionally stable by
		means of choosing $\sfact > 0$. This is a result of the cusp at $1$ in
		$\mathcal{D}$ (see Figure~\ref{Ex_SkewSym_L_diagram}).
		If, for instance, the imaginary
		part of $\mu \in \Lambda( \mat{A}_0, \mat{L} )$ is larger
		(in absolute value) than its real
		part, then the scaled eigenvalue (see Remark~\ref{Rmk:DependenceOnSigma})
		$1+ \sfact^{-1} \mu \in \Lambda( \sfact \mat{A}_0, \mat{L} )$
		will never enter $\mathcal{D}$, regardless of the value of $\sfact$.
	\end{enumerate}
	We highlight these insights with the following simple example.

\begin{figure}[htb!]
\centering	
\includegraphics[width = .98\textwidth]{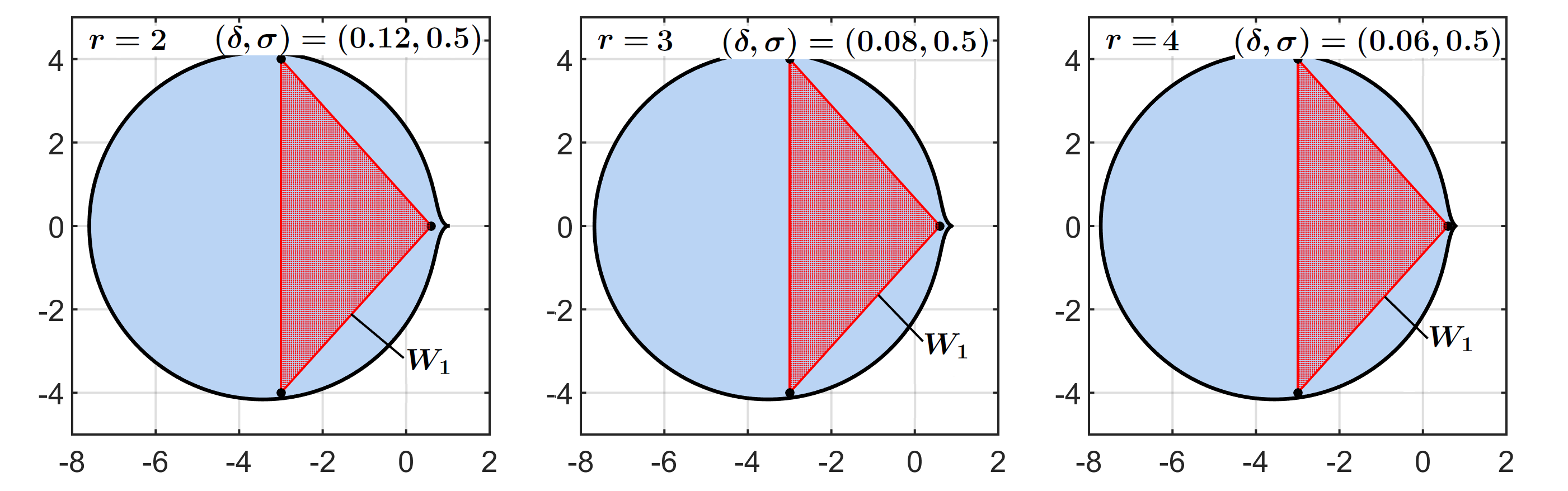}
\caption{Unconditionally stable schemes for \eqref{Ex_SkewSym_L}, and
	orders (left to right) $r \in \{2, 3, 4\}$.
	The figure shows that $W_1(\sfact\mat{A}_0, \mat{B}) \subseteq \mathcal{D}$
	when $\sfact = 0.5$, guaranteeing the (SC) for unconditional stability.
	Values are (left to right)
	$\param \in \{0.12, 0.08, 0.06\}$.
	The chosen $(\param, \sfact)$-values are guided by
	Recipe~\ref{Recipe:ChooseBothImexAndSplit}, and are almost optimal,
	however other values are also possible.}
\label{Ex_SkewSym_L_diagram2}
\end{figure}

\begin{example} (A non-symmetric $\mat{L}$)
	Consider the following non-symmetric matrix $\mat{L}$ and choice of
	matrix $\mat{A}_0$:
	\begin{align}\label{Ex_SkewSym_L}
		\mat{L} = \left(\begin{array}{ccc}
			-0.2 & \phantom{-}0 & \phantom{-}0 \\
			 \phantom{-}0  & -2 & \phantom{-}2 \\
			 \phantom{-}0 &  -2 & -2
		\end{array}\right),\quad
		\mat{A}_0 = -\left(\begin{array}{ccc}
			 1 & 0 & 0 \\
			 0 & 1 & 0 \\
			 0 &  0 & 1
		\end{array}\right).
	\end{align}
	The generalized eigenvalues are
	$\Lambda( \mat{A}_0, \mat{L} ) =
	\{ -0.2, -2 + 2\imath, -2 - 2\imath\}$, and
	$W_{1}( \mat{A}_0, \mat{L}) = \mathrm{conv}\{ -0.2, -2 + 2\imath, -2 - 2\imath\}$
	is a triangle consisting of the convex hull\footnote{The matrix $\mat{L}$
	is normal, which results in a simple expression for the range
	$W_1(\mat{A}_0, \mat{L})$.} of the eigenvalues.
	Note that, for simplicity, we have chosen an example in which
    $\mat{A}$ and $\mat{B}$	commute. Hence, Condition~\ref{Cond:CommMatrices}
	may be used.  We also plot
	$W_1(\mat{A}, \mat{B})$ to illustrate how to apply
	Condition~\ref{Cond:NoCommMatrices} (which is stronger than
	Condition~\ref{Cond:CommMatrices}) when one is faced with
	matrices $\mat{A}, \mat{B}$ that do not commute.
	Figure~\ref{Ex_SkewSym_L_diagram} visualizes that SBDF1 can be made unconditionally
	stable with $\sfact = 2.5$, while SBDF2 cannot be made unconditionally
	stable by only varying $\sfact > 0$. However, high order schemes (i.e.\ $r \geq 2$) that
	are unconditionally stable for (\ref{Ex_SkewSym_L}) \emph{are} possible by varying
	both $(\param, \sfact)$, as seen in Figure~\ref{Ex_SkewSym_L_diagram2}.
\end{example}


\begin{figure}[htb!]
\centering
\includegraphics[width = .98\textwidth]{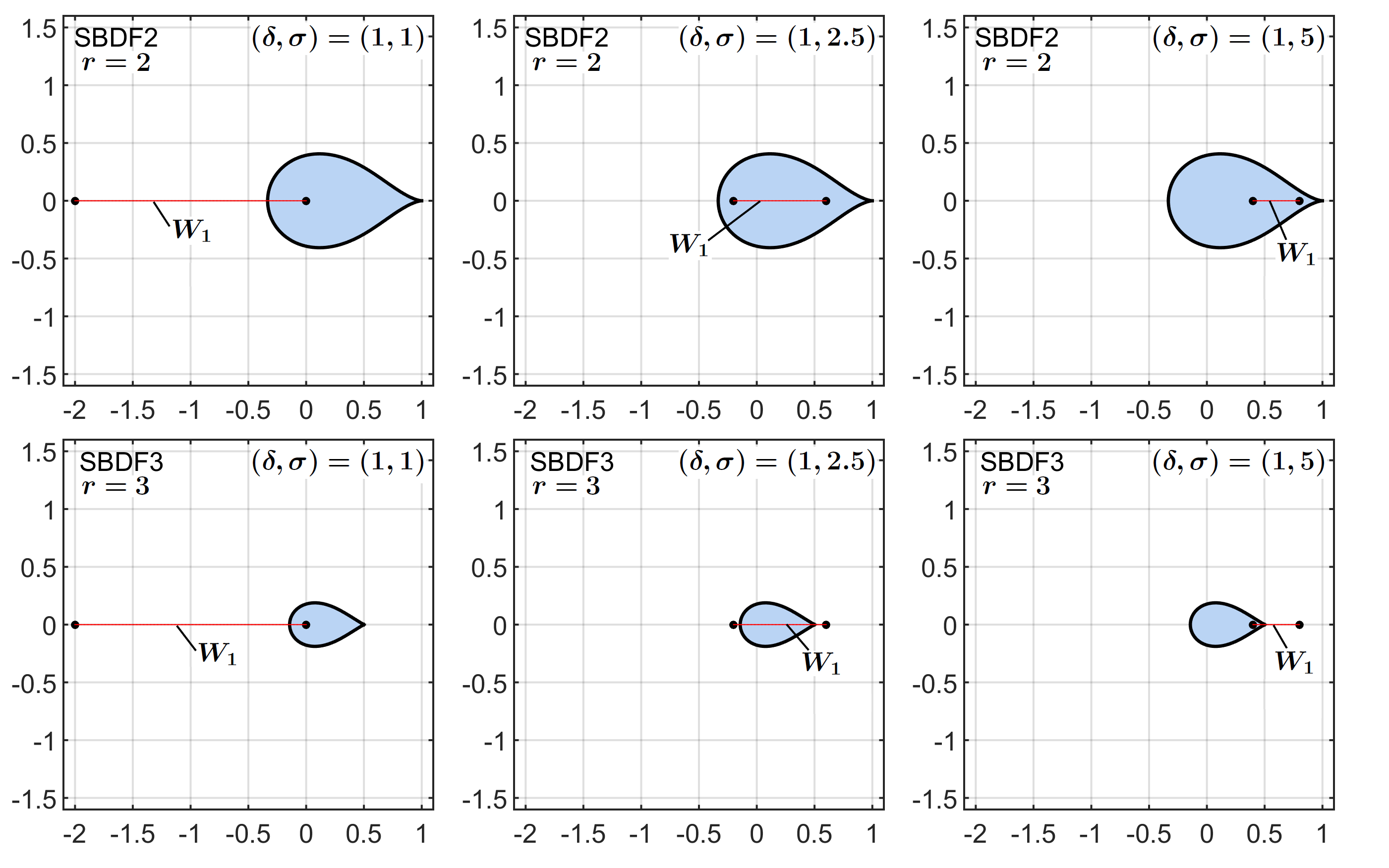}
\caption{Example for symmetric $\mat{L}$ in \eqref{Ex_Sym_L}.
   The figures show the SBDF2 (top row) and SBDF3 (bottom row) stability
   diagrams $\mathcal{D}$ (blue shaded region) in relation to the
   sets $W_1(\sfact \mat{A}_0, \mat{B})$ (red region, abbreviated $W_1$)
   and $\Lambda(\sfact \mat{A}_0, \mat{B})$ (black dots) for
   (left to right) $\sfact \in \{1, 2.5, 5\}$.  Note that
   $W_1(\sfact \mat{A}_0, \mat{B}) \subseteq \mathcal{D}$
   for SBDF2 with $\sfact \in \{2.5, 5\}$, guaranteeing unconditional stability.
   The bottom row highlights the fundamental limitation for SBDF3:
   no $\sfact > 0$ exists that can ensure
   $\Lambda(\sfact \mat{A}_0, \mat{B}) \subseteq \mathcal{D}$. }
\label{Ex_Sym_L_diagram}
\end{figure}

\begin{figure}[htb!]
\centering
\includegraphics[width = .98\textwidth]{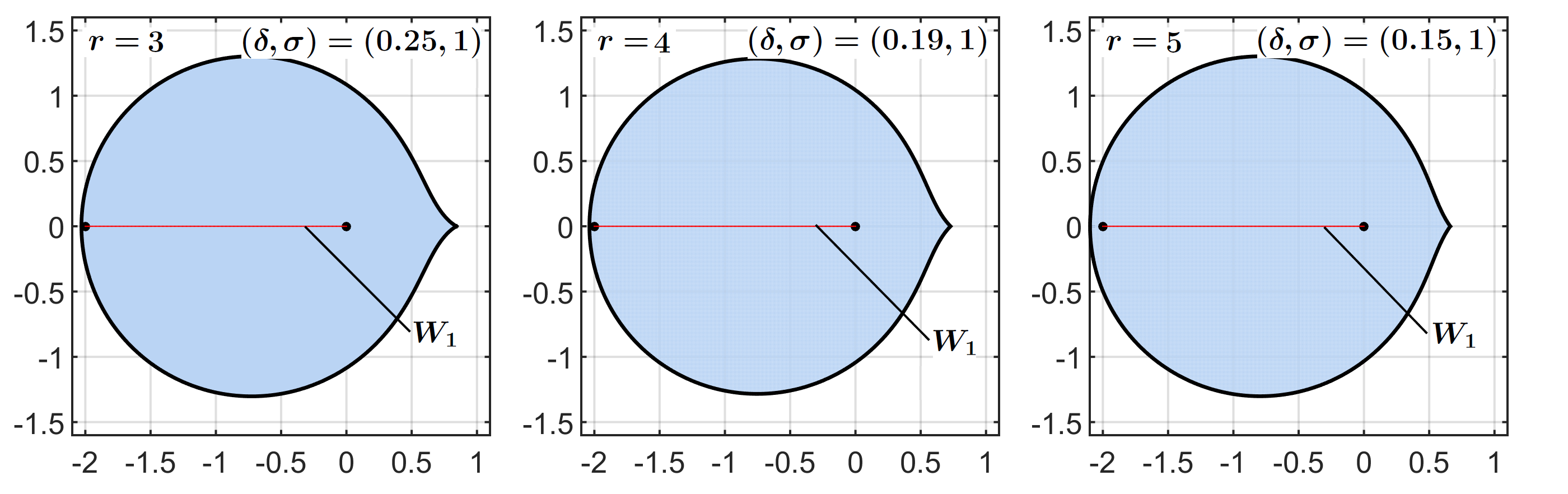}
\caption{Unconditionally stable schemes for (\ref{Ex_Sym_L}), and
    orders (left to right) $r \in \{3, 4, 5\}$.
	The figure shows that $W_1(\sfact \mat{A}_0, \mat{B}) \subseteq \mathcal{D}$
	provided $\sfact = 1$
	and (left to right) $\param \in \{0.25, 0.19, 0.15\}$.
	The chosen $(\param, \sfact)$-values are guided by
	Recipe~\ref{Recipe:ChooseBothImexAndSplit}, however
	other values are also possible.}
\label{Ex_Sym_L_diagram2}
\end{figure}


\medskip
\noindent\textbf{Case 2: $\mat{L}$ symmetric.}\quad
Let $\mat{L}$ be a symmetric negative definite matrix. Assume
    now that $\mat{A}_0$ is such that the range of $W_1(\mat{A}_0, \mat{L})$
    and eigenvalues $\Lambda(\mat{A}_0, \mat{L})$ are
    real and strictly negative.
Such a situation arises for instance in the discretization of a purely
parabolic (gradient flow) problem. The following transition occurs between
second- and third-order schemes for the splittings $(\sfact \mat{A}_0, \mat{B})$:
\begin{enumerate}
	\item SBDF2 can always be made unconditionally stable, by choosing $\sfact$
	suitably large. This is due to the fact that the right-most point of
	$\mathcal{D}$ for SBDF2 is $1$, and $W_1(\mat{A}_0,\mat{L})$ is real and
	negative, so one can always rescale
	(see Remark~\ref{Rmk:DependenceOnSigma}) $W_1(\sfact\mat{A}_0,\mat{B})$
	into $\mathcal{D}$.
	\item SBDF3 can, in general, \emph{not} be made unconditionally stable
	by means of choosing $\sfact > 0$. This is because the
	right-most point of $\mathcal{D}$ is $1/2$ (instead of $1$), so that a
	negative real $\Lambda(\mat{A}_0,\mat{L})$ may be impossible to contain within
	$\mathcal{D}$ via the choice of $\sfact$.
\end{enumerate}
Unconditional stability limitations of SBDF, applied to splittings
$(\sfact \mat{A}_0, \mat{B})$, may be overcome by simultaneously
choosing $(\param, \sfact)$, i.e.~by following
Recipe~\ref{Recipe:ChooseBothImexAndSplit}.


\begin{example}	(A symmetric $\mat{L}$)
	Consider the following symmetric matrices:
	\begin{align}\label{Ex_Sym_L}
		\mat{L} = \left(\begin{array}{ccc}
			-2 & \phantom{-}1  \\
			 \phantom{-}1  & -2
		\end{array}\right),\quad
		\mat{A}_0 = -\left(\begin{array}{ccc}
			 1 & 0 \\
			 0 & 1
		\end{array}\right).
	\end{align}
	Then
	$\Lambda( \mat{A}_0, \mat{L} ) = \{ -3, -1\}$, and
	$W_1(\mat{A}_0, \mat{L}) = [ -3, -1 ]$ is an interval along the real axis.
	Figure~\ref{Ex_Sym_L_diagram} shows that SBDF2 can be made unconditionally
	stable with $\sfact = 2.5$, while SBDF3 cannot be made unconditionally
	stable by only varying $\sfact > 0$. In contrast, third to fifth order schemes that
	are unconditionally stable for \eqref{Ex_Sym_L} \emph{are} possible by varying
	both $(\param, \sfact)$, as seen in Figure~\ref{Ex_Sym_L_diagram2}.
\end{example}



\section{Examples from diffusion PDEs}\label{Sec_NumExperiments}

In this section we apply the unconditional stability theory from \Srm\ref{Subsec_PDEDetails} and Recipe~\ref{Recipe:ChooseBothImexAndSplit} to PDE diffusion problems with spatially varying, and even non-linear diffusion coefficients.
The presented methodology highlights how one can avoid a stiff time step restriction (here: of diffusive type $\delt \propto h^2$, where $h$ is the smallest grid size) --- or any time step restriction for that matter --- while inverting only simple constant coefficient matrices. This allows one to leverage \emph{fast solvers} where an implicit treatment of $\mat{A}_h$ (i.e.~using the fast Fourier transform) can be carried out much more rapidly than a fully implicit treatment of $\mat{L}_h$ (i.e.~that contains all the stiff diffusive terms).
The new \imex coefficients (see \Srm\ref{Sec_SBDF}) enable high order time stepping beyond what is possible using only SBDF methods.

\subsection{Numerical discretization}\label{Subsec:discretizations}
We start by providing numerical details for the one-dimensional Fourier spectral methods used. Computations in three dimensions are then conducted by naturally extending the one-dimensional approach via Cartesian products.
We use a periodic computational domain $\varOmega = [0, 1]$; discretize space
using a uniform grid with an even number of grid points $\sizeN$;
and approximate the function $u(x)$ at $x_j$ by $u_j \approx u(x_j)$, where:
\[ x_j = j h, \quad h = \frac{1}{\sizeN}, \quad
	\vec u = \begin{pmatrix}
		u_1, u_2, \ldots, u_{\sizeN}
	\end{pmatrix}^T \in \mathbb{R}^{\sizeN}.
\]

Because our analysis is based on the matrices $\mat{A}_{0,h}$ that are written in terms of Fourier transforms, it is useful to introduce notation for the discrete Fourier transform (DFT) matrix $\mat{F}$, and for the spectral differentation matrix $\mat{D}$ --- even though in practice one will never use those matrices, but rather use the \emph{fast Fourier transform} (FFT) to compute $\mat{F}\vec{u} = \mathrm{fft}(\vec{u})$. The DFT matrix $\mat{F}$ has the coefficients:
\begin{align*}
	\mat F_{j\ell} = \omega^{(j-1) \times (\ell-1)},
\quad \omega = e^{-\frac{2\pi \imath }{\sizeN}}, \quad \text{so that} \quad
	(\mat{F}\vec{u})_j = \sum_{\ell = 1}^N u_\ell \; \omega^{(j-1)\times (\ell-1)}.
\end{align*}
The (spectral) differentiation of a function defined on the uniform grid amounts to a scalar multiplication in Fourier space, i.e.~$(\mat{D}\vec{u} )_j \approx u_x(x_j)$.
Hence, the matrix $\mat{D}$ takes the form: $\mat D = \imath \mat F^{-1} \, \text{diag}( \vec{\xi} ) \, \mat F$, where $\diag(\vec{\xi})$ denotes the matrix with diagonal entries of the vector:
\begin{align}\label{Wavenumbers}
	\vec{\xi} = \begin{pmatrix}
		\xi_1, \xi_2, \ldots,  \xi_{\sizeN}
	\end{pmatrix}^T \in \mathbb{R}^{\sizeN}, 	
	\quad \text{where} \quad
\xi_j = \left\{
\begin{array}{ll}
2\pi (j-1) & \text{if } 1 \leq j \leq \frac{\sizeN}{2},\\
2\pi (j - \sizeN) & \text{if } \frac{\sizeN}{2}+1 < j \leq \sizeN,\\
\sizeN \pi & \text{if } j = \frac{\sizeN}{2}+1.
\end{array} \right.
\end{align}

Since $\mat{F}^{-1} = \sizeN^{-1} \mat{F}^{\dag}$, the matrix $\mat{D}^{\dag} = -\mat{D}$
is skew-Hermitian and the matrix $\mat{D}^2$ is Hermitian.
If $\mat{A}_h$ is diagonalized by $\mat{F}$, then solving
for $\vec{u}_{n + r}$
in the implicit step of the evolution \eqref{Eq:fulltimestepping},
i.e.~$(a_r \mat I - \delt c_r \mat A_h ) \vec{u}_{n+r} = \text{RHS}$,
is done via two FFTs.

\subsection{An FFT-based treatment for the variable coefficient diffusion equation}\label{subsection_vardiff}\label{Sec_VarDiffusion}
We now devise unconditionally stable \imex schemes for the variable coefficient diffusion equation (with diffusion coefficient $d(x) > 0$)
\begin{align}\label{Eq:VarDiffPDE}
	u_t = \big( d(x) u_x \big)_x + f(x, t),
	\quad \quad \text{on } \varOmega \times (0, T],
\end{align}
that make use of an FFT-based treatment of the implicit matrix $\mat{A}_h$.
The choice of splitting $(\mat{A}_h, \mat{B}_h)$ is guided by \Srm\ref{Subsec_PDEDetails}, and the choice of parameters $(\param, \sfact)$ by Recipe~\ref{Recipe:ChooseBothImexAndSplit}.
To ensure a high spatial accuracy, we adopt a spectral discretization of equation \eqref{Eq:VarDiffPDE} and set:
\[
	\mat{L}_h = \mat{D} \!\big( \diag(\vec{d}) \big)\! \mat{D}, \quad
	\text{where} \quad
	\vec{d} =
	\begin{pmatrix}
		d(x_1), d(x_2), \ldots, d(x_{\sizeN})
	\end{pmatrix}^T.
\]
Note that $\mat{L}_h$ is a dense matrix and (due to the $x$-dependence of
$d(x)$) is not diagonalized via the DFT matrix $\mat{F}$. To seek
an \imex splitting of $\mat{L}_h$, we follow the guidelines in
Remark~\ref{Rmk:GuidelinesA0}: the matrix $\mat{L}_h$ has two factors of
$\mat{D}$ and hence the implicit matrix $\mat{A}_h$ should have two factors
of $\mat{D}$ as well. This motivates the following matrix splitting:
\begin{align} \label{VarDiffusion_Splitting}
	\mat{A}_h = \sfact \mat D^2,
	\quad \quad	
	\mat{B}_h = \mat D \; \Big( \diag( \vec{d} ) -
	\sfact \mat I  \Big) \; \mat D,
\end{align}
i.e.~$\mat{A}_h \vec{u} \approx \sigma u_{xx}$ and
$\mat{B}_h \vec{u} \approx \big( (d(x)-\sigma)u_x\big)_x$.

Our goal is to determine, following Recipe~\ref{Recipe:ChooseBothImexAndSplit}, the parameters $(\param, \sfact)$ that guarantee unconditional stability.
Before doing so, we must discuss a caveat: the matrices $\mat{L}_h$ and $\mat{A}_h$ are not invertible --- which was an assumption in the derivation of the conditions for unconditional stability.
We do not provide a general treatment for when $\mat{A}_h$ is not invertible due to subtleties that may arise (for instance when the null space of $\mat{A}_h$ \emph{interacts} with $\mat{B}_h$ through the \imex evolution). However, for the specific splitting \eqref{VarDiffusion_Splitting}, the unconditional stability theory presented in \Srm\ref{Sec_Review} (and recipes in \Srm\ref{SecRecipe}) can be applied with only a minor adaptation, namely: the definition/computation of the sets $W_p(\mat{A}_h, \mat{B}_h)$ and $\Lambda(\mat{A}_h, \mat{B}_h)$ are done on the subspace $\mathbbm{V}$ where $\mat{A}_h$ is invertible, as follows.

The matrix $\mat{D}$ has a null space spanned by the constant vector $\vec{1} = ( 1,\; 1,\;\ldots,\; 1 )^T$. Hence, $\mat{D}^2$ and $\mat{A}_h$ have the null space $\mathbbm{1}$, and column space (range) $\mathbbm{V}$ where:
\[
	\mathbbm{V} :=
	\{ \vec{u} \in \mathbbm{C}^{\sizeN} : \vec{1}^T \vec{u} = \vec{0} \},
	\quad	
	\mathbbm{1} := \mathrm{span}\{ \vec{1} \},
	\quad
	\text{so that}
	\quad
	\mathbb{C}^{\sizeN} = \mathbbm{V} \oplus \mathbbm{1}.
\]
Using the orthogonal projection $\mat{P} = \mat{I} - \sizeN^{-1} \vec{1}\; \vec{1}^T$ onto $\mathbbm{V}$, and noting that $\mat{A}_h$ (and also $\mat{B}_h$)
satisfies $\vec{1}^T \mat{A}_h = \mat{A}_h\, \vec{1} = \vec{0}$, so that
$\mat{A}_h = \mat{P}\, \mat{A}_h = \mat{A}_h \, \mat{P}$, the evolution equation
\begin{align}\label{Eqn:Dyanamics}
	\vec{u}_t = \mat{A}_h \vec{u} + \mat{B}_h \vec{u}
\end{align}
decouples into separate components that lie in the subspaces $\mathbbm{1}$ and $\mathbbm{V}$ (i.e.~$\mathbbm{1}$ and $\mathbbm{V}$ are invariant subspaces of equation~\eqref{Eqn:Dyanamics}):
\begin{align}\label{Eq:MeanEvolution}
\text{Dynamics in } \mathbbm{1}: \quad\quad
	\big(\vec{1}^T \vec{u} \big)_t &=
	\vec{1}^T \big( \mat{A}_h \vec{u} + \mat{B}_h \vec{u} \big) = 0.\\\label{Eq:DeMeanedEvolution}
\text{Dynamics in } \mathbbm{V}: \quad\quad
	\big( \mat{P} \vec{u} \big)_t &= \mat{P}\big(\mat{A}_h \vec{u} + \mat{B}_h \vec{u}\big)
	= \mat{A}_h \big(\mat{P}\vec{u} \big) + \mat{B}_h \big(\mat{P}\vec{u}\big).
\end{align}
Equation \eqref{Eq:MeanEvolution} shows that the mean of $\vec{u}$, i.e.~$(\vec{1}^T \vec{u})$, remains constant. Any zero-stable \imex scheme (such as the ones we use) applied to \eqref{Eqn:Dyanamics} automatically ensures that $(\vec{1}^T\vec{u})$ evolves according to \eqref{Eq:MeanEvolution} with stable growth factors (independent of $\delt$, given by $a(z) = 0$). Hence, the mean $(\vec{1}^T \vec{u})$ is unconditionally stable.
In turn, equation \eqref{Eq:DeMeanedEvolution} can be viewed as the restriction of equation \eqref{Eqn:Dyanamics} to the space $\mathbbm{V}$. Because $\mat{A}_h$ is invertible on $\mathbbm{V}$, the stability theory outlined in \Srm\ref{Sec_Review} applies to equation \eqref{Eq:DeMeanedEvolution}, where the sets $W_p(\mat{A}_h, \mat{B}_h)$ and $\Lambda(\mat{A}_h, \mat{B}_h)$ are computed on the subspace $\mathbbm{V}$ instead of $\mathbb{C}^{\sizeN}$.
To summarize the results:
\begin{remark}\label{Rmk:ModifiedW}
	(Modification of $W_p(\mat{A}_h, \mat{B}_h)$ and
$\Lambda(\mat{A}_h, \mat{B}_h)$ for a non-invertible $\mat{A}_h$)
	The splitting \eqref{VarDiffusion_Splitting} with the
	discretization in \Srm\ref{Subsec:discretizations}
	leads to a matrix $\mat{A}_h$ that is not invertible.
    This violates the assumptions for the necessary
    and sufficient conditions in \Srm\ref{Sec_Review}.
	Nevertheless, Conditions~\ref{Cond:CommMatrices}--\ref{Cond:NoCommMatrices}
	may be used, provided $W_p(\mat{A}_h, \mat{B}_h)$ is computed
    on the space $\mathbbm{V}$:
	\[
		W_p(\mat{A}_h, \mat{B}_h) = \Big\{
		\langle \vec{x}, (-\mat{A}_h)^{p-1}\mat{B}_h\vec{x}\rangle :
		\langle \vec{x}, (-\mat{A}_h)^{p}\vec{x}\rangle = 1,
		\vec{x} \in \mathbbm{V} \Big\},
	\]
    and likewise, $\mu \in \Lambda(\mat{A}_h, \mat{B}_h)$ are restricted to
	the eigenvalues with corresponding eigenvectors $\vec{v} \in \mathbbm{V}$.
\end{remark}
With a slight abuse of notation, we continue to use $W_p(\mat{A}_h, \mat{B}_h)$
and $\Lambda(\mat{A}_h, \mat{B}_h)$ throughout this section with the understanding
that they are computed only on the subspace $\mathbbm{V}$.

Owing to the simple structure of $\mat{B}_h$ in relation to
$\mat{A}_h$, we can compute (almost exactly) the (modified)
set $W_1(\mat{A}_h, \mat{B}_h)$ described in Remark~\ref{Rmk:ModifiedW},
as well as the minimum and maximum eigenvalues (the eigenvalues in this case
are real)
$\Lambda(\mat{A}_h, \mat{B}_h)$ in terms of the discrete vector
$\vec{d}$ and diffusion coefficient $d(x)$.
To do so, we introduce the notation
\[
	d_\mathrm{min} = \min_{x \in \varOmega} d(x),
	\quad d_\mathrm{max} = \max_{x \in \varOmega} d(x),
\]
as well as the discrete values
\begin{align*}
	\begin{array}{ll}
	d_{2,\mathrm{min}} = \big\{ \text{Second smallest element of } \vec{d} \big\}, &	
	\mu_\mathrm{min} = \min \big\{\mu : \mu \in \Lambda( \mat{A}_h, \mat{B}_h) \big\} , \\
	d_{2,\mathrm{max}} = \big\{ \text{Second largest element of } \vec{d} \big\}, &
	\mu_\mathrm{max} = \max \big\{ \mu : \mu \in \Lambda( \mat{A}_h, \mat{B}_h) \big\}.	
	\end{array}	
\end{align*}
The sets $W_1(\mat{A}_h, \mat{B}_h)$ and max/min values in $\Lambda(\mat{A}_h, \mat{B}_h)$ then satisfy:
\begin{proposition}\label{Prop_NumR_VarDiff}
	The set $W_1(\mat{A}_h, \mat{B}_h)$ for any splitting of the form
	\eqref{VarDiffusion_Splitting} is strictly real and contained
	inside the interval:
	\[
    	1 - \sfact^{-1} d_\mathrm{max} \leq \; W_1(\mat{A}_h, \mat{B}_h) \;
	    \leq 1 -\sfact^{-1}d_\mathrm{min}.
	\]
	Moreover, the generalized eigenvalues $\Lambda(\mat{A}_h, \mat{B}_h)$
	are all real, and are bounded by
	\[
		1 - \sfact^{-1} d_\mathrm{max} \leq \mu_\mathrm{min}
        \leq 1-\sfact^{-1}d_{2,\mathrm{max}}, \quad \quad
		1 - \sfact^{-1}d_{2,\mathrm{min}} \leq \mu_\mathrm{max}
        \leq 1 -\sfact^{-1}d_\mathrm{min}.
	\]
\end{proposition}

\begin{remark}\label{Rmk:IntuitionPropProof}
	(Motivation based on operators)
	The intuition for the proof of Proposition~\ref{Prop_NumR_VarDiff} arises
	at the continuum level of differential operators.
	Roughly speaking, one can write $\mathcal{A} = \frac{d^2}{dx^2}$
	and $\mathcal{B} = \frac{d}{dx} d(x) \frac{d}{dx}$, so
	one may expect $\mathcal{A}^{-\frac{1}{2}} \propto (\frac{d}{dx})^{-1}$.
	This yields the operator product
	$\mathcal{A}^{-\frac{1}{2}} \mathcal{B} \mathcal{A}^{-\frac{1}{2}} = d(x)$,
	which allows for the computation of $W_1(\mathcal{A}, \mathcal{B})$.
	The proof of Proposition~\ref{Prop_NumR_VarDiff} in
	\ref{Sec:AppendixProof} effectively formalizes this
	operator computation at the level of matrices.
	Moreover, due to the continuum nature of the argument,
	Proposition~\ref{Prop_NumR_VarDiff} carries over to
	other spatial discretizations, such
	as other spectral methods, finite differences, etc..
\end{remark}

Proposition~\ref{Prop_NumR_VarDiff}
is useful as it allows the design of unconditionally stable
\imex schemes by choosing $(\param, \sfact)$ so that $W_1(\mat{A}_h, \mat{B}_h) \subseteq \mathcal{D}$. It is significant for two more reasons. First, the bounds on
$W_1(\mat{A}_h, \mat{B}_h)$ and $\Lambda(\mat{A}_h, \mat{B}_h)$
in Proposition~\ref{Prop_NumR_VarDiff}
do not depend on $h$. This allows one to chose a single \imex parameter $\param$
(independent of $h$) to stabilize an entire family
of splittings $(\mat{A}_h, \mat{B}_h)$.  Second, the proposition is almost exact:
\begin{remark}(Proposition~\ref{Prop_NumR_VarDiff} is almost exact)
	Although the formulas in Proposition~\ref{Prop_NumR_VarDiff} are
	inequalities, they are almost exact.
	For smooth functions $d(x)$, the values $d_{2,\mathrm{min}}, d_{2,\mathrm{max}}$ are
	at least $\mathcal{O}(\sizeN^{-1})$ close to $d_\mathrm{min}$ and $d_\mathrm{max}$.
	Hence, the bounds for $\mu_\mathrm{min}$ or $\mu_\mathrm{max}$ are sharp to within
	$\mathcal{O}(\sizeN^{-1})$.
	In a similar fashion, it can be shown that the inequalities
	on the set $W_{1}(\mat{A}_h, \mat{B}_h)$ in Proposition~\ref{Prop_NumR_VarDiff}
	are accurate to within an error $\mathcal{O}(\sizeN^{-1})$. 	
\end{remark}

\begin{figure}[htb!]
	\centering
	\begin{minipage}[b]{0.48\textwidth}
		\includegraphics[width = 0.65\textwidth]{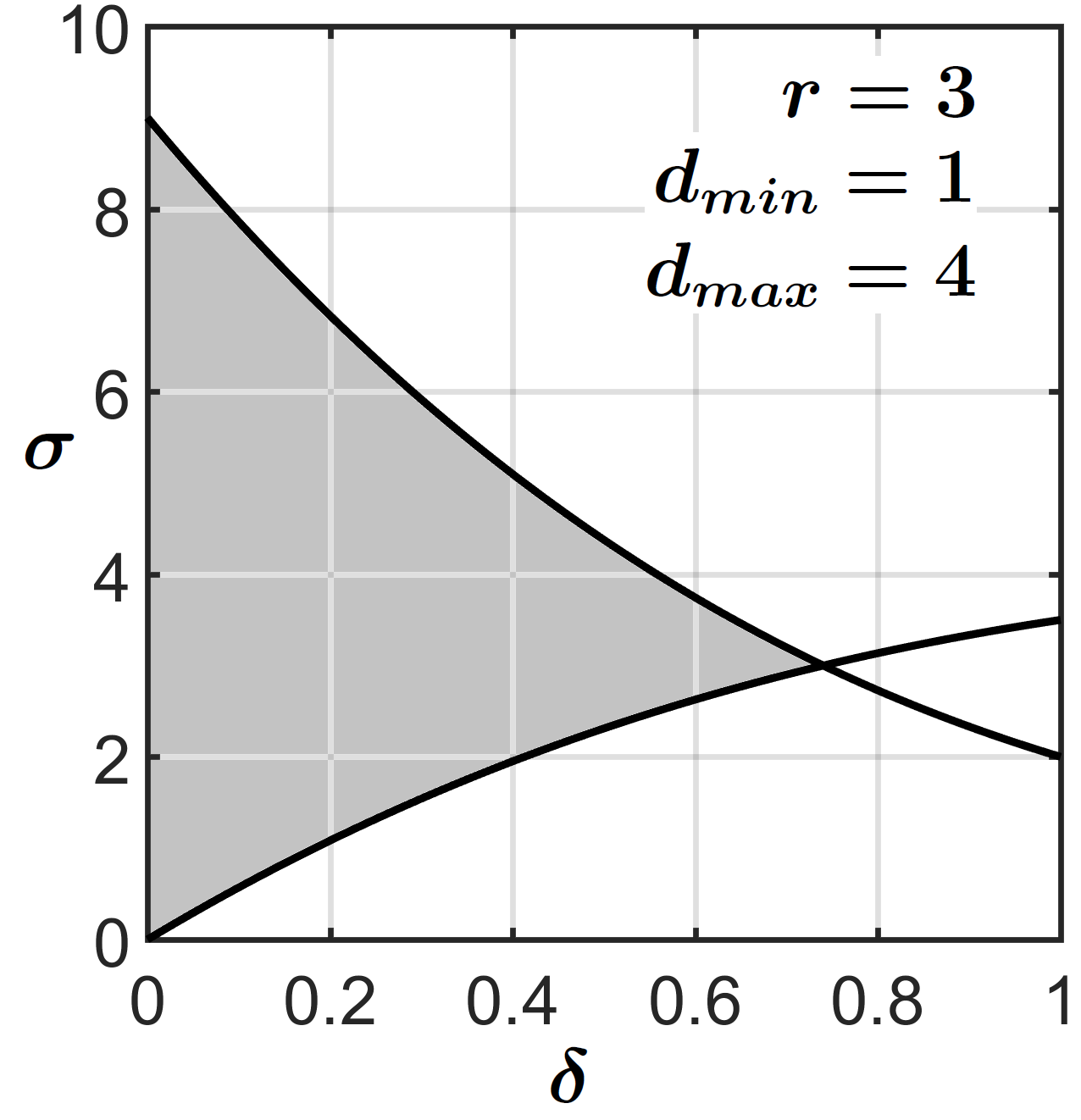}
		\caption{Example of feasible $(\param, \sfact)$-values (shaded region),
		that satisfy (SC) for order $r = 3$ and values $d_\mathrm{min} = 1$,
		$d_\mathrm{max} = 4$. The feasible point that maximizes $\param$ is approximately
		$(\param, \sfact) \approx (0.74, 3)$.}
		\label{Fig:ExampleFeasibleValues}
	\end{minipage}
	  \hfill
	\begin{minipage}[b]{0.48\textwidth}
		\includegraphics[width = 0.65\textwidth]{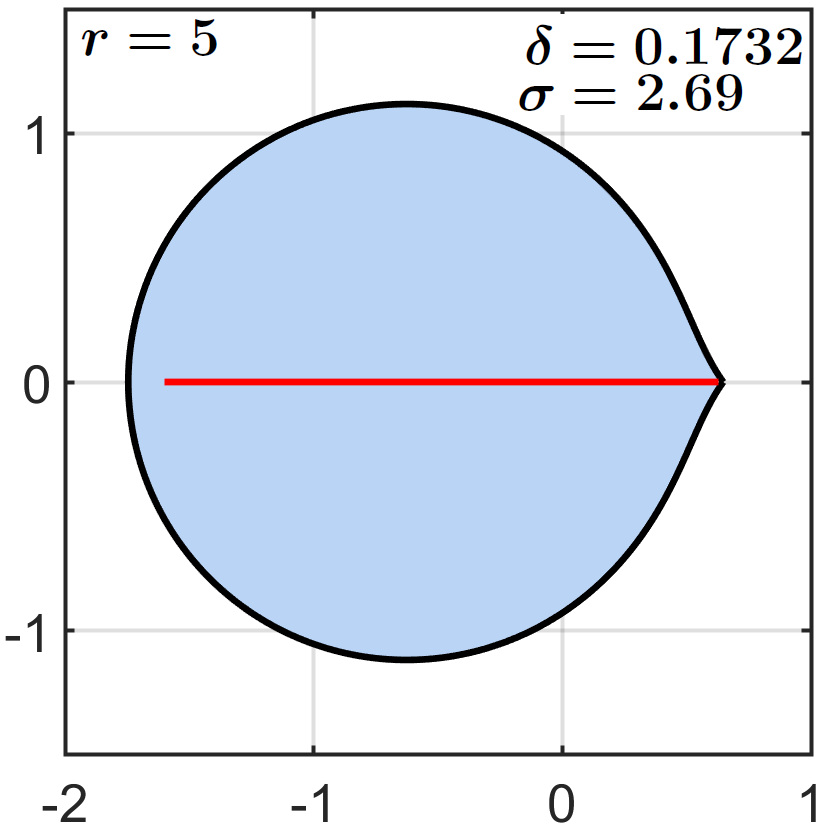}
		\caption{Set $\mathcal{D}$ (blue region) containing
		$W_1(\mat{A}_h, \mat{B}_h)$ (shown in red) with (optimal)
		parameters $(\param, \sfact) = (0.1732, 2.69)$, order $r = 5$, $\sizeN = 64$. \\}
	\label{NumericalRange_DiffEq}
  	\end{minipage}
\end{figure}

We now follow Recipe~\ref{Recipe:ChooseBothImexAndSplit} to
choose both $(\param, \sfact)$ to design an unconditionally stable scheme:
\begin{enumerate}[ 1.]
	\item Retrieve the formulas for $\mathcal{D}$. Since
	both $\Lambda(\mat{A}_h, \mat{B}_h)$ and $W_1(\mat{A}_h, \mat{B}_h)$ are
    real, it is sufficient to use the interval $[m_l, m_r]$ of
	$\mathcal{D}$ on the real line via the formulas (B2).
	\item The second step is heuristic only: establish a range
	of $(\param, \sfact)$-values that ensure the (NC), i.e.\
	$\Lambda(\mat{A}_h, \mat{B}_h) \subseteq \mathcal{D}$. In this case,
	the upper (resp.\ lower) estimate for $\mu_\mathrm{max}$
    (resp.\ $\mu_\mathrm{min}$) agrees exactly with the upper
    (resp.\ lower) estimate on $W_{1}(\mat{A}_h, \mat{B}_h)$. Therefore,
	there is (essentially) no difference in trying to ensure that
	$\Lambda(\mat{A}_h, \mat{B}_h) \subseteq \mathcal{D}$, versus
	$W_1(\mat{A}_h, \mat{B}_h) \subseteq \mathcal{D}$.
	\item Apply the (SC) to determine feasible $(\param, \sfact)$-values.
	Setting $W_1(\mat{A}_h, \mat{B}_h) \subseteq \mathcal{D}$
	requires that the endpoints of $W_1(\mat{A}_h, \mat{B}_h)$ lie
    within $\mathcal{D}$:
	\begin{align}\label{Inequality_Left}
		\text{Left endpoint of } W_1(\mat{A}_h, \mat{B}_h) \text{ in } \mathcal{D}:
        & \quad	
		m_l <  1 - \sfact^{-1} d_\mathrm{max}, \\
        \label{Inequality_Right}
		\text{Right endpoint of } W_1(\mat{A}_h, \mat{B}_h) \text{ in } \mathcal{D}:
        & \quad	
		1 - \sfact^{-1} d_\mathrm{min} < m_r.
	\end{align}
    Equations \eqref{Inequality_Left}--\eqref{Inequality_Right}
    can be rewritten as:
	\begin{align}\label{VarDiff_Inequality}
		(1-m_l)^{-1} d_\mathrm{max} < \sfact
        \quad \text{and} \quad
		\sfact < (1 - m_r)^{-1} d_\mathrm{min}.
	\end{align}	
	The inequalities \eqref{VarDiff_Inequality}, along with
	$\sfact > 0$, $0 < \param \leq 1$, establish the feasible
    points $(\param, \sfact)$ that guarantee unconditional stability.
    This feasible set is always non-empty, because as $\param \rightarrow 0$,
    \eqref{VarDiff_Inequality} yields
    $0 < \sfact < (1 + \cos^{-r}(\pi/r)) d_\mathrm{min}$ for $2 \leq r \leq 5$.
    Hence, one can always achieve unconditional stability,
    by choosing $\param$ small enough.
	\item The last step is to choose a value $(\param, \sfact)$ in the feasible set
	that \emph{maximizes} $\param$ ---
    which is a proxy for minimizing the numerical truncation error.\\[.2em]
	\textbf{Case 1:} $1 \leq r \leq 2$.
    Here the maximum value of $\param^* = 1$ is feasible
    (i.e.,\ one may use SBDF).
    The upper bound inequality \eqref{VarDiff_Inequality}
	is satisfied since $m_r = 1$ yields $\sfact < \infty$.
	The lower bound constraint on $\sfact$ leads to a range of possible
	$(\param^*, \sfact^*)$ values:
	\begin{align}\label{Eq:OptiomalValuesCase1}
		\text{For $r = 1$: }\; \param^* = 1, \ \
        \sfact^* > \tfrac{1}{2} d_\mathrm{max},
        \quad\quad
		\text{For $r = 2$: }\; \param^* = 1, \ \
        \sfact^* > \tfrac{3}{4} d_\mathrm{max}.
	\end{align}
	Generally speaking, choosing $\sfact^*$ large leads to
	large truncation errors. This motivates a choice of $\sfact^*$ close
    to the minimum possible values above.\\[.2em]
	\textbf{Case 2:} $3\leq r \leq 5$. In this case, SBDF may not be
	able to guarantee unconditional stability (see \Srm\ref{SecRecipe})
	when the value $\param^* = 1$ is outside
	the inequalities \eqref{VarDiff_Inequality}. 	
	Figure~\ref{Fig:ExampleFeasibleValues} displays the
    allowable $(\param, \sfact)$-values
	defined by the inequalities in \eqref{VarDiff_Inequality} for some
    representative $d_\mathrm{min}$ and $d_\mathrm{max}$
	values. Note that the optimal point
	(i.e.\ maximum $\param$) occurs at the intersection of the
	inequalities \eqref{VarDiff_Inequality}, and is below $1$.
	To solve for the optimal $(\param^*, \sfact^*)$
	values, we set the upper and lower bounds in \eqref{VarDiff_Inequality}
	\emph{almost} equal to each other. Specifically, introduce a
	\emph{gap parameter} $0 < \eta < 1$, and set the left and right hand
	inequalities for $\sfact$ in \eqref{VarDiff_Inequality} equal to within
	a factor of $(1-\eta)$ to eliminate $\sfact$:
	\begin{align}\label{InEquality:GapEquation}
		\Big(1 - (1-\param/2)^r \Big) d_\mathrm{max} = (1- \eta)
		\Big(1 + (1-\param/2)^r \cos^{-r}(\pi/r) \Big) d_\mathrm{min}.
	\end{align}
	Equation \eqref{InEquality:GapEquation} defines an optimal
    (largest, up to a $(1-\eta)$ error) $\param^*$ value.
    Substituting this $\param^*$ back into \eqref{VarDiff_Inequality}
    yields a range (roughly of size $(1-\eta)$) of feasible $\sfact$ values.
	Among those, we choose (somewhat arbitrarily) $\sfact^*$ as the average
    of the two bounds in \eqref{VarDiff_Inequality}.
	Formulas for the (almost) optimal solutions $(\param^*, \sfact^*)$ are
	given as follows: fix a gap parameter $0 < \eta < 1$
    (smaller values of $\eta$ are more optimal) and order $3 \leq r \leq 5$:
	\begin{align}\label{InequalitySolution}
		&\param^* = 2 - 2 \Big( \frac{1-\kappa}{1+\kappa \cos^{-r}(\pi/r) }
		\Big)^{1/r},	\quad
		\sfact^* = d_\mathrm{min} \Big(1-\frac{1}{2}\eta \Big)
		\frac{1+\cos^{-r}(\pi/r)}{1+\kappa \cos^{-r}(\pi/r)}, \quad
		\\ \nonumber
		&\text{where } \kappa = \frac{d_\mathrm{min}}{d_\mathrm{max}}(1-\eta).
	\end{align}
\end{enumerate}

\begin{remark}\label{Rmk:FailureSBDFDiff}
	(Failure of unconditional stability for SBDF and orders $3 \leq r \leq 5$)
	The necessary conditions for unconditional stability
	require that the generalized eigenvalues
	$\Lambda(\mat{A}_h, \mat{B}_h) \subseteq \mathcal{D}$. The formulas
	from Proposition~\ref{Prop_NumR_VarDiff} lead to the requirement that:
	\[
		(1 - 2^{-r})  \; d_{2,\mathrm{max}} \leq \sfact
        \quad \text{and} \quad
		\sfact \leq \big( 1 + 2^{-r}\cos^{-r}(\pi/r) \big) \; d_{2,\mathrm{min}}.
	\]
	These two inequalities cannot be simultaneously satisfied if
	$d_{2,\mathrm{max}} > D_r \, d_{2,\mathrm{min}}$, where $D_r$ is given by
	$D_3 = 2.1429$, $D_4 = 1.2667$, $D_5 = 1.0931$ for orders $r =3, 4, 5$, respectively.
	Note that $d_{2,\mathrm{max}}/d_{2,\mathrm{min}}$ is
    (up to $\mathcal{O}(\sizeN^{-1})$) a measure
	of the ratio $d_\mathrm{max}/d_\mathrm{min}$.
    As a result, if the ratio between the maximum and
	minimum diffusion coefficient values exceed $D_r$,
    then SBDF cannot provide unconditional
	stability for splittings of the form \eqref{VarDiffusion_Splitting}.
\end{remark}

\begin{remark}\label{Rmk:OvercomingSBDFDiff}
	(Overcoming limitations for SBDF and orders $3 \leq r \leq 5$)
	The formulas \eqref{InequalitySolution} provide a way to overcome
	the unconditional stability limitations encountered with SBDF methods for
	variable coefficient diffusion problems --- regardless of the
	diffusion coefficient $d(x)$.
\end{remark}


To demonstrate that the new approach works in practice, we conduct a convergence
test of equation \eqref{Eq:VarDiffPDE} with the variable diffusion coefficient
\begin{align*}
	d(x) &= 4 + 3 \cos(2\pi x) \quad \Longrightarrow
	\quad d_\mathrm{min} = 1, \quad d_\mathrm{max} = 7.
\end{align*}
This ratio $d_\mathrm{max}/d_\mathrm{min} = 7$ exceeds the value that can
be stabilized by SBDF (see Remark~\ref{Rmk:FailureSBDFDiff}).
Using a value of $\eta = 0.1$ in \eqref{InequalitySolution},
and order $r = 5$ yields the \imex and splitting parameters
$(\param, \sfact) = (0.1732, 2.69)$.
Since the unconditional stability region $\mathcal{D}$
becomes smaller as the order $r$ increases, using $r = 5$
automatically guarantees unconditional stability for all orders
$1 \leq r \leq 5$.
We manufacture the forcing $f(x,t)$ to generate an exact solution
\begin{equation}\label{Eq:VarDiffPDE_solution}
	u^*(x, t) = \sin(20 t) e^{\sin(2\pi x)},
\end{equation}
run up to final time $t_f = 5$.
The multistep scheme is initialized with the exact data:
$u_j = u^*(j \delt)$ for $j = 0, -1, \ldots, -r+1$.
The spatial resolution is $\sizeN = 64$.

Table~\ref{VariableDiff_CvgTest3} shows the error $\|u - u^*\|_{\infty, h}$, using the discrete maximum norm $\| u \|_{\infty, h} = \max_{1\leq j \leq \sizeN} |u(x_j)|$, capped at $10^{-9}$. Convergence rates for $1 \leq r \leq 5$ are reported.
Note that with $\sizeN = 64$, the diffusive time step restriction is $\delt \leq 2^{-18}$. Hence, time steps can be used that are orders of magnitude larger than those required by an explicit scheme. This highlights the benefits of unconditional stability when performing computations with progressively smaller grids.


\begin{table}[htb!]
\centering
{\small
\begin{tabular}{|@{~}c@{~}|@{~}c@{~}||@{~}c@{~}| @{~}c@{~}|| @{~} c@{~} |@{~} c@{~} ||@{~} c@{~} | @{~} c@{~} || @{~} c@{~} |@{~} c@{~} || @{~} c@{~} |@{~} c@{~} |}
\hline
Num.    	& 	$\delt$		&	Error & Rate & Error & Rate & Error & Rate & Error & Rate & Error & Rate  \\ \hline
Steps  		& 					&	$r = 1$ &  & $r = 2$ &  & $r = 3$ &  & $r = 4$ &  & $r = 5$ &  \\ \hline
\phantom{0000}5 	&	$1$	 		& 	  7.9e+00 	 & 	 - &	  4.7e+01 	 & 	- &	  3.4e+02 	 & 	 - &	  9.0e+02 	 & 	- &	  8.8e+02 	 & 	- \\
\phantom{000}10	&	$2^{-1}$	& 	  3.4e+00 	 & 	 \phantom{-}1.2 &	  6.7e+01 	 & 	 -0.5 &	  4.9e+02 	 & 	 -0.5 &	  2.2e+03 	 & 	 -1.3 &	  4.3e+03 	 & 	 -2.3 \\
\phantom{000}20 	&	$2^{-2}$	& 	  4.3e+00 	 & 	 -0.4 &	  2.4e+01 	 & 	 \phantom{-}1.5 &	  5.6e+02 	 & 	 -0.2 &	  3.7e+03 	 & 	 -0.7 &	  6.3e+03 	 & 	 -0.5 \\
\phantom{000}40 	&	$2^{-3}$	& 	  1.3e+00 	 & 	 \phantom{-}1.7 &	  3.5e+01 	 & 	 -0.5 &	  5.4e+02 	 & 	 \phantom{-}0.1 &	  6.3e+03 	 & 	 -0.8 &	  5.8e+04 	 & 	 -3.2 \\
\phantom{000}80 	&	$2^{-4}$	& 	  6.9e-01 	 & 	 \phantom{-}1.0 &	  7.1e+00 	 & 	 \phantom{-}2.3 &	  1.3e+01 	 & 	 \phantom{-}5.4 &	  7.4e+02 	 & 	 \phantom{-}3.1 &	  6.0e+03 	 & 	 \phantom{-}3.3 \\
\phantom{00}160 	&	$2^{-5}$	& 	  2.7e-01 	 & 	 \phantom{-}1.4 &	  1.0e+00 	 & 	 \phantom{-}2.8 &	  1.1e+01 	 & 	 \phantom{-}0.2 &	  5.3e+01 	 & 	 \phantom{-}3.8 &	  5.7e+01 	 & 	 \phantom{-}6.7 \\
\phantom{00}320	&	$2^{-6}$	& 	  2.2e-01 	 & 	 \phantom{-}0.3 &	  6.0e-01 	 & 	 \phantom{-}0.8 &	  2.5e+00 	 & 	 \phantom{-}2.2 &	  2.8e+00 	 & 	 \phantom{-}4.2 &	  7.1e+00 	 & 	 \phantom{-}3.0 \\
\phantom{00}640	&	$2^{-7}$	& 	  2.9e-01 	 & 	 -0.4 &	  5.3e-01 	 & 	 0.2 &	  6.3e-01 	 & 	 \phantom{-}2.0 &	  1.5e-01 	 & 	 \phantom{-}4.3 &	  4.0e-01 	 & 	 \phantom{-}4.1 \\
\phantom{0}1280			&	$2^{-8}$	& 	  2.5e-01 	 & 	 \phantom{-}0.2 &	  2.2e-01 	 & 	 \phantom{-}1.3 &	  5.0e-02 	 & 	 \phantom{-}3.7 &	  3.6e-02 	 & 	 \phantom{-}2.1 &	  2.5e-02 	 & 	 \phantom{-}4.0 \\
\phantom{0}2560			&	$2^{-9}$	& 	  1.6e-01 	 & 	 \phantom{-}0.6 &	  5.6e-02 	 & 	 \phantom{-}2.0 &	  4.9e-03 	 & 	 \phantom{-}3.4 &	  3.5e-03 	 & 	 \phantom{-}3.4 &	  2.8e-04 	 & 	 \phantom{-}6.4 \\
\phantom{0}5120			&	$2^{-10}$	& 	  9.1e-02 	 & 	 \phantom{-}0.8 &	  1.2e-02 	 & 	 \phantom{-}2.2 &	  8.5e-04 	 & 	 \phantom{-}2.5 &	  2.0e-04 	 & 	 \phantom{-}4.1 &	  1.0e-05 	 & 	 \phantom{-}4.8 \\
1.0e+04			&	$2^{-11}$	& 	  4.8e-02 	 & 	 \phantom{-}0.9 &	  2.8e-03 	 & 	 \phantom{-}2.1 &	  1.3e-04 	 & 	 \phantom{-}2.7 &	  1.1e-05 	 & 	 \phantom{-}4.2 &	  3.8e-07 	 & 	 \phantom{-}4.7 \\
2.0e+04			&	$2^{-12}$	& 	  2.5e-02 	 & 	 \phantom{-}1.0 &	  6.7e-04 	 & 	 \phantom{-}2.1 &	  1.8e-05 	 & 	 \phantom{-}2.9 &	  6.1e-07 	 & 	 \phantom{-}4.2 &	  1.3e-08 	 & 	 \phantom{-}4.9 \\
4.1e+04			&	$2^{-13}$	& 	  1.2e-02 	 & 	 \phantom{-}1.0 &	  1.6e-04 	 & 	 \phantom{-}2.0 &	  2.4e-06 	 & 	 \phantom{-}2.9 &	  3.6e-08 	 & 	 \phantom{-}4.1 &	  1.1e-09 	 & 	 \phantom{-}3.5 \\
8.2e+04			&	$2^{-14}$	& 	  6.3e-03 	 & 	 \phantom{-}1.0 &	  4.0e-05 	 & 	 \phantom{-}2.0 &	  3.0e-07 	 & 	 \phantom{-}3.0 &	  2.2e-09 	 & 	 \phantom{-}4.0 &	  1.4e-09 	 & 	 -\\
1.6e+05			&	$2^{-15}$	& 	  3.1e-03 	 & 	 \phantom{-}1.0 &	  9.8e-06 	 & 	 \phantom{-}2.0 &	  3.8e-08 	 & 	 \phantom{-}3.0 &	  2.3e-10 	 & 	 \phantom{-}3.3	&	  2.8e-09 	 & 	 - \\
 \hline
 \end{tabular}
 }
	\caption{Errors $\|u - u^*\|_{\infty,h}$ for variable coefficient diffusion test case
	\eqref{VarDiffusion_Splitting} with \eqref{Eq:VarDiffPDE_solution}, using
	\imex and splitting parameters $(\param, \sfact) = (0.1732, 2.69)$, final time
	$t_f = 5$, and $\sizeN = 64$ Fourier modes.
	The range of fourth and fifth order convergence is capped due to round-off errors, 
	amplified by the problem's conditioning.}
	\label{VariableDiff_CvgTest3}
\end{table}


\subsection{A nonlinear example: Diffusion in porous media and anomalous diffusion rates}

Thus far, the unconditional stability theory has been applied exclusively to linear problems. Now we use the linear theory as a guide for choosing $(\param, \sfact)$ in nonlinear problems, and numerically demonstrate that the new concepts work.
In spirit, the presented methodology shares some similarities with Rosenbrock methods (chapter VI.4, \cite{WannerHairer1991}) in that it also avoids nonlinear implicit terms by means of a properly chosen linear implicit term. A key difference --- other than the fact that Rosenbrock methods are multistage schemes --- is that we do not compute Jacobian matrices (which can be dense and time-dependent), but rather always invert a simple constant coefficient matrix determined from the theory. Hence, the new approach offers more flexibility for choosing efficiency-based implicit terms.

We consider a nonlinear model for a gas diffusing into a porous medium
\cite{Leibenzon1945, Muskat1937}:
\begin{align*}
	\rho_t + \nabla \cdot \big(  \vec{V} \rho \big) = 0
	\quad \text{(Conservation of mass)}, \quad \quad\quad
	\vec{V} = -\frac{\tilde{\kappa}}{\tilde{\mu}} \nabla p
	\quad\text{(Darcy's law)}.
\end{align*}
Here $\tilde{\kappa}$ is the intrinsic permeability of the medium, and $\tilde{\mu}$ is
the effective viscosity. Combined with the equation of state
$p = p_0 \, \rho^{\tilde{\gamma}}$, where $\tilde{\gamma}$ is the adiabatic
constant ($\tilde{\gamma} = 5/3$ for an
ideal monatomic gas), the porous media equation takes the form:
\begin{align}\label{Eq:NonlinearDiffusion}
	\rho_t = a \, \nabla \cdot \big( \rho^{\tilde{\gamma}} \, \nabla \rho \big),
	\quad \text{on } \varOmega \times (0, T].
\end{align}
The constant $a = \kappa p_0 \tilde{\gamma} \; \tilde{\mu}^{-1}$ may, without loss of generality, be set to any positive value by re-scaling time.

We discretize \eqref{Eq:NonlinearDiffusion} in three space dimensions using $N^3$
Fourier modes on the periodic domain $\varOmega = [0, 1]^3$.
Our goal is to achieve unconditional stability by choosing the discrete matrix $\mat{A}_h \approx \sfact \nabla^2$ proportional to the constant coefficient Laplacian (which is easy to treat implicitly). This approach then avoids an implicit treatment of nonlinear terms, thereby bypassing the need for nonlinear solvers.
Due to the nonlinearity in the diffusion coefficient, our choice of $(\param, \sfact)$ via the formulas \eqref{InequalitySolution} requires estimates for the maximum and minimum values of the solution $\rho(x,y,z,t)$ over the simulation. At first glance, it may seem troubling to require time stepping parameters based on the solution; however this is not unusual --- numerical simulations for nonlinear PDEs often require choosing a time step $\delt$ that may depend on the solution.

To test the approach for unconditional stability and accuracy, we perform convergence tests of \eqref{Eq:NonlinearDiffusion} with $\tilde{\gamma} = 5/3$ and $a = 1$, using the manufactured solution
\begin{align}\label{Eq:ExactNonlinearSolution}
	\rho^*(x,y,z,t) &= 2e + e^{\sin(4\pi x)} \cos(2 \pi y) \cos(2\pi z) \cos(t).
\end{align}
We initialize the \imex scheme with the exact initial data $\rho_j(x,y,z) = \rho^*(x,y,z, j \delt)$ for $j = 0, -1, \ldots, -r+1$. We estimate the maximum and minimum value of the nonlinear diffusion coefficient:
\[
	\max_{\vec{x}\in\varOmega, t \in \mathbb{R}} \;
	\rho^*(x,y,z,t)^{\tilde{\gamma}} \leq (3e)^{5/3},
	\quad \quad \quad
	\min_{\vec{x}\in\varOmega, t \in \mathbb{R}} \;
	\rho^*(x,y,z,t)^{\tilde{\gamma}} \geq e^{5/3}.
\]
Using formulas \eqref{InequalitySolution} with $d_\mathrm{max} = (3e)^{5/3}$, $d_\mathrm{min} = e^{5/3}$, $\eta = 0.1$, and $r = 5$ (so that the resulting scheme is stable for all orders $1$ through $5$), yields: $(\param, \sfact) = (0.19166, 13.8)$.

Table~\ref{Table:NonlinearDiff_CvgTest1} shows the numerical error $\|\rho - \rho^*\|_{\infty, h}$, using the discrete norm $\| u \|_{\infty, h} = \max_{\vec{x} \in \text{grid}} |u(\vec{x})|$ evaluated at the final time $t_f = 1$, using $64^3$ grid points ($\sizeN = 64$). In addition to confirming the convergence orders, the table demonstrates that the scheme is stable for $\delt$-values far larger than required by a fully explicit scheme. Moreover, we have confirmed the observations using other values $\tilde{\gamma} \neq 5/3$ and other manufactured solutions (not shown here).

\begin{table}[htb!]
\centering
{\small
\begin{tabular}{ |@{~} c@{~}|@{~} c@{~}||@{~} c@{~}| @{~} c@{~}|| @{~} c@{~} |@{~} c@{~}
||@{~} c@{~} | @{~} c@{~} || @{~} c@{~} |@{~} c@{~} || @{~} c@{~} |@{~} c@{~} |}
\hline
 Num. 		 & $\delt$ & Error & Rate & Error & Rate & Error & Rate & Error & Rate & Error & Rate  \\ \hline
 Steps 		 &         & $r = 1$ &  & $r = 2$ &  & $r = 3$ &  & $r = 4$ &  & $r = 5$ &  \\ \hline

\phantom{000}8	&	$2^{-3} $	 & 	  1.0e+00 	 & 	 - &	  8.3e-01 	 & 	 - &	  6.4e-02 	 & 	 - &	  9.7e-02 	 & 	 - &	  3.4e-04 	 & 	 - \\
\phantom{00}16	&	$2^{-4} $	 & 	  7.7e-01 	 & 	 0.4 &	  3.8e-01 	 & 	 1.1 &	  3.4e-02 	 & 	 0.9 &	  2.2e-02 	 & 	 2.1 &	  8.1e-04 	 & 	 -1.2 \\
\phantom{00}32	&	$2^{-5} $	 & 	  5.0e-01 	 & 	 0.6 &	  8.3e-02 	 & 	 2.2 &	  8.6e-03 	 & 	 2.0 &	  1.9e-03 	 & 	 3.6 &	  1.2e-04 	 & 	\phantom{-}2.7 \\
\phantom{00}64	&	$2^{-6} $	 & 	  2.6e-01 	 & 	 0.9 &	  1.5e-02 	 & 	 2.4 &	  1.4e-03 	 & 	 2.6 &	  1.2e-04 	 & 	 4.0 &	  7.6e-06 	 & 	\phantom{-}4.0 \\
\phantom{0}128	&	$2^{-7} $	 & 	  1.3e-01 	 & 	 1.0 &	  3.6e-03 	 & 	 2.1 &	  1.9e-04 	 & 	 2.8 &	  6.6e-06 	 & 	 4.2 &	  3.0e-07 	 & 	\phantom{-}4.7 \\
\phantom{0}256	&	$2^{-8} $	 & 	  6.4e-02 	 & 	 1.0 &	  8.6e-04 	 & 	 2.1 &	  2.5e-05 	 & 	 2.9 &	  3.8e-07 	 & 	 4.1 &	  1.3e-08 	 & 	\phantom{-}4.5 \\
\phantom{0}512	&	$2^{-9} $	 & 	  3.2e-02 	 & 	 1.0 &	  2.1e-04 	 & 	 2.0 &	  3.2e-06 	 & 	 3.0 &	  2.2e-08 	 & 	 4.1 &	  6.2e-09 	 & 	 - \\
1024			&	$2^{-10}$ 	 & 	  1.6e-02 	 & 	 1.0 &	  5.2e-05 	 & 	 2.0 &	  4.0e-07 	 & 	 3.0 &	  9.8e-10 	 & 	 4.5 &	  1.2e-08 	 & 	 - \\
2048			&	$2^{-11}$ 	 & 	  7.8e-03 	 & 	 1.0 &	  1.3e-05 	 & 	 2.0 &	  5.0e-08 	 & 	 3.0 &	  6.9e-10 	 & 	 - &	  2.4e-08 	 & 	 - \\
 \hline
 \end{tabular}
 }
	\caption{Errors $\|\rho - \rho^*\|_{\infty, h}$
	for three-dimensional nonlinear diffusion coefficient test case with
	manufactured solution \eqref{Eq:ExactNonlinearSolution} and
	$(\param, \sfact) = (0.19166, 13.8)$.
	Errors are computed at the final time $t_f = 1$, with
	$64^3$ ($\sizeN = 64$) Fourier modes.
	Cancellation errors, amplified by the problem's conditioning, 
	limit the
	observed range of fourth and fifth order convergence.}
	\label{Table:NonlinearDiff_CvgTest1}
\end{table}

Now we conduct a test in which the nonlinear behavior is natural to equation \eqref{Eq:NonlinearDiffusion}: a decaying/spreading profile without forcing. We use a Gaussian $\rho(x, y, z, 0) = 1 + e^{-\|\vec{x} - (0.5,0.5,0.5) \|^2 /0.15^2}$ as initial data, which is not exactly periodic, but is sufficiently resolved in space using $128^3$ Fourier modes ($N = 128$) to carry out temporal convergence studies. We choose $\tilde{\gamma} = 5/3$, and set $a = 2^{-4}$, which for the given initial data will lead to dynamics that evolve on an $\mathcal{O}(1)$ time scale.
Using \eqref{InequalitySolution} with $d_\mathrm{max} = 2^{\tilde{\gamma}}$, $d_\mathrm{min} = 1$, $\eta = 0.1$, and $r = 3$, we obtain $(\param, \sfact) = (0.794, 2.616)$.
To generate the initial data required to start the high-order multistep methods, we use the low order ($r = 1$ and $r = 2$) unconditionally stable schemes with many subgrid time steps.

Table~\ref{Table:NonlinearDiff_CvgDecayGaussian} provides a verification for the test, by computing the convergence rate estimate
\begin{align}\label{Def:Ratio}
	R_{\delt} := \log_2\big( \| \rho_{4\delt}(t_f) - \rho_{2\delt}(t_f)\|_{\infty,h}
	/\| \rho_{2\delt}(t_f) - \rho_{\delt}(t_f)\|_{\infty,h} \big),
\end{align}
where $\rho_{k}(t_f)$ denotes the discrete solution at $t_f$ computed using time step $k$. Table~\ref{Table:NonlinearDiff_CvgDecayGaussian} confirms that at small $\delt$ values, $R_{\delt}$ converges to the order of the scheme.

\begin{table}[htb!]
\centering
{\small
\begin{tabular}{ | c || c |  c | c | c | c | c | c | c | c |}
\hline
$\delt$ &	$2^{-4}$	 &	$2^{-5}$	&	$2^{-6}$	&	$2^{-7}$	&	$2^{-8}$	&	$2^{-9}$	&	$2^{-10}$	&	$2^{-11}$	& $2^{-12}$ \\ \hline \hline
$r = 1$,	$R_{\delt}$  &	1.22	&	1.58	&	1.52	&	1.28	&	1.14	&	1.07	&	1.03	&	1.02	& 1.01 \\	\hline
$r = 2$,	$R_{\delt}$  &	0.39	&	 5.12			&	2.06	&	2.02	&	1.96	&	1.95	&	1.96	&	1.97	& 1.98 \\	\hline
$r = 3$,	$R_{\delt}$  &	1.38	&	3.54	&	1.63	&	8.52	&	2.83	&	2.74	&	2.81	&	2.88	& 2.93	\\	\hline
 \end{tabular}
 }
	\caption{Convergence test for a decaying Gaussian initial data, reporting the values
	$R_{\delt}$ from \eqref{Def:Ratio} at $t_f = 1$, and $128^3$ grid points.}
	\label{Table:NonlinearDiff_CvgDecayGaussian}
\end{table}

\begin{figure}[htb!]
	\centering
	\includegraphics[width = \textwidth]{%
		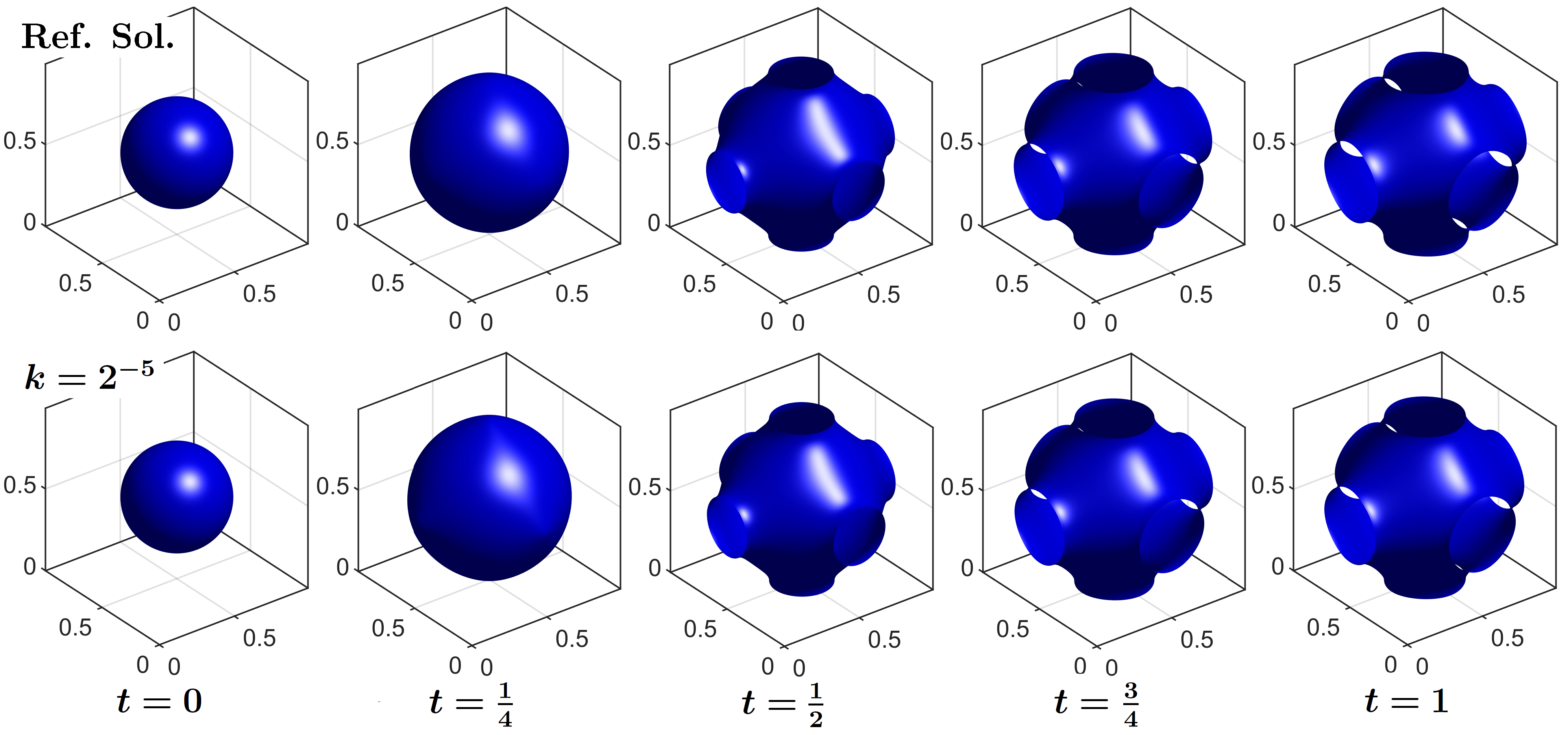} \hfill
		\caption{Evolution of the solution towards $t_f = 1$ for a Gaussian
		initial data. Snapshots are of the level sets
		$\{\vec{x} : \rho(\vec{x},t) = \bar{\rho}\}$ for
		the reference solution (top), compared to
		the solution computed using order $r = 3$ with $(\param, \sfact) = (0.794, 2.616)$
		and a large time step $\delt = 2^{-5}$ (bottom).
		Here $\bar{\rho}$ is the discrete average value of the solution
		(which is a conserved quantity).} 
	\label{Fig:DecayGaussian_Snapshots}
\end{figure}

Finally, guided by the data in Table~\ref{Table:NonlinearDiff_CvgDecayGaussian}, we choose a time step $\delt = 2^{-5}$, and compute the solution towards $t_f =1$ (i.e.~merely 32 time steps). Figure~\ref{Fig:DecayGaussian_Snapshots} visualizes level sets of the solution for different times and compares them to a reference solution (obtained using the fully explicit, second-order scheme with $\delta = 1$, i.e.\ Adams-Bashforth, with $\delt = 2^{-16}$, i.e.\ 65536 time steps). The unconditionally stable method is successfully capturing the solution, albeit using very large time steps.

\begin{figure}[htb!]
	\centering
	\includegraphics[width = 0.65\textwidth]{%
		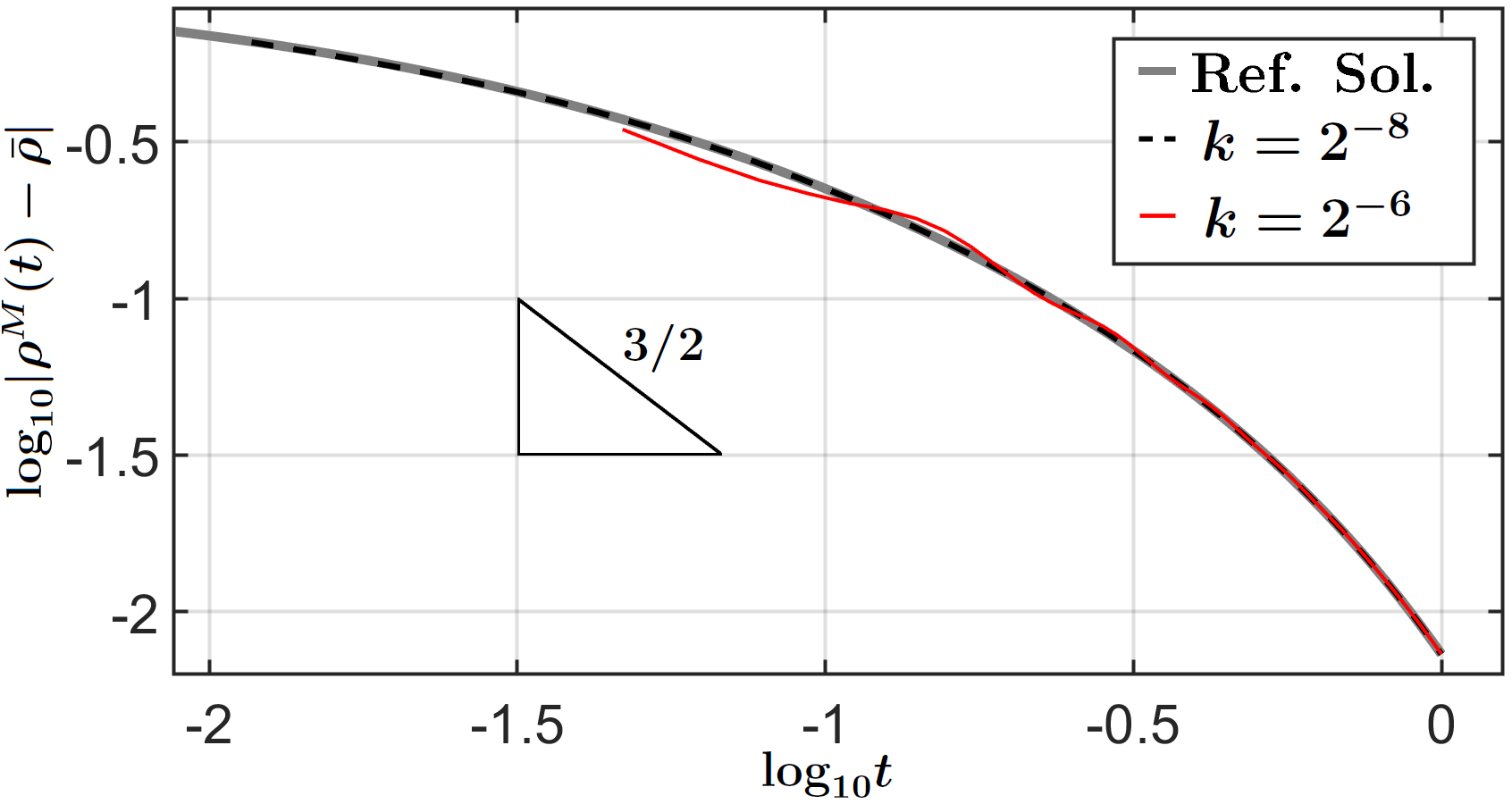} \hfill
		\caption{Decay of the maximum solution
		$\rho^M(t)$ to a constant. For reference,
		the slope of $-3/2$ is shown --- which is the decay rate for
		Gaussian initial data obtained with a constant coefficient diffusion
		equation. The simulation uses the order
		$r = 3$ scheme with $(\param, \sfact) = (0.794, 2.616)$,
		and $\delt = 2^{-6}$ (64 time steps) as well as $\delt = 2^{-8}$ (256 time steps).
        Both are much larger than the
		restriction $\delt = 2^{-16}$ (65536 time steps) required by the
		fully explicit reference solution.} 
	\label{Fig:DecayGaussian_Peak}
\end{figure}

Meanwhile, Figure~\ref{Fig:DecayGaussian_Peak} highlights the nonlinear effect (anomalous diffusion) in the solution with plots of the decay rate of the peak $\rho^M(t) = \|\rho(\vec{x}, t)\|_{\infty, h}$ relative to the mean value $\bar{\rho} := h^3 \sum_{\vec{x} \in \text{grid}} \rho(\vec{x},0) \approx \int_{0}^1\int_0^1\int_0^1 \rho(\vec{x},0) \du \vec{x}$. For reference, the decay rate for a linear constant coefficient diffusion problem ($-3/2$) is shown as well. The plot shows, against the reference solution, the \imex schemes with $\delt = 2^{-6}$ (64 time steps, a small error) and $\delt = 2^{-8}$ (256 time steps, visually indistinguishable).

\section{Incompressible channel flow}\label{Sec:NavierStokes}

In this section we perform a study that zooms in on the
question of unconditional stability for \imex splittings that arise in fluid
dynamics problems. Specifically, for the time-dependent Stokes equation
in a channel geometry, we devise unconditionally stable \imex schemes that treat the pressure explicitly and viscosity implicitly.  High-order \imex schemes that are provably unconditionally stable for the incompressible Navier-Stokes equations are notoriously difficult to attain (see for instance \cite{GuermondMinev2015, LiuLiuPego2007, LiuLiuPego2010}).  This study highlights some peculiar challenges that arise with \imex schemes for incompressible flows, for instance that unconditional stability may depend on model parameters such as the shape and size of the domain.
The new unconditional stability theory may provide new ways to stabilize operator splittings in fluid dynamics applications that might otherwise be unstable.

It is worth mentioning that in fluid flow, the alternative to unconditional stability may not necessarily be detrimental (in fact, the situation here is of that type). For instance, when unconditional stability is not attained, an \imex approach may still provide competitive stability benefits, by incurring a time step restriction that is $\mathcal{O}(1)$ (i.e.\ independent on the spatial mesh), and thus not stiff. This type of stability restriction has been recently referred to as \emph{quasiunconditional stability} \cite{BrunoCubillos2016, BrunoCubillos2017} (and has been applied to the compressible Navier-Stokes equations in an alternate direction implicit (ADI) setting). Nonetheless the question of whether unconditionally stable schemes can be devised is of interest, as in other situations the lack of unconditional stability may not be as forgiving as the quasiunconditional stability scenario.

We focus on reformulations of the incompressible Navier-Stokes equations \cite{Henshaw1994, JohnstonLiu2002, JohnstonLiu2004, LiuLiuPego2010, ShirokoffRosales2010}, that take the form of a pressure Poisson equation (PPE) system (sometimes also referred to as an \emph{extended Navier-Stokes} system). They replace the divergence-free constraint by a non-local pressure operator (defined via the solution of a Poisson equation), and thus allow for the application of time-stepping schemes without having to worry about constraints. In contrast to projection methods \cite{JohnstonLiu2004, LiuLiuPego2010,Henshaw1994, ShirokoffRosales2010}, PPE systems are not based on fractional steps, and thus allow in principle for arbitrary order in time. In turn, the use of \imex schemes allows for an explicit treatment of the pressure --- which (in contrast to fully implicit time-stepping) avoids large saddle-point problems in which velocity and pressure are coupled together. The challenge is that the explicit pressure term may
become stiff (because it can be recast as a function of the viscosity term
\cite{JohnstonLiu2004, ShirokoffRosales2010}).

For simplicity, we restrict this presentation to the linear Navier-Stokes equations (i.e.\ without the advection terms), because these equations already capture the key challenges arising from the interaction between viscosity and pressure. One could also investigate (unconditional) stability for incompressible flows with advection terms; however we do not pursue this here. For a two dimensional domain $\varOmega \subset \mathbb{R}^2$, we use the PPE reformulation by Johnston and Liu \cite{JohnstonLiu2004, JohnstonLiu2002} for problems with no-slip boundary conditions:

\begin{equation} \label{JL_PPE1}
   \left.
   \begin{array}{rcll}
     u_t      & = & u_{xx} + u_{yy} - p_x + f_1
        \quad & \textrm{for}\;\; \mathbf{x} \in \phantom{\partial}\/\varOmega\/,
        \\ \rule{0ex}{2.5ex}
     v_t & = & v_{xx} + v_{yy} - p_y + f_2
        \quad & \textrm{for}\;\; \mathbf{x} \in \phantom{\partial}\/\varOmega\/,
        \\ \rule{0ex}{2.5ex}
       &  & u = v = 0
        \quad & \textrm{for}\;\; \mathbf{x} \in \partial\/\varOmega\/,
   \end{array}
   \right\}
\end{equation}
where $p$ is the solution of
\begin{equation} \label{JL_PPE2}
   \left.
   \begin{array}{rcll}
     p_{xx} + p_{yy}      & = & (f_1)_x + (f_2)_y
        \quad & \textrm{for}\;\; \mathbf{x} \in \phantom{\partial}\/\varOmega\/,
        \\ \rule{0ex}{2.5ex}
      \vec n \cdot \nabla p & = &  -\vec n \cdot (\nabla \times \nabla \times \vec u ) + \vec n \cdot \vec f
        \quad & \textrm{for}\;\; \mathbf{x} \in \partial \/\varOmega\/, \\ \rule{0ex}{2.5ex}
        \int_{\Omega} p(\vec{x}) \du \vec{x} & = & 0
   \end{array}
   \right\}.
\end{equation}
Equation \eqref{JL_PPE1} is the standard momentum equation for the velocity field $(u, v)$, while equation \eqref{JL_PPE2} is the PPE reformulation for the pressure, acting to keep the flow incompressible. Note that the last requirement in equation \eqref{JL_PPE2} is added to uniquely define the pressure. Here, $\vec n$ is the outward facing normal on $\partial \varOmega$, and $\vec{f} = (f_1, f_2)$ is the body force. The viscosity terms $u_{xx} + u_{yy}$ and $v_{xx} + v_{yy}$ are the stiff terms that are treated implicitly (matrix $\mat{A}$), while the pressure terms  $p_x, p_y$ are linear functions of the velocity $(u, v)$ that are treated explicitly (matrix $\mat{B}$). We consider a channel geometry, $\varOmega = [0, L_x]\times(0, 1)$, that is periodic in the $x$-direction.

The simple geometry allows us to solve for the pressure $p$ analytically, and convert equations \eqref{JL_PPE1}--\eqref{JL_PPE2} into a non-local PDE for the velocity $u$. This simplifies the computation of the set $W_p(\mat{A}, \mat{B})$ and provides fundamental insight into why existing \imex splittings that treat the viscosity terms implicitly and pressure terms explicitly may become unstable. In more general problems, one would of course need to conduct a full spatial discretization of equations \eqref{JL_PPE1}--\eqref{JL_PPE2}, and apply the recipes described in \Srm\ref{SecRecipe} to the resulting large ODE system. It is important to note that, while theoretical insights are less clear in that situation, there is no fundamental problem with applying the methodology.

To derive the non-local PDE for the velocity $u$, we start off by setting $\vec{f} = 0$ (unconditional stability does not depend on the forcing). Because $\varOmega$ is periodic in the $x$-direction, we conduct a Fourier expansion in the $x$-direction and set~$(u, v) = ( \imath u(y,t; \xi), v(y,t; \xi) ) e^{\imath \xi x}$
and~$p = p(y; \xi)\/ e^{\imath \xi x}$. The system \eqref{JL_PPE1}--\eqref{JL_PPE2}
then becomes:
\begin{align} \label{JL_SemiAnalytic_1}
	\begin{pmatrix}
		u \\
		v
	\end{pmatrix}_t \!= \Big(\frac{\partial^2}{\partial y^2} - \xi^2 \Big)
	\begin{pmatrix}
		u \\
		v
	\end{pmatrix} -
	\begin{pmatrix}
		\xi p \\
		p_y
	\end{pmatrix} \ \ \text{on~} 0 < y < 1\;, \quad
	&\text{and~}
	\begin{pmatrix}
		u \\
		v
	\end{pmatrix} = \vec{0} \quad \text{on} \ \
	y = \{0, 1\}, \\
    \label{JL_SemiAnalytic_2}
	p_{yy} - \xi^2 p = 0 \ \ \text{on~} 0 < y < 1\;, \quad
    &\text{and~} \
	\frac{dp}{dy} = \xi u_y \ \ \text{on~} y = \{0, 1\}.
\end{align}
The allowable wave numbers are given by $\xi = \pm 2\pi L_x^{-1}\/ n_{\xi}$ and natural numbers $n_{\xi} \in \mathbb{N}$. The pressure equation \eqref{JL_SemiAnalytic_2} is uniquely solvable for all $\xi \neq 0$; while for $\xi = 0$, the integral constraint in \eqref{JL_PPE2}, together with \eqref{JL_SemiAnalytic_2}, fixes $p(y; 0) = 0$. For $\xi \neq 0$, equation \eqref{JL_SemiAnalytic_2} can be solved analytically to obtain:
\begin{align} \label{PressureSolved}
	&p(y; \xi) = -\frac{ \cosh\big(\xi(y - 1)\big)}{\sinh(\xi)} \; u_y(0) +
	\frac{ \cosh(\xi y)}{\sinh(\xi)} \;  u_y(1). \quad
	\xi \neq 0.
\end{align}
The pressure can then be substituted back into equation \eqref{JL_SemiAnalytic_1}
to yield a non-local PDE for the horizontal velocity $u = u(y, t; \xi)$ (for $\xi \neq 0$):
\begin{equation}
\label{SingleComponentChannelFlow}
\begin{split}
	u_t &= \Big( \frac{\partial}{\partial y^2} -\xi^2\Big) u + \xi
	\, \frac{ \cosh\big(\xi(y - 1)\big)}{\sinh(\xi)} \; u_y(0) -
	\xi \, \frac{ \cosh(\xi y)}{\sinh(\xi)} \;  u_y(1)	\\
	&\textrm{Boundary conditions: }\quad u = 0, \textrm{ on } y = \{0, 1\}.
\end{split}
\end{equation}
Solving equation \eqref{SingleComponentChannelFlow} for $u(y, t; \xi)$ at every wave number $\xi$ then allows one to reconstruct $u(x, y, t)$.  In a similar fashion, one can use \eqref{PressureSolved} to reconstruct $p(x, y)$. Once either $u(x, y, t)$ or $p(x,y)$ is known, the vertical velocity $v(x,y,t)$ can then be obtained by solving either (i) the $v$-component of equation \eqref{JL_SemiAnalytic_1} with the pressure as a prescribed forcing, or (ii) using the fact that the PPE reformulation automatically enforces the divergence constraint so that $v(y;\xi)_y = -\xi\/ u(y;\xi)$. Collectively, the solutions $(u,v,p)$ from \eqref{PressureSolved}--\eqref{SingleComponentChannelFlow} solve the original PDEs \eqref{JL_PPE1}--\eqref{JL_PPE2} with $\vec{f} = 0$. Thus, devising unconditionally stable \imex splittings for \eqref{SingleComponentChannelFlow} can be used as a guide for stabilizing discretizations of the full equations \eqref{JL_PPE1}--\eqref{JL_PPE2}.

\subsection{Numerical discretization and \imex splitting for equation \eqref{SingleComponentChannelFlow}}

In line with prior examples, we seek \imex splittings of equation \eqref{SingleComponentChannelFlow} in which the pressure terms are treated explicitly, while a portion of the viscosity is treated implicitly:
\begin{align}\label{PPE_MatrixSplitting}
	\mat A_{h,\xi}\, \vec u \approx
	\sfact \/ \Big( \frac{\partial}{\partial y^2} -\xi^2 \Big) u,
	\quad
	\mat B_{h,\xi} \, \vec u \approx (1-\sfact)\/
	\Big( \frac{\partial}{\partial y^2} -\xi^2 \Big) u - \xi p\;.
\end{align}
Although equation \eqref{SingleComponentChannelFlow} is a non-local
PDE, the highest derivative degree is $2$ so that the splitting
\eqref{PPE_MatrixSplitting} still adheres to the guidelines in
Remark~\ref{Rmk:GuidelinesA0}. As with prior discussions, the
inclusion of the splitting parameter $\sfact$ in \eqref{PPE_MatrixSplitting}
provides additional flexibility in devising stable schemes,
compared to many existing splittings that effectively fix $\sfact = 1$.

To obtain the matrices $(\mat{A}_{h,\xi}, \mat{B}_{h,\xi})$ we discretize
the $y$-direction using $\sizeN_y$ equispaced grid points $y_j = j h$, and
spacing $h = (\sizeN_y + 1)^{-1}$, so that
$\vec{u}_j = u(y_j) \in \mathbb{R}^{\sizeN_y}$, for $1 \leq j \leq \sizeN_y$.
A standard 3-point finite difference stencil for $\partial_{yy}$ (with
Dirichlet boundary conditions $u(0) = u(1) = 0$ at $y \in \{0, 1\}$) leads to
the following discretization: $\mat{A}_{h,\xi} = \sfact \/ \mat{A}_{0, h, \xi}$, where
\begin{align}\label{JL_Mat_A}
	\mat{A}_{0, h, \xi} = \frac{1}{h^2} \begin{pmatrix}
		 -2 & \phantom{-}1 &  &  &  &  \\
		 \phantom{-}1& -2 & \phantom{-}1 &  &  & \\ 		
		  & \ddots & \ddots & \ddots \\
		  &  & \phantom{-}1 & -2 & \phantom{-}1 \\
  		  &  &   & \phantom{-}1 & -2 \\
	\end{pmatrix} - \xi^2 \mat I.
\end{align}	
In a similar fashion, we set
$\mat{B}_{h, \xi} = (1 - \sfact) \mat{A}_{0, h, \xi} + \mat{Q}_{h, \xi}$
where $\mat{Q}_{h, \xi} \vec{u} \approx -\xi p$ is a matrix that computes the pressure.
The matrix $\mat{Q}_{h, \xi}$ is built using equation \eqref{PressureSolved} as
\begin{align}\label{PressureMatrix}
	\mat {Q}_{h, \xi} &= \vec{a}(\xi) \, \vec{d}_1^T + \vec{b}(\xi) \, \vec{d}_2^T
	\quad \text{for} \quad
	\xi \neq 0,
	\quad \text{and} \quad
	\mat{Q}_{h, 0} = \mat{0}
	\quad \text{for} \quad \xi = 0.
\end{align}
Here the vectors
$\vec{d}_{1} = h^{-1}\begin{pmatrix} 1, & 0, & \ldots, & 0 \end{pmatrix}^T$ and
$\vec{d}_{2} = h^{-1}\begin{pmatrix} 0, & \ldots , & 0, -1 \end{pmatrix}^T$
approximate the derivatives $\vec{d}_{1}^T \/ \vec{u} \approx u_y(0)$ and $\vec{d}_{2}^T \/ \vec{u} \approx u_y(1)$, and thus encode the boundary conditions $u(0) = u(1) = 0$.
The vectors $\mat{a}(\xi), \mat{b}(\xi) \in \mathbb{R}^{\sizeN_y}$ are discretizations of the functions
\[
	\vec a(\xi)_j = \xi \; \textrm{csch}(\xi) \cosh\big(\xi(y_j - 1)\big),  \quad	
	\vec b(\xi)_j = -\xi \; \textrm{csch}(\xi)\cosh(\xi y_j).
\]
Note that for each fixed value of $\xi$, the matrix $\mat{Q}_{h, \xi}$ is the sum of two rank-$1$ matrices.

\subsection{Applying the unconditional stability theory to determine $(\param, \sfact)$}

Now we follow the guidelines in Remark~\ref{Rmk:GuidelinesA0} to determine $(\param, \sfact)$. We directly focus on computing an appropriate set $W_p(\mat{A}_{h, \xi}, \mat{B}_{h, \xi})$ to be used in the stability theory.

Unlike prior examples, for the channel flow application considered here, the choice $p=2$ is most useful for the analysis of the $W_p$ sets. This is because the maximum size and shape of the sets $W_2(\mat{A}_{h, \xi}, \mat{B}_{h, \xi})$ are effectively independent of $h$. In contrast, numerical experiments show that the sets $W_1(\mat{A}_{h, \xi}, \mat{B}_{h, \xi})$ arising from \eqref{PPE_MatrixSplitting} tend to grow as $h\rightarrow 0$. It is worth noting that $p = 2$ is also motivated by \cite{JohnstonLiu2004}, in which a version of the set $W_2(\mat{A}_{h, \xi}, \mat{B}_{h, \xi})$ was studied to prove that SBDF1 is unconditionally stable when applied to the channel flow PDEs \eqref{JL_PPE1}--\eqref{JL_PPE2}. Here, we apply the full new unconditional stability theory to systematically investigate stability for high order schemes.

To compute $W_2(\mat{A}_{h, \xi}, \mat{B}_{h,\xi})$,
we use the scaling property from Remark~\ref{Rmk:DependenceOnSigma} to write
\begin{align}\label{NumRangeNS}
	W_2( \sfact \mat{A}_{0, h, \xi}, \mat{B}_{h, \xi} )
	= 1 - \sfact^{-1} + \sfact^{-1} \, W_2( \mat{A}_{0, h, \xi}, \mat{Q}_{h, \xi}),
\end{align}
and use \textsf{Chebfun}'s numerical range (field of values)
routine \cite{DriscollHaleTrefethen2014} to compute
$W_2( \mat{A}_{0, h, \xi}, \mat{Q}_{h, \xi})$; from which
$W_2( \sfact \mat{A}_{0, h, \xi}, \mat{B}_{h, \xi} )$ is obtained via a shift
and re-scaling. Note that \textsf{Chebfun} employs an
algorithm due to Johnson \cite{Johnson1978} that reduces the computation of
$W_2( \sfact \mat{A}_{0, h, \xi}, \mat{B}_{h, \xi} )$ to a collection of
eigenvalue computations.

A theoretical study of the set $W_2(\mat{A}_{0, h, \xi}, \mat{Q}_{h, \xi})$
where continuum operators $\mathcal{A} = (-\Delta \vec{u})$ and $\mathcal{B} = (-\nabla p)$ were used
instead of discrete matrices $\mat{A}_{0, h, \xi}, \mat{Q}_{h, \xi}$ (note that the set
$W_2$ is still defined using operators), was carried out
in part of the work \cite{JohnstonLiu2004}.  They showed that the set
$W_2(\mathcal{A}, \mathcal{B})$, using continuum operators, was real and contained
in the interval $[0, 1]$. Our numerical computations of
$W_2( \sfact \mat{A}_{0, h, \xi}, \mat{B}_{h, \xi} )$
also show that for each fixed value of $\xi$, the sets are within a discretization error
(at most $\mathcal{O}(h)$) of the interval $[0, 1]$.

\begin{figure}[htb]
\centering
\includegraphics[width = \textwidth]{%
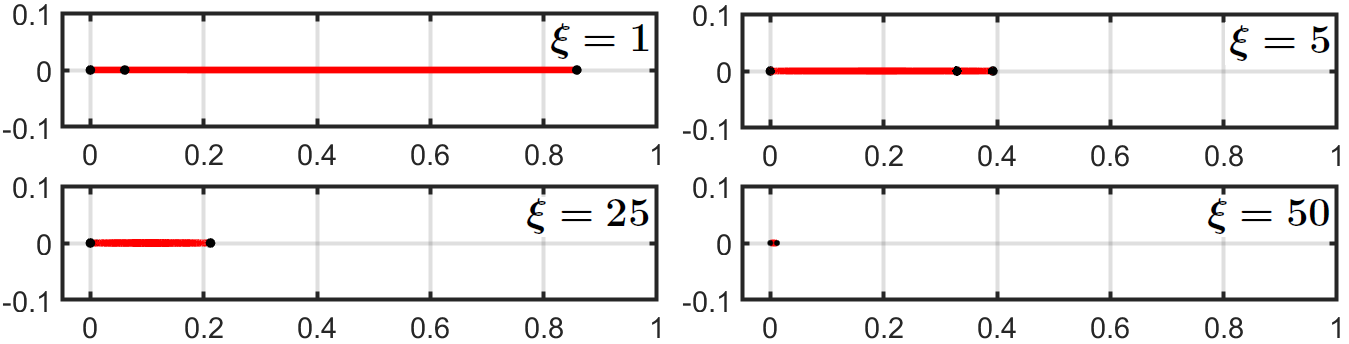}	
\caption{Sets
	$W_2( \mat{A}_{0, h, \xi}, \mat{Q}_{h, \xi})$ (red) in \eqref{NumRangeNS} for
	$\sizeN_y = 32$ and wave numbers $\xi \in \{1, 5, 25, 50\}$.  The sets are computed
	numerically and are confined to the real axis. The sets shrink in size
	as $\xi$ increases. }
\label{NumericalRange_W2_PPE_sets}
\end{figure}

\begin{figure}[htb]
\centering
\includegraphics[width = 0.8\textwidth]{%
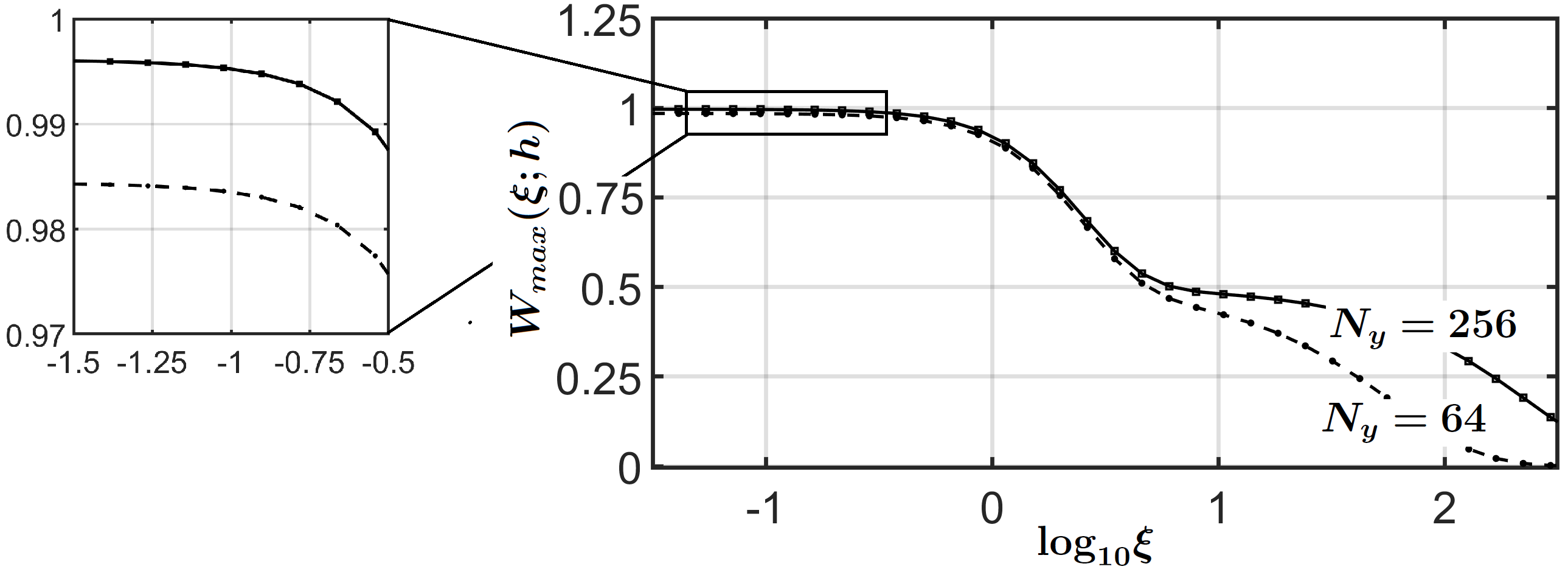}	
\caption{Plot of $W_{max}(\xi; h)$ versus wave
	number $\xi$ for $\sizeN_y = 64$ (dashed)
	and $\sizeN_y = 256$ (solid). The sets
	$W_2( \mat{A}_{0, h, \xi}, \mat{Q}_{h, \xi})$ (and consequently
	$W_2(\sfact \mat{A}_{0, h,\xi}, \mat{B}_{h,\xi})$)
	are contained in the interval $[0, 1]$ along the real axis. }	
\label{NumericalRange_PPE_maxW2}
\end{figure}

To quantify the region that
$W_2(\mat{A}_{0, h, \xi}, \mat{Q}_{h, \xi})$ occupies along the real axis, let
\[
	W_{max}(\xi; h) = \max \textrm{Re}\Big(\/ W_2(\mat{A}_{0, h, \xi}, \mat{Q}_{h, \xi}) \Big),
	\quad
	W_{min}(\xi; h) = \min \textrm{Re}\Big(\/ W_2(\mat{A}_{0, h, \xi}, \mat{Q}_{h, \xi}) \Big).
\]
Figure~\ref{NumericalRange_W2_PPE_sets} plots the sets
$W_2(\mat{A}_{0, h, \xi}, \mat{Q}_{h, \xi})$ for $\xi \in \{1, 5, 25, 50\}$ and grid spacing
$\sizeN_y = 32$.
Figure~\ref{NumericalRange_PPE_maxW2} plots $W_{max}(\xi; h)$ for different
wave numbers $\xi$
(where $\xi = 0$ is a special case and excluded from the plot) and grids
$\sizeN_y \in \{64, 256\}$. The plot shows that $W_{max}(\xi; h)$ decreases
(monotonically) with increasing $\xi$ --- which is important
for the simultaneous stabilization of all wave numbers that arise
in a channel geometry.
The minimum value $W_{min}(\xi; h)$ is always zero, i.e.\ $W_{min}(\xi; h) = 0$.

Since the sets $W_2( \sfact \mat{A}_{0, h, \xi}, \mat{B}_{h, \xi} )$ lie along the real axis,
the procedure for choosing $(\param^*, \sfact^*)$ parallels that of the diffusion example in \Srm\ref{subsection_vardiff}.
For instance, we can use equations \eqref{Eq:OptiomalValuesCase1} (for $r \in \{1, 2\}$) and \eqref{InequalitySolution} (for $r \in \{3, 4, 5\}$) via the substitutions $d_{\max} \rightarrow (1 - W_{min}(\xi; h) ) = 1$ and $d_{\min} \rightarrow (1 - W_{max}(\xi;h))$ to ensure that $W_2( \sfact \mat{A}_{0, h, \xi}, \mat{B}_{h, \xi} ) \subseteq \mathcal{D}$.
This approach determines the parameters $(\param^*, \sfact^*)$ that ensure unconditional stability for the PDE \eqref{PressureMatrix} for one fixed wave number $\xi$. In practice, however, stabilizing equation \eqref{SingleComponentChannelFlow} for a channel geometry requires that it be unconditionally stable for all allowable wave numbers $\xi = \pm 2\pi \/ L_x^{-1} n_x$. The following remark shows that it suffices to ensure unconditional stability for the smallest non-zero wave number (which then stabilizes all others).

\begin{remark}\label{Rmk:StabilizeSmallestXi} (Choosing one $(\param^*, \sfact^*)$ that works for all wave numbers $\xi$)
Unconditional stability for equations \eqref{JL_PPE1}--\eqref{JL_PPE2} corresponds to ensuring that \eqref{SingleComponentChannelFlow} is stable for all modes $\xi = \pm 2\pi \/L_x^{-1} n_{\xi}$. Here we argue why it suffices to ensure that only the smallest mode is stable, i.e., to require that $W_2(\sfact \mat{A}_{0,h, \xi_1}, \mat{B}_{h,\xi_1}) \subseteq \mathcal{D}$, where $\xi_1 = 2\pi L_x^{-1}$ is the smallest positive wave number.

To show this, we use the simple property that the sets $W_2(\cdot, \cdot)$
may be written as a numerical range.	
	First, one can view the simultaneous solution of equation
	\eqref{SingleComponentChannelFlow} over all allowable $\xi$,
	as solving one very large
	system of equations with matrices $\mat{A}_{0,h}$ and $\mat{B}_{h}$
	that are written as direct sums.  Specifically, write
	$\mat{A}_{0,h} = \bigoplus_{\xi} \mat{A}_{0, h, \xi}$ and
	$\mat{B}_{h} = \bigoplus_{\xi} \mat{B}_{h, \xi}$, where the direct
	sum is over the wave numbers $\xi = \pm 2\pi L_x^{-1} n_x$
	and natural numbers $n_{\xi} \in \mathbb{N}$ (i.e.\ $(\mat{A}_{0,h}, \mat{B}_h)$
	are infinite block
	diagonal matrices with each block being $\mat{A}_{0, h, \xi}$ or $\mat{B}_{h, \xi}$).
	Now use the fact that the set
	$W_2( \sfact \mat{A}_{0, h, \xi}, \mat{B}_{h, \xi} )$ may be written
	as a numerical range (see Remark~\ref{Rmk:NumRange}) --- and that
	the numerical range of the direct sum of two matrices is the
	convex hull of their two numerical ranges, i.e.\
	$W(\mat{X} \bigoplus \mat{Y} ) = \mathrm{conv}\{ W(\mat{X}), W(\mat{Y}) \}$.
	As a result, the set $W_2(\sfact \mat{A}_{0,h}, \mat{B}_{h} )$ is the convex hull
	of the sets $W_2(\sfact \mat{A}_{0,h, \xi}, \mat{B}_{h, \xi} )$ over all
	allowable wave numbers.
	Lastly, we observe that the set $W_2(\sfact \mat{A}_{0,h, \xi_1}, \mat{B}_{h, \xi_1} )$
	contains (up to at most an error $\mathcal{O}(h)$)
	each of the sets $W_2(\sfact \mat{A}_{0,h, \xi}, \mat{B}_{h, \xi} )$
	for all $\xi = \pm 2\pi \/ L_x^{-1} \/n_x$ and $n_x \in \mathbb{N}$.
	This shows that the convex hull of the sets
	$W_2(\sfact \mat{A}_{0,h, \xi}, \mat{B}_{h, \xi} )$ over all $\xi$ is approximately
	equal to the set
	$W_2(\sfact \mat{A}_{0,h, \xi_1}, \mat{B}_{h, \xi_1} )$. Specifically:
	\begin{itemize}
		\item For the mode $\xi = 0$, we have $\mat{Q}_{h} = \mat{0}$. Hence
		$W_2( \sfact \mat{A}_{0, h, 0}, \mat{B}_{h, 0} ) = \{ 1-\sfact^{-1} \}$
		is a single point contained in the set
		$W_2( \sfact \mat{A}_{0, h, \xi_1}, \mat{B}_{h, \xi_1} )$.
		\item The matrices $\mat{Q}_{h, -\xi} = \mat{Q}_{h, \xi}$ and
		$\mat{A}_{0, h, -\xi}= \mat{A}_{0, h, \xi}$ are even functions of $\xi$.  Hence
		the sets $W_2( \sfact \mat{A}_{0, h, \xi}, \mat{B}_{h, \xi} )$ are the same
		for both $\pm \xi$ values.
		\item Figure~\ref{NumericalRange_PPE_maxW2} shows that the value
        $W_{max}(\xi; h)$
		is a decreasing function of $\xi$. Hence, using formula \eqref{NumRangeNS}
		along with the definitions of $W_{max}(\xi;h)$
        (and the fact that $W_{min}(0;h) = 0$),
		one has
		$W_{2}( \sfact \mat{A}_{0, h, \xi}, \mat{B}_{h, \xi} ) \subseteq
		W_2(\sfact \mat{A}_{0, h, \xi_1}, \mat{B}_{h, \xi_1} )$ whenever $\xi > \xi_1$.		
	\end{itemize}
	Hence we have that
	$W_2(\sfact \mat{A}_{0,h}, \mat{B}_{h} ) \approx W_2(\sfact \mat{A}_{0,h, \xi_1}, \mat{B}_{h, \xi_1} )$, where the $\approx$ sign (as opposed to an $=$ sign)
	denotes the fact that there may be an $\mathcal{O}(h)$ error.
\end{remark}

\begin{SCfigure}
	\centering
	\includegraphics[width = 0.30\textwidth]{%
	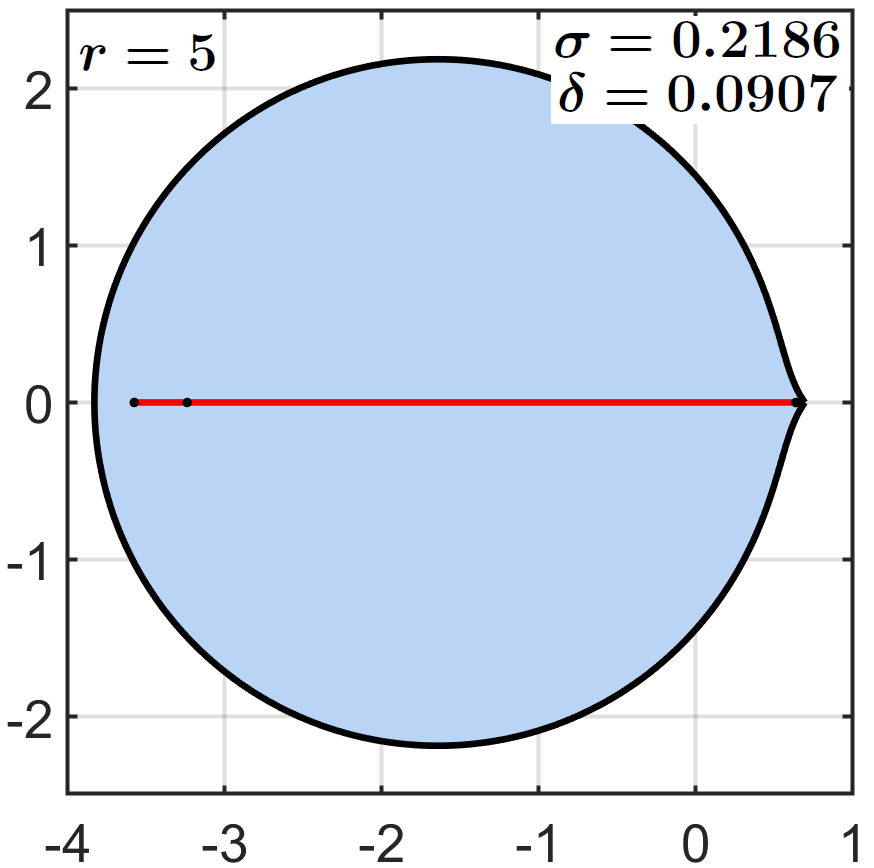}
	\hfill
\caption{Unconditional stability for equation \eqref{SingleComponentChannelFlow}
	and channel length $L_x = 2\pi$.
	The set $W_2(\mat{A}_{h}, \mat{B}_{h})$ (shown in red)
	is contained within the unconditional stability
	region (blue set) when $(\param, \sfact) = (0.0907, 0.2186)$.}
\label{NumericalRange_NS_Channel}
\end{SCfigure}

Based on this important insight, we now choose $(\param^*, \sfact^*)$ for a channel geometry of length $L_x = 2\pi$ and smallest positive wave number $\xi_{1} = 1$.
For grid size $\sizeN_y = 256$, the set $W_2(\sfact \mat{A}_{0, h, \xi_1}, \mat{B}_{h, \xi_1} )$ is obtained by inserting the values $W_{min}(1;h) = 0$ and $W_{max}(1; h) = 0.93$ (see Figure~\ref{NumericalRange_PPE_maxW2}) into equation \eqref{NumRangeNS}.
A crucial observation is that the largest and smallest generalized eigenvalues $\Lambda( \mat{A}_{0, \xi_1; h}, \mat{Q}_{\xi_1; h})$ equal the maximum and minimum values of $W_{max}(\xi_1; h)$, see Figure~\ref{NumericalRange_PPE_maxW2} --- and the gap between these generalized eigenvalues exceeds the unconditional stability capabilities of SBDF3 (see Remark~\ref{Rmk:FailureSBDFDiff}). Hence, the new coefficients (i.e.\ $\param<1$) must be used to achieve unconditional stability. Using a gap parameter of $\eta = 0.1$, order $r = 5$, and values $d_{min} \rightarrow 1-0.93 = 0.07$, $d_{max} \rightarrow 1$ in equation \eqref{InequalitySolution}, leads to $(\param, \sfact) = (0.0907, 0.2186)$. Figure~\ref{NumericalRange_NS_Channel} verifies that this choice in fact satisfies the sufficient conditions for unconditional stability.

We conclude this section with a few important observations.  For a fixed value of
$\xi$, the maximum value $W_{max}(\xi; h)$ remains bounded below $1$ as
$h \rightarrow 0$. This implies that one
value of $(\param, \sfact)$ may be used to stabilize an entire family
of splittings. On the other hand, the sets
$W_2(\sfact\mat{A}_{0, h, \xi}, \mat{B}_{\xi, h})$ do depend on $\xi$ --- and
Figure~\ref{NumericalRange_PPE_maxW2} implies that the sets
$W_2(\sfact\mat{A}_{0, h, \xi}, \mat{B}_{\xi, h})$  (and also the generalized eigenvalues
$\Lambda(\sfact\mat{A}_{0, h, \xi}, \mat{B}_{\xi, h})$) become large when $\xi_1 \rightarrow 0$.
This observation is important because it implies that designing unconditionally stable
schemes requires a choice of $(\param, \sfact)$ that depends on the domain
size $L_x$, and also on the fact that
the domain is a channel geometry.

A natural question then arises: what values can $W_2(\mat{A}_{h}, \mat{B}_{h})$ take for a general fluid dynamics problem? Because the set may depend on the geometry shape (for instance whether $\varOmega$ has corners, see \cite{CozziPego2011}) and size of the computational domain, one would expect that numerical computations may be required to determine or estimate $W_2(\mat{A}_{h}, \mat{B}_{h})$. For instance, one may perform a few rapid computations of $W_2(\mat{A}_{h}, \mat{B}_{h})$ using a coarse mesh (i.e.\ large $h$) and thus small matrices $\mat{A}_{h}$ and $\mat{B}_h$, to obtain a guide for determining the parameters $(\param, \sfact)$ for the fully resolved problem.


\section{Conclusions and outlook}\label{Sec_Conclusions}

With this work on unconditionally stable \imex multistep methods we wish to stress two key messages: first, we advocate to conduct the selection of the \imex splitting and selection of the time-stepping scheme in a simultaneous fashion; and second, it is often possible to achieve unconditional stability in significantly more general settings than one might think at first glance.

The examples and applications discussed herein may serve as a blueprint for how to approach many other types of problems, by using the new stability theory and new \imex schemes to determine feasible and optimal parameters $(\param, \sfact)$ that characterize the scheme and splitting, respectively.

The theoretical foundations of this work establish necessary and sufficient conditions for unconditional stability, resulting from unconditional stability diagrams (that depend only on the scheme) and computable matrix quantities (that depend only on the \imex splitting). This analysis is then used to explain fundamental limitations of the popular SBDF schemes. In particular, it is shown why SBDF can frequently not be extended beyond first or second order --- and how the new schemes can overcome this barrier.

The variable coefficient and nonlinear diffusion examples highlight the practical impact that the new methodology can bring: being able to treat problems whose stiff terms are challenging to invert, without stiff time step restrictions \emph{and} without having to conduct challenging solves. In addition, the theory serves to provide some fundamental insight into unconditional stability (or breakdown thereof) in incompressible fluid flow simulations.

A key limitation of this work is the formal restriction to positive definite matrices $\mat{A}$. This excludes many splittings that would be warranted for intrinsically non-symmetric problems, such as advection or dispersion. Regarding this limitation, it should first be noted that much of the theory persists when the assumptions on $\mat{A}$ are relaxed (for instance, $\mat{A}$ may be a normal matrix with eigenvalues $\lambda$ having complex arguments $|\mathrm{arg}(-\lambda)|$ that are not too large). Second, the extension of the theory to truly non-symmetric $\mat{A}$ is an important subject of future work.


%% file: SecAppendix.tex
\section{Proof of Proposition~\ref{Prop_NumR_VarDiff}}\label{Sec:AppendixProof}

Throughout this section we suppress the subscript $h$ on the matrices $(\mat{A}_h, \mat{B}_h)$, and simply write $(\mat{A}, \mat{B})$. We start with computing the set
\begin{align}\label{Def:ModifiedWp}
	W_1(\mat{A}, \mat{B}) &= \Big\{ \langle \vec v, \mat B \vec v\rangle :
          \langle \vec v, (-\mat A) \vec v\rangle = 1,
          \vec{v} \in \mathbbm{V} \Big\}\/ \\ \label{W1:Comp_line2}
        &= \Big\{ \langle \mat{D}\vec v, \big(\sfact \mat{I} - \diag(\vec{d}) \big)
        \mat{D}\vec v\rangle :
        \langle \mat{D} \vec v,  \mat{D} \vec v\rangle = \sfact^{-1}, \vec{v} \in \mathbbm{V} \Big\}\/
\end{align}
In the expression in \eqref{W1:Comp_line2}, we have used the fact that the derivative
matrix is skew-symmetric
$\mat{D}^{\dag} = \overline{\mat{D}}^T = -\mat{D}$.  Note that $\mat{D}$
is invertible on the space $\mathbbm{V}$ (i.e.~$\mathbbm{V}$ is orthogonal to
$\vec{1}$ --- which is
the nullspace of $\mat{D}$). Making the change of variables
$\vec{y} = \sfact^{\frac{1}{2}} \mat{D} \vec{v}$ in \eqref{W1:Comp_line2}, we observe
that as $\vec{v}$ varies over $\mathbbm{V}$, $\vec{y}$ varies over $\mathbbm{V}$. This yields:
 \begin{align*}
	W_1(\mat{A}, \mat{B}) &= \Big\{ \langle \vec{y}, \big(\mat{I}
        -\sfact^{-1} \diag(\vec{d}) \big)
        \vec{y} \rangle :
        \| \vec{y} \| = 1,  \vec{y} \in \mathbbm{V} \Big\}\/
        = 1 - \sfact^{-1} \/ \sum_{j = 1}^{\sizeN} d(x_j) \/ |y_j|^2,
\end{align*}
where $\vec{y} = \begin{pmatrix}
y_1, \/ y_2, \ldots, y_{\sizeN}
\end{pmatrix}^{T}$.
Since $\| \vec{y} \|^2 = 1$, each value $|y_j|^2$ is real and confined to the region $0 \leq |y_j|^2 \leq 1$ (note that because $\vec{y} \in \mathbbm{V}$, not all vectors $\vec{y}$ are allowed --- only those having zero mean). Combining the results leads to the following inequality:
\[
	d_\mathrm{min} = d_\mathrm{min} \sum_{j = 1}^{\sizeN}
	|y_j|^2 \leq \sum_{j = 1}^{\sizeN} d(x_j) |y_j|^2  \leq d_\mathrm{max} \sum_{j = 1}^{\sizeN}
	|y_j|^2 = d_\mathrm{max}.
\]
Hence, the set $W_1(\mat{A}, \mat{B})$ is real and bounded by:
\[
	1 - \sfact^{-1} d_\mathrm{max} \leq W_1(\mat{A}, \mat{B}) \leq 1 - \sfact^{-1} d_\mathrm{min}.
\]
This concludes the first part of the proof.
To prove the eigenvalue bounds on $\Lambda(\mat{A}, \mat{B})$, the upper bound estimates (i.e.\ the bounds overestimating the largest $\mu$ and underestimating the smallest $\mu$) follow directly from using the established bounds on $W_1(\mat{A}, \mat{B})$ with the fact that $\Lambda(\mat{A}, \mat{B}) \subseteq W_1(\mat{A}, \mat{B})$. Thus, it suffices to prove only the lower bounds. We are interested in bounding the eigenvalues of $\mu (-\mat{A}) \vec{v} = \mat{B}\vec{v}$ with eigenvectors $\vec{1}^T \/ \vec{v}= 0$ restricted to $\mathbbm{V}$. Substituting $\mat{A}$ and $\mat{B}$ into the eigenvalue equation yields:
	\[
		\mu \/\mat{D}^2 \vec{v} = \Big(\mat{D} \mat{S} \mat{D} \Big) \vec{v},
		\quad \textrm{where} \quad
		\mat{S} := (\mat{I} - \sfact^{-1} \diag(\vec{d})).
	\]
As before, let $\vec{y} = \mat{D} \vec{v}$, which is an invertible transformation on $\mathbbm{V}$. Then
	\[
		\mat{D} \Big( \mat{S}\vec{y} - \mu\vec{y} \Big) = 0, \quad
		\vec{y}^T \vec{1} = 0,
	\]
or alternatively
	\begin{align}\label{ConstrainedProblem}
		(\mat{S}- \mu \mat{I}) \vec{y} = \alpha \vec{1},
		\quad \vec{y}^T \vec{1} = 0.
	\end{align}
for some $\alpha \in \mathbb{C}$.
We now solve equation \eqref{ConstrainedProblem} for two separate cases:

\medskip
\noindent\textbf{Case 1:} $\mat{S} \/ \vec{1} = \vec{0}$.
Here $\sfact = d(x_j)$ for all $d(x_j)$. This is only possible if $d(x_j) = d_0$ for
all $j$, i.e.\ one has a constant coefficient diffusion, and thus $\mat{S} = \mat{0}$ identically.
Dotting \eqref{ConstrainedProblem} through by $\vec{1}$ further shows that $\alpha = 0$, thereby forcing all $\mu = 0$. Hence, Proposition~\ref{Prop_NumR_VarDiff} is satisfied (trivially) because $\mu_\mathrm{max} \geq 1 - \sfact^{-1} d_{2,\mathrm{min}} = 0$ and $\mu_\mathrm{min} \leq 1 - \sfact^{-1} d_{2,\mathrm{max}} = 0$.

\medskip
\noindent\textbf{Case 2:} $\mat{S} \/ \vec{1} \neq \vec{0}$.
In this case, we solve equations \eqref{ConstrainedProblem} by first writing the components of $\vec{y}$ in terms of the unknown eigenvalue $\mu$:
	\[
		y_j = \frac{\alpha}{1-\sfact^{-1} \/ d(x_j) - \mu}.
	\]
	Applying the constraint $\vec{1}^T \vec{y} = 0$ to the vector $\vec{y}$,
	shows that the eigenvalues $\mu$ are roots to the following equation:
	\[
		g(\mu) = 0, \quad \textrm{where} \quad
		g(\mu) := \sum_{j = 1}^{\sizeN}
		\big(1 - \sfact^{-1} d(x_j) - \mu\big)^{-1}.
	\]
	Ordering the poles of $g(\mu)$ along the real axis from smallest to
	largest shows that there is at least one root of $g(\mu)$ between the smallest
	two values (or largest two values) of $(1-\sfact^{-1}d(x_j))$. Hence,
	Proposition~\ref{Prop_NumR_VarDiff} follows.	
	

\clearpage

\section{Formulas for the ImEx coefficients}
\label{Sec_ImexCoefficientsTables}
\vspace{-1em}

\begin{table}[!htb]
\begin{small}
\setlength{\mylength}{0.191\textwidth}
\noindent
\begin{tabular}{ |@{~}>{\centering\arraybackslash}p{0.08\textwidth} @{~} |@{~}>{\centering\arraybackslash} p{0.03\textwidth} @{~} || @{~}>{\centering\arraybackslash}p{\mylength} @{~}| @{~}>{\centering\arraybackslash}p{\mylength} @{~}|@{~} >{\centering\arraybackslash}p{\mylength} @{~} |@{~} >{\centering\arraybackslash}p{\mylength} @{~} |}
\hline
	Order &  & $j = 3$ & $j = 2$ & $j = 1$ & $j = 0$  \\ \hline
	1 & $a_j$  & . &  . & $\param$ & $-\param$ \\
	  & $c_j$  & . &  . & 1 & ($\param$-1) \\
	  & $b_j$  & . &  . & 0 & $\param$ \\ \hline
\end{tabular}

\bigskip
\noindent
\begin{tabular}{ |@{~}>{\centering\arraybackslash}p{0.08\textwidth} @{~} |@{~}>{\centering\arraybackslash} p{0.03\textwidth} @{~} || @{~}>{\centering\arraybackslash}p{\mylength} @{~}| @{~}>{\centering\arraybackslash}p{\mylength} @{~}|@{~} >{\centering\arraybackslash}p{\mylength} @{~} |@{~} >{\centering\arraybackslash}p{\mylength} @{~} |}
\hline
	Order &  & $j = 3$ & $j = 2$ & $j = 1$ & $j = 0$  \\ \hline
    2 & $a_j$  & . & $2\param - \frac{1}{2}\param^2$ & $-4\param + 2\param^2$ & $2\param- \frac{3}{2}\param^2$ \\
      & $c_j$  & . & 1 & $2(\param - 1)$ & $(\param - 1)^2$ \\
      & $b_j$  & . & 0 & $2\param$ & $(\param - 1)^2 - 1$ \\ \hline
\end{tabular}

\bigskip
\noindent
\begin{tabular}{ |@{~}>{\centering\arraybackslash}p{0.08\textwidth} @{~} |@{~}>{\centering\arraybackslash} p{0.03\textwidth} @{~} || @{~}>{\centering\arraybackslash}p{\mylength} @{~}| @{~}>{\centering\arraybackslash}p{\mylength} @{~}|@{~} >{\centering\arraybackslash}p{\mylength} @{~} |@{~} >{\centering\arraybackslash}p{\mylength} @{~} |}
\hline
	Order &  & $j = 3$ & $j = 2$ & $j = 1$ & $j = 0$  \\ \hline
    3 & $a_j$  & $3\param - \frac{3}{2}\param^2 + \frac{1}{3}\param^3$ & $-9\param + \frac{15}{2}\param^2 -\frac{3}{2}\param^3$ & $9\param - \frac{21}{2}\param^2 + 3\param^3$ & $-3\param + \frac{9}{2}\param^2 - \frac{11}{6}\param^3$ \\
      & $c_j$  & 1 & $3(\param-1)$ & $3(\param - 1)^2$ & $(\param - 1)^3$ \\
      & $b_j$  & 0 & $3\param$ & $-6\param + 3\param^2$ & $(\param - 1)^3 + 1$ \\ \hline
\end{tabular}

\setlength{\mylength}{0.26\textwidth}
\bigskip
\noindent
\begin{tabular}{ |@{~}>{\centering\arraybackslash}p{0.08\textwidth} @{~} |@{~} >{\centering\arraybackslash}p{0.03\textwidth} @{~} || @{~}>{\centering\arraybackslash} p{\mylength} @{~}| @{~}>{\centering\arraybackslash}p{\mylength} @{~}| @{~}>{\centering\arraybackslash}p{\mylength} @{~}|}
\hline
	Order &        &  		 & $j = 4$   & $j = 3$ \\ \hline
	4     & $a_j$  &  	.	 & $4\param - 3\param^2 + \frac{4}{3}\param^3 - \frac{1}{4}\param^4$ & $-16\param + 18\param^2 - \frac{22}{3}\param^3 + \frac{4}{3}\param^4$ \\
	      & $c_j$  &  	.	 &  1 		 & $4(\param-1)$ \\
	      & $b_j$  &  	.    &  0 		 & $4\param$ \\ \hline \hline

	      &        & $j = 2$ & $j = 1$ & $j = 0$ \\ \hline
	      & $a_j$  & $24\param - 36\param^2 + 18\param^3 -3\param^4$ & $-16\param + 30\param^2 - \frac{58}{3}\param^3 + 4\param^4$ & $4\param - 9\param^2 + \frac{22}{3}\param^3 - \frac{25}{12}\param^4$ \\
      	  & $c_j$  & $6(\param - 1)^2$ & $4(\param - 1)^3$ & $(\param - 1)^4$ \\
      	  & $b_j$  & $-12\param + 6\param^2$ & $12\param - 12\param^2 + 4\param^3$  & $(\param - 1)^4  - 1$ \\ \hline	
\end{tabular}

\setlength{\mylength}{0.399\textwidth}
\bigskip
\noindent
\begin{tabular}{ |@{~}>{\centering\arraybackslash} p{0.08\textwidth} @{~} |@{~} >{\centering\arraybackslash}p{0.03\textwidth} @{~} || @{~} >{\centering\arraybackslash}p{\mylength} @{~}| @{~}>{\centering\arraybackslash}p{\mylength} @{~}|}
\hline
	Order &        & $j = 5$ & $j = 4$   \\ \hline
    5     & $a_j$  & $5\param - 5\param^2 + \frac{10}{3}\param^3 - \frac{5}{4}\param^4 + \frac{1}{5}\param^5$ & $-25\param + 35\param^2 - \frac{65}{3}\param^3 + \frac{95}{12}\param^4 - \frac{5}{4}\param^5$\\
          & $c_j$  & 1 & $5(\param - 1)$\\
          & $b_j$  & 0 & $5\param$ \\ \hline \hline
	 	  &        & $j = 3$ & $j = 2$   \\ \hline
          & $a_j$  &  $50\param - 90\param^2 + \frac{190}{3}\param^3 - \frac{65}{3}\param^4 + \frac{10}{3}\param^5$ & $- 50\param + 110\param^2 - \frac{280}{3}\param^3+ 35\param^4 -5\param^5  $\\
          & $c_j$  & $10(\param - 1)^2$ & $10(\param - 1)^3$ \\
          & $b_j$  & $-20\param + 10\param^2$  & $30\param + 10\param^3 - 30\param^2$  \\ \hline \hline
	      &        & $j = 1$ & $j = 0$   \\ \hline
          & $a_j$  & $25\param - 65\param^2 + \frac{200}{3}\param^3 - \frac{365}{12}\param^4 + 5\param^5$ & $- 5\param + 15\param^2 - \frac{55}{3}\param^3 + \frac{125}{12}\param^4 - \frac{137}{60}\param^5 $ \\
          & $c_j$  & $5(\param - 1)^4$ & $(\param - 1)^5$ \\
          & $b_j$  & $-20\param  + 30\param^2 - 20\param^3+ 5\param^4$ & $(\param - 1)^5 + 1$\\ \hline
    \end{tabular}
	\caption{\imex coefficients for orders 1--5 as functions of $\param$. To use, choose an order and determine a value $0 < \param \leq 1$ small enough to ensure the splitting of choice $(\mat{A}, \mat{B})$ is unconditionally stable. Substitute this value $\param$ into the table to obtain the time stepping coefficients. Coefficients reduce to SBDF when $\param = 1$.}
	\label{Table:ImexCoeff}
\end{small}
\end{table}

\clearpage

%% file: Imex_Practice_Builder_Generic.bbl
\begin{thebibliography}{10}

\bibitem{AbdulleMedovikov2001}
{\sc A.~Abdulle and A.~A. Medovikov}, {\em Second order {C}hebyshev methods
  based on orthogonal polynomials}, Numer. Math., 90 (2001), pp.~1--18.

\bibitem{Akrivis2013}
{\sc G.~Akrivis}, {\em Implicit-explicit multistep methods for nonlinear
  parabolic equations}, Math. Comput., 82 (2012), pp.~45--68.

\bibitem{AkrivisCrouzeixMakridakis1998}
{\sc G.~Akrivis, M.~Crouzeix, and C.~Makridakis}, {\em Implicit-explicit
  multistep finite element methods for nonlinear parabolic problems}, Math.
  Comput., 67 (1998), pp.~457--477.

\bibitem{AkrivisCrouzeixMakridakis1999}
\leavevmode\vrule height 2pt depth -1.6pt width 23pt, {\em Implicit-explicit
  multistep methods for quasilinear parabolic equations}, Numer. Math, 82
  (1999), pp.~521--541.

\bibitem{AkrivisKarakatsani2003}
{\sc G.~Akrivis and F.~Karakatsani}, {\em Modified implicit-explicit {BDF}
  methods for nonlinear parabolic equations}, BIT Numer. Math., 43 (2003),
  pp.~467--483.

\bibitem{AnitescuLaytonPahlevani2004}
{\sc M.~Anitescu, W.~Layton, and F.~Pahlevani}, {\em Implicit for local
  effects, explicit for nonlocal is unconditionally stable}, ETNA, 18 (2004),
  pp.~174--187.

\bibitem{AscherRuuthWetton1995}
{\sc U.~Ascher, S.~J. Ruuth, and B.~Wetton}, {\em Implicit-explicit methods for
  time dependent partial differential equations}, SIAM J. Numer. Anal., 32
  (1995), pp.~797--823.

\bibitem{BadalassiCenicerosBanerjee2003}
{\sc V.~Badalassi, H.~Ceniceros, and S.~Banerjee}, {\em Computation of
  multiphase systems with phase field models}, J. Comput. Phys., 190 (2003),
  pp.~371--397.

\bibitem{BertozziJuLu2011}
{\sc A.~Bertozzi, N.~Ju, and J.-W. Lu}, {\em A biharmonic-modified forward time
  stepping method for fourth order nonlinear diffusion equations}, Discrete and
  Continuous Dynamical Systems, 29 (2011), pp.~1367--1391.

\bibitem{BrunoCubillos2016}
{\sc O.~P. Bruno and M.~Cubillos}, {\em Higher-order in time
  quasi-unconditionally stable {ADI} solvers for the compressible
  {N}avier-{S}tokes equations in {2D} and {3D} curvilinear domains}, J. Comput.
  Phys., 307 (2016), pp.~476--495.

\bibitem{BrunoCubillos2017}
\leavevmode\vrule height 2pt depth -1.6pt width 23pt, {\em On the
  quasi-unconditional stability of {BDF}-{ADI} solvers for the compressible
  {N}avier-{S}tokes equations}, SIAM J. Numer. Anal., 55 (2017), pp.~892--922.

\bibitem{BrunoJimenez2014}
{\sc O.~P. Bruno and E.~Jimenez}, {\em Higher-order linear-time unconditionally
  stable {ADI} methods for nonlinear convection-diffusion {PDE} systems}, J.
  Fluids Eng., 136 (2014), pp.~060904--060904--7.

\bibitem{BrunoLyon2010I}
{\sc O.~P. Bruno and M.~Lyon}, {\em High-order unconditionally stable {FC}-{AD}
  solvers for general smooth domains {I}. basic elements}, J. Comput. Phys.,
  229 (2010), pp.~2009--2033.

\bibitem{BrunoLyon2010II}
\leavevmode\vrule height 2pt depth -1.6pt width 23pt, {\em High-order
  unconditionally stable {FC}-{AD} solvers for general smooth domains {II}.
  {E}lliptic, parabolic and hyperbolic {PDE}s; theoretical considerations}, J.
  Comput. Phys., 229 (2010), pp.~3358--3381.

\bibitem{Ceniceros2002}
{\sc H.~Ceniceros}, {\em A semi-implicit moving mesh method for the focusing
  nonlinear {S}chroedinger equation}, Comm. on Pure and Appl. Anal., 1 (2002),
  pp.~1--14.

\bibitem{ChristliebJonesPromislowWettonWilloughby2014}
{\sc A.~Christlieb, J.~Jones, K.~Promislow, B.~Wetton, and M.~Willoughby}, {\em
  High accuracy solutions to energy gradient flows from material science
  models}, J. Comput. Phys., 257 (2014), pp.~193--215.

\bibitem{CozziPego2011}
{\sc E.~Cozzi and R.~L. Pego}, {\em On optimal estimates for the
  {L}aplace-{L}eray commutator in planar domains with corners}, Proc. Amer.
  Math. Soc., 139 (2011), pp.~1691--1706.

\bibitem{Crouzeix1980}
{\sc M.~Crouzeix}, {\em Une m\'{e}thode multipas implicite-explicite pour
  l'approximation des \'{e}quations d'\'{e}volution paraboliques}, Numer. Math,
  35 (1980), pp.~257--276.

\bibitem{DouglasDupont1970b}
{\sc J.~Douglas and T.~Dupont}, {\em Alternating-direction {G}alerkin methods
  on rectangles}, in Numerical Solution of Partial Differential Equations,
  B.~Hubbard, ed., vol.~II, College Park, Md., 1971, SYNSPADE-1970, Univ. of
  Maryland, Academic Press, New York, pp.~133--213.

\bibitem{DriscollHaleTrefethen2014}
{\sc T.~A. Driscoll, N.~Hale, and L.~N. Trefethen}, eds., {\em Chebfun Guide},
  Pafnuty Publications, Oxford, 2014.

\bibitem{DucheminEggers2014}
{\sc L.~Duchemin and J.~Eggers}, {\em The explicit-implicit-null method:
  removing the numerical instability of {PDE}s}, J. Comput. Phys., 263 (2014),
  pp.~37--52.

\bibitem{ElseyWirth2013}
{\sc M.~Elsey and B.~Wirth}, {\em A simple and efficient scheme for phase field
  crystal simulation}, M2AN, 47 (2013), pp.~1413--1432.

\bibitem{Eyre1998}
{\sc D.~Eyre}, {\em Unconditionally gradient stable time marching the
  {C}ahn-{H}illiard equation}, in Computational and Mathematical Models of
  Microstructural Evolution, J.~W. Bullard, R.~Kalia, M.~Stoneham, and L.~Chen,
  eds., vol.~53, Warrendale, PA, USA, 1998, Materials Research Society,
  pp.~1686--1712.

\bibitem{FrankHundsdorferVerwer1997}
{\sc J.~Frank, W.~Hundsdorfer, and J.~Verwer}, {\em On the stability of {IMEX}
  {LM} methods}, Appl. Numer. Math., 25 (1997), pp.~193--–205.

\bibitem{GlasnerOrizaga2016}
{\sc K.~Glasner and S.~Orizaga}, {\em Improving the accuracy of convexity
  splitting methods for gradient flow equations}, J. Comput. Phys., 315 (2016),
  pp.~52--64.

\bibitem{GuanLowengrubWangWise2014}
{\sc Z.~Guan, J.~Lowengrub, C.~Wang, and S.~Wise}, {\em Second-order convex
  splitting schemes for periodic nonlocal {C}ahn-{H}illiard and {A}llen-{C}ahn
  equations}, J. Comput. Phys., 277 (2014), pp.~48--71.

\bibitem{GuermondMinev2015}
{\sc J.-L. Guermond and P.~Minev}, {\em High-order time stepping for the
  incompressible {N}avier-{S}tokes equations}, SIAM J. Sci. Comput., 36 (2015),
  pp.~A2656--A2681.

\bibitem{HairerNorsettWanner1987}
{\sc E.~Hairer, S.~P. N{\o}rsett, and G.~Wanner}, {\em Solving ordinary
  differential equations {I}: Nonstiff problems}, Springer-Verlag, Berlin,
  second~ed., 1987.

\bibitem{WannerHairer1991}
{\sc E.~Hairer and G.~Wanner}, {\em Solving ordinary differential equations
  {II}: Stiff and differential-algebraic problems}, Springer-Verlag, Berlin,
  second~ed., 1991.

\bibitem{HeisterOlshanskiiRebholz2017}
{\sc T.~Heister, M.~A. Olshanskii, and L.~G. Rebholz}, {\em Decoupled,
  unconditionally stable, higher order discretizations for {MHD} flow
  simulation}, J. Sci. Comput., 71 (2017), pp.~21--43.

\bibitem{HeisterOlshanskiiRebholz2016}
\leavevmode\vrule height 2pt depth -1.6pt width 23pt, {\em Unconditional
  long-time stability method for the 2{D} {N}avier-{S}tokes equations}, Numer.
  Math., 135 (2017), pp.~143--167.

\bibitem{Henshaw1994}
{\sc W.~D. Henshaw}, {\em A fourth-order accurate method for the incompressible
  {N}avier-{S}tokes equations on overlapping grids}, J. Comput. Phys., 113
  (1994), pp.~13--25.

\bibitem{HornJohnson1991}
{\sc A.~Horn and C.~Johnson}, {\em Topics in matrix analysis}, Cambridge
  University Press, 1991.

\bibitem{HundsdorferVerwer2003}
{\sc W.~Hundsdorfer and J.~Verwer}, {\em Numerical solution of time-dependent
  advection-diffusion-feaction equations}, Springer Series in Comput. Math. 33,
  Springer, 2003.

\bibitem{JeltschNevanlinna1981}
{\sc R.~Jeltsch and O.~Nevanlinna}, {\em Stability of explicit time
  discretizations for solving initial value problems}, Numer. Math., 37 (1981),
  pp.~61--91.

\bibitem{JeltschNevanlinna1982}
\leavevmode\vrule height 2pt depth -1.6pt width 23pt, {\em Stability and
  accuracy of time discretizations for initial value problems}, Numer. Math.,
  40 (1982), pp.~245--296.

\bibitem{Johnson1978}
{\sc C.~R. Johnson}, {\em Numerical determination of the field of values of a
  general complex matrix}, SIAM J. Numer. Anal., 15 (1978), pp.~595--602.

\bibitem{JohnstonLiu2002}
{\sc H.~Johnston and J.-G. Liu}, {\em A finite difference method for
  incompressible flow based on local pressure boundary conditions}, J. Comput.
  Phys., 180 (2002), pp.~120--154.

\bibitem{JohnstonLiu2004}
\leavevmode\vrule height 2pt depth -1.6pt width 23pt, {\em Accurate, stable and
  efficient {N}avier-{S}tokes solvers based on explicit treatment of the
  pressure term}, J. Comput. Phys., 199 (2004), pp.~221--259.

\bibitem{JuZhangZhuDu2014}
{\sc L.~Ju, J.~Zhang, L.~Zhu, and Q.~Du}, {\em Fast explicit integration factor
  methods for semilinear parabolic equations}, J. Sci. Comput., 62 (2015),
  pp.~431--455.

\bibitem{KarniadakisIsraeliOrszag1991}
{\sc G.~Karniadakis, M.~Israeli, and S.~A. Orszag}, {\em High-order splitting
  methods for the incompressible {N}avier-{S}tokes equations}, J. Comput.
  Phys., 97 (1991), pp.~414--443.

\bibitem{KassamTrefethen2005}
{\sc A.-K. Kassam and L.~N. Trefethen}, {\em Fourth-order time-stepping for
  stiff {PDE}s}, SIAM J. Sci. Comput., 26 (2005), pp.~1214--1233.

\bibitem{KimMoin1985}
{\sc J.~Kim and P.~Moin}, {\em Application of a fractional step method to
  incompressible {N}avier-{S}tokes equations}, J. Comput. Phys., 59 (1985),
  pp.~308--323.

\bibitem{Koto2009}
{\sc T.~Koto}, {\em Stability of implicit-explicit linear multistep methods for
  ordinary and delay differential equations}, Front. Math. China, 4 (2009),
  pp.~113--129.

\bibitem{LaytonTrenchea2012}
{\sc W.~Layton and C.~Trenchea}, {\em Stability of two {IMEX} methods, {CNLF}
  and {BDF2}-{AB2}, for uncoupling systems of evolution equations}, Appl.
  Numer. Math., 62 (2012), pp.~112–--120.

\bibitem{Leibenzon1945}
{\sc L.~S. Leibenzon}, {\em General problem of the movement of a compressible
  fluid in a porous media}, Izv. Akad. Nauk SSSR, Geography and Geophysics, 9
  (1945), pp.~7--10.
\newblock (Russian).

\bibitem{LeVeque2007}
{\sc R.~J. LeVeque}, {\em Finite difference methods for ordinary and partial
  differential equations: {S}teady-state and time-dependent problems}, SIAM,
  Philadelphia, first~ed., 2007.

\bibitem{LiuLiuPego2007}
{\sc J.-G. Liu, J.~Liu, and R.~L. Pego}, {\em Stability and convergence of
  efficient {N}avier-{S}tokes solvers via a commutator estimate}, Comm. Pure
  Appl. Math., 60 (2007), pp.~1443--1487.

\bibitem{LiuLiuPego2010}
\leavevmode\vrule height 2pt depth -1.6pt width 23pt, {\em Stable and accurate
  pressure approximation for unsteady incompressible viscous flow}, J. Comput.
  Phys., 229 (2010), pp.~3428--3453.

\bibitem{MilewskiTabak1999}
{\sc P.~A. Milewski and E.~G. Tabak}, {\em A pseudo-spectral algorithm for the
  solution of nonlinear wave equations}, SIAM J. Sci. Comput., 21 (1999),
  pp.~1102--1114.

\bibitem{Minion2003}
{\sc M.~L. Minion}, {\em Semi-implicit spectral deferred correction methods for
  ordinary differential equations}, Commun. Math Sci., 1 (2003), pp.~471--500.

\bibitem{Muskat1937}
{\sc M.~Muskat}, {\em The flow of homogeneous fluids through porous media},
  McGrawHill, New York, 1937.

\bibitem{RosalesSeiboldShirokoffZhou2017}
{\sc R.~Rosales, B.~Seibold, D.~Shirokoff, and D.~Zhou}, {\em Unconditional
  stability for multistep {I}m{E}x schemes -- {T}heory}, SIAM J. Numer. Anal.,
  55 (2017), pp.~2336--2360.

\bibitem{ShengWangDuWangLiuChen2010}
{\sc G.~Sheng, T.~Wang, Q.~Du, K.~Wang, Z.~Liu, and L.~Q. Chen}, {\em
  Coarsening kinetics of a two phase mixture with highly disparate diffusion
  mobility}, Commun. Comput. Phys., 8 (2010), pp.~249--264.

\bibitem{ShinLeeLee2017}
{\sc J.~Shin, H.~Lee, and J.-Y. Lee}, {\em Unconditionally stable methods for
  gradient flow using convex splitting {R}unge-{K}utta scheme}, J. Comput.
  Phys., 347 (2017), pp.~367--381.

\bibitem{ShirokoffRosales2010}
{\sc D.~Shirokoff and R.~R. Rosales}, {\em An efficient method for the
  incompressible {N}avier-{S}tokes equations on irregular domains with no-slip
  boundary conditions, high order up to the boundary}, J. Comput. Phys., 230
  (2011), pp.~8619--8646.

\bibitem{Smereka2003}
{\sc P.~Smereka}, {\em Semi-implicit level set methods for curvature and
  surface diffusion motion}, J. Sci. Comput., 19 (2003), pp.~439--456.

\bibitem{TrefethenBau1997}
{\sc L.~N. Trefethen and D.~Bau}, {\em Numerical Linear Algebra}, SIAM,
  Philadelphia, 2000.

\bibitem{Trenchea2016}
{\sc C.~Trenchea}, {\em Second order implicit for local effects and explicit
  for nonlocal effects is unconditionally stable}, Romai J., 12 (2016),
  pp.~163--178.

\bibitem{Varah1980}
{\sc J.~M. Varah}, {\em Stability restrictions on second order, three level
  finite difference schemes for parabolic equations}, SIAM J. Numer. Anal., 17
  (1980), pp.~300--309.

\bibitem{YanChenWangWise2015}
{\sc Y.~Yan, W.~Chen, C.~Wang, and S.~Wise}, {\em A second-order energy stable
  {BDF} numerical scheme for the {C}ahn-{H}illiard equation}, Commun. Comput.
  Phys., 23 (2018), pp.~572--602.

\end{thebibliography}
